\documentclass[12pt]{article}

\usepackage{amssymb}
\usepackage{amsfonts}
\usepackage{amsmath}
\usepackage{amsthm}
\usepackage[all]{xy}

\usepackage{cite}

\newcommand{\Z}{\mathbb{Z}}
\newcommand{\F}{\mathbb{F}}
\newcommand{\Q}{\mathbb{Q}}
\newcommand{\C}{\mathbb{C}}
\newcommand{\Te}{{\rm Teich}}

\theoremstyle{plain}
\newtheorem{thm}{Th{\`e}or{\'e}me}[section]

\newtheorem{assertion}[thm]{Assertion} 

\newtheorem{cor}[thm]{Corollaire}
\newtheorem{lem}[thm]{Lemme}
\newtheorem{proposition}[thm]{Proposition}

\theoremstyle{definition} 
\newtheorem{defi}{D\'{e}finition}[section]

\newtheorem{notation}[defi]{Notation}

\theoremstyle{remark}

\newtheorem{rem}[thm]{Remarque}

\newenvironment{dem}{{\bf D\'emonstration : }}{}
\def\findem{\hfill \mbox{\ $\Box$}}

\begin{document}

\title{Formule de Trace pour les Anneaux de Witt.}
\author{Benali BENZAGHOU, Siham MOKHFI.\\
\ \ \ \ \ \ \ \ \ \ \ \ \ \ \ \ \ \ \ \ \ \\
USTHB, Facult\'{e} des Math\'{e}matiques, \\
B. P. 32, El Alia, 16111, Alger, Alg\'{e}rie\\
USDB BP 270, Route Soumaa Blida, Alg\'{e}rie \\
s\_mokhfi@usdb-blida.dz}
\maketitle

\section{Anneaux de Witt}

\subsection{Le foncteur $\textbf{S}$ des suites }

\begin{notation} \label{anneauS}
Soit $A$ un anneau commutatif. L'anneau $\textbf{S}(A)$ est l'ensemble $A^{\mathbb{N} }$ de toutes les suites d'\'el\'ements de $A$ muni de la structure d'anneau produit. \end{notation} 
On voit donc que la multiplication de $\textbf{S}(A)$ est la multiplication terme \`{a} terme, appel\'ee multiplication de Hadamard par certains auteurs \cite{Barsky},\cite{Benzaghou}. Pour des raisons de clart\'{e} typographique, {\`a} l'instar de \cite{Pulita}, nous repr\'{e}sentons la suite d'\'{e}l\'{e}ments de $A$ de terme g\'{e}n\'{e}ral $a_n$ par la notation $\langle a_n \rangle_{n \in \mathbb{N}}$ quand nous la consid\'{e}rons comme \'{e}l\'{e}ment de l'anneau $\textbf{S}(A)$. 
\begin{notation} \label{projection} Soit $A$ un anneau et $m$ un entier naturel. Nous noterons $\kappa_m$ l'application de $\textbf{S}(A)$ dans $A$ d\'{e}fini par $$ \forall a = \langle a_n \rangle_{n \in \mathbb{N}} \in \textbf{S}(A), \qquad \kappa_m(a) = a_m \; .$$\end{notation} 
Il est clair que ces applications $\kappa_m$ sont des morphismes d'anneaux de $\textbf{S}(A)$ dans $A$.
\begin{notation} \label{morphismeS}
Soit $A$ et $B$ deux anneaux commutatifs et $\rho : A \to B$ un morphisme d'anneaux. Le morphisme $\textbf{S}(\rho)$ de $\textbf{S}(A)$ {\`a}  $\textbf{S}(B)$ est l'application telle que 
$$ \forall \langle a_n \rangle_{n \in \mathbb{N}} \in \textbf{S}(A), \qquad \textbf{S}(\rho) \left( \langle a_n \rangle_{n \in \mathbb{N}}\right) = \langle \rho(a_n) \rangle_{n \in \mathbb{N}} \; .$$
\end{notation} 

On v\'{e}rifie imm\'{e}diatement que les notations \ref{anneauS} et \ref{morphismeS} permettent de d\'{e}finir un foncteur $\textbf{S}$ de la cat\'{e}gorie des anneaux commutatifs dans elle-m\^eme. 

\begin{notation} Pour un id\'{e}al $I$ de $A$, on note $\textbf{S}(I)$ l'ensemble des suites $\langle a_n \rangle_{n \in \mathbb{N}}$ dont tous les termes sont \'{e}l\'{e}ments de $I$.  \end{notation}

Il est imm\'{e}diat de voir que, si $I$ est id\'{e}al de $A$, alors $\textbf{S}(I)$ est un id\'{e}al de $\textbf{S}(A)$. D'autre part, lorsque $I$ et $J$ sont deux id\'{e}aux de $A$, on a l'inclusion $\textbf{S}(I) \textbf{S}(J) \subseteq \textbf{S}(IJ)$. En particulier, pour tout couple $(m, n)$ d'entiers naturels, et pour tout id\'{e}al $I$ de $A$, on a l'inclusion $\textbf{S}(I^m) \textbf{S}(I^n) \subseteq \textbf{S}(I^{m+n})$. Ceci signifie que la suite d'id\'{e}aux $\left( \textbf{S}(I^n) \right)_{n \in \mathbb{N}}$ est une filtration d\'{e}croissante de l'anneau $\textbf{S}(A)$. 

\begin{defi} \label{Stopologie} Soit $A$ un anneau et $I$ un id\'{e}al de $A$. La topologie d\'{e}finie par la filtration $\left( \textbf{S}(I^n) \right)_{n \in \mathbb{N}}$ de $\textbf{S}(A)$ est appel\'{e}e la $I$-topologie de $\textbf{S}(A)$. \end{defi}

\begin{rem} Les morphismes $\kappa_m$ de la notation \ref{projection} sont \'{e}videmment continus quand on munit $\textbf{S}(A)$ de la $I$-topologie et $A$ de la topologie $I$-adique. \end{rem} 

\begin{proposition} \label{SCauchy} 
Soit $A$ un anneau et $I$ un id\'{e}al de $A$. Si l'anneau $A$ est s\'{e}par\'{e} et complet pour la topologie $I$-adique, alors l'anneau $\textbf{S}(A)$ est s\'{e}par\'{e} et complet pour la $I$-topologie. \end{proposition} 
\begin{dem} L'anneau $A$ est s\'{e}par\'{e} pour la topologie $I$-adique si et seulement si l'intersection de tous les id\'{e}aux $I^n$ de $A$ est r\'{e}duite {\`a} $\lbrace 0 \rbrace$. Ceci \'{e}quivaut \'{e}videmment {\`a} dire que l'intersection de tous les id\'{e}aux $\textbf{S}(I^n)$ de l'anneau $\textbf{S}(A)$ est r\'{e}duite {\`a} $\lbrace 0 \rbrace$, c'est-{\`a}-dire que $\textbf{S}(A)$ est s\'{e}par\'{e} pour la $I$-topologie. 

Soit d'autre part une suite $\left( \langle a_{k,n} \rangle_{n \in \mathbb{N}}\right)_{k \in \mathbb{N}}$ (avec $a_{k,n} \in A$ pour tout couple $(k,n)$ d'entiers naturels) d'\'{e}l\'{e}ments de $\textbf{S}(A)$ qui est de Cauchy pour la $I$-topologie. On v\'{e}rifie imm\'{e}diatement que, pour tout entier naturel $n$, la suite $(a_{k,n})_{k \in \mathbb{N}}$ est de Cauchy pour la topologie $I$-adique. Si $A$ est suppos\'{e} complet pour la topologie $I$-adique, il en r\'{e}sulte l'existence d'un \'{e}l\'{e}ment $b_n$ de $A$ qui est limite de $(a_{k,n})_{k \in \mathbb{N}}$ pour la topologie $I$-adique. De plus, pour tout entier $h \in \mathbb{N}$, il existe un entier naturel $k_0$ tel que, pour tout couple $(k,n)$ d'entier naturels tel que $k \ge k_0$, on ait $a_{k,n}-b_n \in I^h$ : ceci r\'{e}sulte du fait que la suite $\left( \langle a_{k,n} \rangle_{n \in \mathbb{N}}\right)_{k \in \mathbb{N}}$ est de Cauchy pour la $I$-topologie de $\textbf{S}(A)$, et que $I^h$ est un ferm\'{e} pour la topologie $I$-adique de $A$.  On voit donc que notre suite converge vers la suite $\langle b_n \rangle_{n \in \mathbb{N}}$ pour la $I$-topologie de $\textbf{S}(A)$. \end{dem} \findem

\begin{proposition} \label{Sevaluation} 
Soit $A$ un anneau, $I$ un id\'{e}al de $A$ et $a = \langle a_n \rangle_{n \in \mathbb{N}}$ un \'{e}l\'{e}ment de $\textbf{S}(I)$. On suppose que l'anneau $A$ est s\'{e}par\'{e} et complet pour la topologie $I$-adique. Il existe un unique morphisme continu ${\rm sv}_{a,I}$ de l'anneau $\textbf{S}(A)\lbrack \! \lbrack t \rbrack \! \rbrack$ des s\'{e}ries enti{\`e}res formelles {\`a} coefficients dans $\textbf{S}(A)$, muni de la topologie $t$-adique, sur l'anneau $\textbf{S}(A)$, muni de la $I$-topologie, tel que : \\
- la restriction de ${\rm sv}_{a,I}$ au sous-anneau $\textbf{S}(A)$ est l'identit\'{e} ; \\
- ${\rm sv}_{a,I}(t) = a$. 
\end{proposition} 
\begin{dem} Les conditions donn\'{e}es caract\'{e}risent au plus un morphisme continu, puisque le sous-anneau $\textbf{S}(A) \lbrack t \rbrack $ de l'anneau $\textbf{S}(A)\lbrack \! \lbrack t \rbrack \! \rbrack$ engendr\'{e} par la partie $\textbf{S}(A) \cup \lbrace t \rbrace$ est dense pour la topologie $t$-adique de $\textbf{S}(A)\lbrack \! \lbrack t \rbrack \! \rbrack$. 

Pour \'{e}tablir l'existence du morphisme ${\rm sv}_{a,I}$, il suffit de construire d'abord sa  restriction au sous-anneau $\textbf{S}(A) \lbrack t \rbrack$, puis de montrer que cette restriction est uniform\'{e}ment continue pour la distance $t$-adique. La premi{\'e}re construction est \'{e}videmment r\'{e}alis\'{e}e en posant 
$$ \forall r \in \mathbb{N}, \forall (a_0, \ldots, a_r) \in \textbf{S}(A)^{r+1}, \qquad {\rm sv}_{a,I}\left(\sum_{j=0}^r a_j t^j \right) = \sum_{j=0}^r a_j a^j \in \textbf{S}(A) \; .$$ 
Pour v\'{e}rifier que cette derni{\`e}re formule d\'{e}finit une application uniform\'{e}ment continue de l'anneau $\textbf{S}(A)\lbrack t \rbrack$, muni de la distance $t$-adique, dans l'anneau $\textbf{S}(A)$, muni de la distance associ\'{e}e {\`a} la filtration de ce dernier anneau par les id\'{e}aux $\textbf{S}(I^n)$, il suffit d'observer que la validit\'{e} de l'implication 
$$ \sum_{j=0}^r a_j t^j \in t^n \textbf{S}(A) \lbrack \! \lbrack t \rbrack \! \rbrack \Rightarrow \sum_{j=0}^r a_j a^j  \in \textbf{S}(I)^n \subseteq \textbf{S}(I^n)\; $$
r\'{e}sulte du fait que $a \in \textbf{S}(I)$. \end{dem} \findem

On justifie de fa{\c con analogue que, si $I$ est un id\'{e}al de $A$ tel que $A$ est s\'{e}par\'{e} et complet pour la topologie $I$-adique, alors, pour tout \'{e}l\'{e}ment $a \in I$, il existe un unique morphisme continu $\varepsilon_a$, dit de sp\'{e}cialisation, de l'anneau $A\lbrack \! \lbrack T \rbrack \! \rbrack$, muni de la topologie $T$-adique, dans l'anneau $A$, muni de la topologie $I$-adique, tel que $\varepsilon_a(b) = b$ pour tout $b \in A$ et $\varepsilon_a(T) = a$.

On d\'{e}finit aussi un unique morphisme $S_T$, continu pour les topologies $T$-adiques, de l'anneau $A \lbrack \! \lbrack T \rbrack \! \rbrack$ des s\'{e}ries enti{\`e}res formelles en une ind\'{e}termin\'{e}e $T$ {\`a} coefficients dans $A$ dans l'anneau $\textbf{S}(A)\lbrack \! \lbrack T \rbrack \! \rbrack$ tel que $S_T(T) = T$ et $S_T(a) = \langle a \rangle_{n \in \mathbb{N}}$ pour tout \'{e}l\'{e}ment $a \in A$. 

La formule suivante peut dans certains cas simplifier le calcul du morphisme ${\rm sv}_{a,I}$. 
\begin{proposition} \label{Sevalformule} Dans les hypoth{\`e}ses et avec les notations de la proposition \ref{Sevaluation}, on a pour tout entier naturel $n$ la relation 
$$ \varepsilon_{\kappa_n(a)} = \kappa_n \circ {\rm sv}_{a,I} \circ S_T \; .$$
\end{proposition} 
\begin{dem} 
Les deux membres de la formule {\`a} montrer sont deux morphismes continus de $A\lbrack \! \lbrack T \rbrack \! \rbrack$, muni de la topologie $T$-adique, dans l'anneau $A$, muni de la topologie $I$-adique, et ils donnent {\`a} tout \'{e}l\'{e}ment $b \in A$, ainsi qu'{\`a} l'ind\'{e}termin\'{e}e $T$, la m\^eme image. Puisque le sous-anneau $A\lbrack T \rbrack$ engendr\'{e} par $A \cup \lbrace T \rbrace$ est dense dans l'anneau $A\lbrack \! \lbrack T \rbrack \! \rbrack$ pour la topologie $T$-adique, on en conclut {\`a} l'\'{e}galit\'{e} d\'{e}sir\'{e}e. \end{dem} \findem

\begin{notation} \label{FrobS} Pour tout anneau commutatif $A$, on note $f_A$ l'application de $\textbf{S}(A)$ dans lui-m\^{e}me d\'{e}finie par $$\forall \langle a_n \rangle_{n \in \mathbb{N}}, \qquad f_A \left( \left\langle a_n \right\rangle_{n\in \mathbb{N}} \right)= \left\langle a_{n+1}\right\rangle_{n\in \mathbb{N}}. $$ \end{notation}
On v\'{e}rifie facilement qu'on a ainsi d\'{e}fini un endomorphisme fonctoriel $f$ du foncteur $\textbf{S}$. 

Alors on peut consid\'erer l'anneau des suites $\textbf{S}(A)$ comme un module sur lui-m\^{e}me de deux fa\c{c}ons diff\'erentes : pour la multiplication de Hadamard, ou bien pour la multiplication $\star$ d\'efinie par 
$$ \forall (a, b) \in \textbf{S}(A) \! \times \! \textbf{S}(A),    \qquad a \star b = f_A(a) b .$$  

\begin{notation} \label{decalageS} On note $v_A$ l'application de $\textbf{S}(A)$ dans lui-m\^{e}me telle que $$
v_A \left( \left\langle a_n\right\rangle_{n\in \mathbb{N}}\right) = \left\langle 0, p.a_{0}, p.a_{1}, ...\right\rangle. $$ \end{notation} On remarque que $v_A$ est un morphisme $\textbf{S}(A)$-lin\'eaire quand on munit l'ensemble de d\'epart de la multiplication $\star$ et l'ensemble d'arriv\'ee de la multiplication de Hadamard. 

\subsection{Polyn\^omes de Witt et application fant\^ome}

Soit $p$ un nombre premier.

\begin{defi}
Pour tout entier $n\geq0$, on appelle $n$-i\'{e}me polyn\^{o}me de Witt l'\'{e}l\'{e}ment ${\rm fant}_n$ de l'anneau $\mathbb{Z}[X_0, ..., X_n]$ des polyn\^{o}mes \`{a} coefficients entiers en $n+1$ ind\'{e}termin\'{e}es d\'{e}fini par  $${\rm fant}_n \left(X_0, ..., X_n\right )=\sum _{i=0}^n p^{i} X_{i}^{p^{n-i}}=X_{0}^{p^{n}}+pX_{1}^{p^{n-1}}+ \ldots + p^{n}X_n . $$ 
\end{defi}

{\bf Formulaire pour les polyn\^omes de Witt}

$${\rm fant}_0 \left(X_0\right )=X_0 ; $$
$${\rm fant}_{n+1} \left(X_0, ..., X_{n+1}\right ) = {\rm fant}_{n} \left( X_{0}^{p}, ..., X_{n}^{p} \right) +p^{n+1}X_{n+1} ; $$
$${\rm fant}_{n+1}\left(X_0, ..., X_{n+1}\right ) =X_o ^{p ^{n+1}}+p \, {\rm fant}_n\left(X_1 , ..., X_{n+1}\right). $$



\begin{defi} \label{fantome} Pour un anneau commutatif $A$, l'application fant\^{o}me de $A$ est l'application de l'ensemble $A^\mathbb{N}$ dans lui-m\^{e}me qui associe \`{a} la suite $a = (a_n)_{n \in \mathbb{N}}$ la suite $${\rm fant}_A(a) = ({\rm fant}_n(a_0, \ldots, a_n))_{n \in \mathbb{N}}. $$ On dit que l'\'{e}l\'{e}ment ${\rm fant}_n(a_0, \ldots, a_n)$ de $A$ est la composante fant\^ome d'indice $n$ de la suite $a = (a_n)_{n \in \mathbb{N}}$. \end{defi}

Le lemme suivant caract\'erise l'image de l'application fant\^ome quand il existe un endomorphisme de l'anneau $A$ qui rel\`{e}ve l'endomorphisme de Frobenius de l'anneau quotient $A/pA$. 

\begin{lem}\label{lemmefantome} 
 {\it Soit $A$ un anneau muni d'un endomorphisme $\sigma$ v\'{e}rifiant pour tout $a\in A$ la congruence $$\sigma\left(a\right)\equiv a^{p} \pmod {pA}  . $$ 
Soit $(u_{n})_{n \in \mathbb{N}}$ un \'{e}l\'{e}ment de $A^\mathbb{N}$. Les conditions suivantes sont \'{e}quivalentes :\\
 a) il existe une suite $(a_n)_{n \in \mathbb{N}}$ d'\'el\'ements de $A$ telle que $u_n={\rm fant}_n\left( a_{0}, ..., a_{n}\right)$ pour tout entier $n \in \mathbb{N}$ ;  \\
 b) la congruence $\sigma\left(u_{n-1}\right)\equiv u_n \pmod {p^{n}A}$ est v\'erifi\'ee pour tout entier $n \in \mathbb{N}$. 

Si de plus $p.1$ n'est pas diviseur de z\'ero dans $A$, alors la suite  $(a_n)_{n \in \mathbb{N}}$ de l'assertion a) est unique. }
\end{lem}
\begin{dem}  Voir \cite[page AC IX.3]{Bourbaki}. \end{dem} \findem
 
 \subsection{Construction de Polyn\^{o}mes}

 
  \begin{proposition}
Soit $k \ge 1$ un entier naturel. Pour tout polyn\^ome $\Phi$ \`{a} coefficients entiers rationnels en k ind\'etermin\'ees, il existe une unique suite $(\varphi_n)_{n \in \mathbb{N}}$ de polyn\^omes \`{a} coefficients entiers en k suites d'ind\'etermin\'ees $(X_{n,j})_{n \in \mathbb{N}, 1 \le j \le k}$ telle que, pour tout entier naturel $n$, on ait  $$ {\rm fant}_n(\varphi_0, \ldots, \varphi_n) = \Phi({\rm fant}_n(X_{0,1}, \ldots, X_{n,1}), \cdots, {\rm fant}_n(X_{0,k}, \ldots, X_{n,k})) .$$ 
De plus, pour tout entier naturel $n$, le polyn\^ome $\varphi_n$ est \'el\'ement de $\mathbb{Z}[X_{m,j}]_{0 \le m \le n, 1 \le j \le k}$. \end{proposition}

\begin{dem} 
Cf \cite[Theorem 5.2 p. 332]{Hazewinkel} ou \cite[Th\'{e}or{\`e}me 1.5]{Christol2}. 
 \end{dem}\findem

 
 
 

Pour les choix $\Phi(U,V) = U+V$, $\Phi(U,V) = UV$, $\Phi(U) = -U$,  la proposition pr\'ec\'edente montre l'existence des suites de polyn\^omes \`{a} coefficients entiers $\textbf{S}=\left(\textbf{S}_{n}\right)_{n\in \mathbb{N}}$, $\textbf{P}=\left(\textbf{P}_{n}\right)_{n\in \mathbb{N}}$, $\textbf{I}=\left(\textbf{I}_{n}\right)_{n\in \mathbb{N}}$ telles que, pour tout entier naturel $n$, on ait :
   $${\rm fant}_{n}\left(\textbf{S}_{0}, ..., \textbf{S}_{n}\right)={\rm fant}_{n}\left(X_{0},..., X_{n}\right)+{\rm fant}_{n}\left(Y_{0},...,Y_{n}\right) ; $$
   $${\rm fant}_{n}\left(\textbf{P}_{0}, ..., \textbf{P}_{n}\right)={\rm fant}_{n}\left(X_{0},..., X_{n}\right){\rm fant}_{n}\left(Y_{0},...,Y_{n}\right) ; $$
    $${\rm fant}_{n}\left(\textbf{I}_{0}, ...,\textbf{ I}_{n}\right)=-{\rm fant}_{n}\left(X_{0},..., X_{n}\right).$$ 
    
  On a par exemple : 
    $$\textbf{S}_{0}=X_{0}+Y_{0}, \qquad \textbf{S}_{1}=X_{1}+Y_{1}-\sum^{p-1}_{i=1}\frac{1}{p}\binom{p}{i} X^{i}_{0}Y^{p-i}_{0}, $$
    $$\textbf{P}_{0}=X_{0}Y_{0}, \qquad \textbf{P}_{1}=pX_{1}Y_{1}+X^{p}_{0}Y_{1}+X_{1}Y^{p}_{0} .$$
    Pour un nombre premier $p>2$ et pour tout entier naturel $n$, on a $\textbf{I}_{n}=-X_{n}$. Si au contraire $p=2$, on a $\textbf{I}_{0}=-X_{0}, \quad \textbf{I}_{1}=-\left(X^{2}_{0}+X_{1}\right), \quad \textbf{I}_2 = -X_0^{4} -X_0^{2}X_1 -X_1^{2}-X_2$. 
    
 \begin{proposition} \label{polynomeFrob}
  Il existe une unique suite $(\textbf{F}_n)_{n \in \mathbb{N}}$ de polyn\^{o}mes \`{a} coefficients entiers en une suite d'ind\'etermin\'ees $(X_n)_{n \in \mathbb{N}}$ telle que, pour tout entier naturel $n$, on ait 
     $${\rm fant}_{n}\left(\textbf{F}_{0}, ..., \textbf{F}_{n}\right)={\rm fant}_{n+1}\left(X_{0},..., X_{n+1}\right). $$
     De plus, pour tout entier naturel $n$, le polyn\^{o}me $\textbf{F}_n$ est \'el\'ement de $\mathbb{Z}[X_0, X_1, \ldots, X_{n+1}]$. 
\end{proposition}     
     \begin{dem} Soit $A_1 = \mathbb{Z}\left[(X_n)_{n \in \mathbb{N}}\right]$. On observe que l'endomorphisme $\theta_1$ de l'anneau $A_1$ tel que $\theta_1(X_n) = X_n^p$ est tel que la congruence $\theta_1(a_1) \equiv a_1^p \pmod{pA_1}$ est satisfaite pour tout \'el\'ement $a_1$ de $A_1$. Pour un entier naturel $n$, on pose $u_n = {\rm fant}_{n+1}\left(X_{0},..., X_{n+1}\right) \in A_1$.  On a 
     $$ \theta_1(u_{n-1}) = {\rm fant}_{n}\left(X_{0}^p,..., X_{n}^p\right)$$ 
     D'apr\`{e}s la formule 2 du formulaire, on a 
     $$ \theta_1(u_{n-1}) = u_n - p^{n+1} X_{n+1} \equiv u_n \pmod{p^nA_1} .$$
     Il suffit d'appliquer le lemme \ref{lemmefantome} pour montrer l'existence et l'unicit\'e d'une suite
     $(\textbf{F}_n)_{n \in \mathbb{N}}$ de polyn\^{o}mes v\'erifiant l'identit\'e 
     $${\rm fant}_{n}\left(\textbf{F}_{0}, ..., \textbf{F}_{n}\right)=u_n $$
     pour tout entier naturel $n$. La derni\`{e}re assertion de la proposition r\'esulte d'une r\'ecurrence sur l'entier $n$, en remarquant que $u_n \in \mathbb{Z}[X_0, \ldots, X_{n+1}]$. \hfill $\Box$\end{dem}
     
    On a par exemple $u_0 = X_{0}^{p}+pX_{1}$ et $u_1 = X_0^{p^2} + p X_1^p + p^2 X_2$, d'o\`{u} $\textbf{F}_{0}=X_{0}^{p}+pX_{1}$ et 
    $\textbf{F}_{1}=X_{1}^{p}+pX_{2}-\sum^{p-1}_{i=0}\binom{p}{i}p^{p-i-1}X_{0}^{pi}X^{p-i}_{1}. $ 
      
\begin{rem} \label{1.8}
Comme ${\rm fant}_{n}\left(\textbf{F}_{0},..., \textbf{F}_{n}\right)\equiv{\rm fant}_{n}\left(X_{0}^{p}, \ldots, X^{p}_{n}\right) \pmod{p^{n+1}A_1}$, 
on peut d\'emontrer  par r\'ecurrence sur l'entier naturel $n$ la congruence
$$\textbf{F}_{n}\equiv X_{n}^{p}\, \pmod{pA_1} . $$
\end{rem}

\subsection{L'anneau $\textbf{W}\!\left(A\right)$ des vecteurs de Witt}

\subsubsection{Op\'{e}rations de Witt}
\begin{defi}
Soit $A$ un anneau commutatif. On d\'{e}finit deux lois de composition internes, appel\'ees addition de Witt et multiplication de Witt, et not\'ees respectivement $+_W$ et $\times_W$, sur l'ensemble $A^{\mathbb{N}}$ des suites d'\'el\'ements de $A$ comme suit:
$$\forall \textbf{a} = \left(a_{n}\right)_{n\in \mathbb{N}}, \forall \textbf{b}= \left(b_{n}\right)_{n\in \mathbb{N}} \qquad \textbf{a} +_W \textbf{b}=\left(S_{n}\left(a_{0}, \ldots, a_n, b_{0}, \ldots, b_n \right) \right)_{n \in \mathbb{N}}. $$
et pareillement 
$$\forall \textbf{a} = \left(a_{n}\right)_{n\in \mathbb{N}}, \forall \textbf{b}= \left(b_{n}\right)_{n\in \mathbb{N}} \qquad \textbf{a} \times_W \textbf{b}=\left(P_{n}\left(a_{0}, \ldots, a_n, b_{0}, \ldots, b_n \right) \right)_{n \in \mathbb{N}}. $$
On note $\textbf{W}(A)$ l'ensemble $A^\mathbb{N}$ des suites d'\'el\'ements de $A$ muni de l'addition de Witt et de la multiplication de Witt. \end{defi}

\begin{thm}
{\it Soit $A$ un anneau commutatif. Le triplet $(\textbf{W}(A), +_W, \times_W)$ est un anneau commutatif. Son \'{e}lement nul est le vecteur $\textbf{0}$ dont toutes les composantes sont nulles. Son \'{e}lement unit\'{e} est le vecteur $\textbf{1}$ dont toutes les composantes sont nulles sauf celle d'indice $0$ qui vaut $1$}. 

L'oppos\'{e} d'un \'{e}lement $\textbf{a}$ de $\textbf{W}\left(A\right)$ est  le vecteur $- \textbf{a} = \left( \textbf{I}_n(a_0, \ldots, a_n) \right)_{n \in \mathbb{N}}$. \end{thm}

Les \'{e}l\'{e}ments de l'anneau $\textbf{W}(A)$ sont appel\'es les vecteurs de Witt sur $A$. 

\begin{proposition}
Pour tout anneau commutatif $A$, l'application ${\rm fant}_{A}$ de $\textbf{W}(A)$ dans $\textbf{S}(A)$ d\'efinie par 
$$ \forall \textbf{a}= \left( a_n \right)_{n \in \mathbb{N}} \in \textbf{W}(A), \qquad {\rm fant}_{A}(\textbf{a}) = \left\langle {\rm fant}_n(a_0, a_1, \ldots, a_n) \right\rangle_{n \in \mathbb{N}}$$ 
est un morphisme d'anneaux. Elle est bijective si $p \cdot 1$ est inversible dans $A$, et injective si le groupe additif de $A$ n'a pas de $p$-torsion. \end{proposition}


\subsubsection{Fonctorialit\'{e}}
 
\begin{defi} 
Soit $\rho \colon B\to A$ un morphisme d'anneaux.
D\'{e}finissons l'application $\textbf{W}\! (\rho)$ de $\textbf{W}(A)$ vers  $\textbf{W}(B)$ telle que $$\forall \:\textbf{a}=\left(a _{n}\right)_{n\in\mathbb{N}}\in \textbf{W}(A), \qquad \textbf{W}\! (\rho)\left(a\right)=\left(\rho\left(a_{n}\right)\right)_{n\in \mathbb{N}}. $$
\end{defi}
\begin{proposition}
{\it L'application $\textbf{W}\left(\rho\right)$ est  un morphisme d'anneaux de $\textbf{W}(A)$ dans $\textbf{W}(B)$. De plus, on a ${\rm fant}_{B} \circ \textbf{W}\! (\rho ) = \textbf{S}(\rho) \circ {\rm fant}_{A}$. } \
\end{proposition}
\begin{proposition} 
L'association de l'anneau $\textbf{W}(A)$ {\`a} l'anneau $A$ et du morphisme $\textbf{W}(\rho)$ au morphisme $\rho : A \to B$, d\'{e}finit un foncteur de la cat\'{e}gorie des anneaux dans elle-m\^{e}me. De plus, l'association de l'application ${\rm fant}_A$ {\`a} l'anneau $A$ d\'{e}finit un morphisme fonctoriel du foncteur $\textbf{W}$ dans le foncteur $\textbf{S}$. 
\end{proposition} 

\subsubsection{La transformation naturelle $\tau$}
\begin{notation} Soit $A$ un anneau. Notons par $\boldsymbol{\tau}_A$ l'application de $A$ dans $\textbf{W}(A)$ qui {\`a} chaque \'{e}l\'{e}ment $a$ de $A$  associe le  vecteur $\boldsymbol{\tau}_A(a) = (a,0,0,\ldots)$.  \end{notation}
\begin{assertion} \label{taunaturelle} La transformation $\tau$ est naturelle, c'est-{\`a}-dire que, pour tout morphisme d'anneaux $\rho : A \to B$, on a $\textbf{W}(\rho) \circ \tau_A = \tau_B \circ \rho$. \end{assertion} 
\begin{proposition} Soient $a,b$ des \'{e}l\'{e}ments de l'anneau $A$ et $\textbf{x}=(x_{n})_{n\in \mathbb{N}}$ un  \'{e}l\'{e}ment de $\textbf{W}(A).$ On a : $${\rm fant}_A(\boldsymbol{\tau}_A (a))= \left\langle a^{p^{n}} \right\rangle_{n\in \mathbb{N}} \; ; $$
$$\boldsymbol{\tau}_A(a)\times\textbf{x}=(a^{p^{n}}x_{n})_{n\in \mathbb{N}} \; ; $$
$$\boldsymbol{\tau}_A(ab)=\boldsymbol{\tau}_A(a)\times \boldsymbol{\tau}_A(b) \; .$$\end{proposition}

\subsubsection{Les homomorphismes de Frobenius et de d\'ecalage}
Soit $ A $ un anneau commutatif. Rappelons que, d'apr\'{e}s la proposition \ref{polynomeFrob}, la suite de polyn\^{o}mes $(\textbf{F}_n)_{n \in \mathbb{N}}$ est d\'{e}finie par la relation ${\rm fant}_n(\textbf{F}_0, \ldots, \textbf{F}_n) = {\rm fant}_{n+1}(X_0, \ldots , X_{n+1})$. 

\begin{defi}
On d\'{e}finit deux applications ${\rm Frob}_A$ et $V_A$ de ${\textbf{W}\left(A\right)}$ dans lui-m\^{e}me en posant : 
$$ \forall \textbf{a}=\left(a_{n}\right)_{n\in \mathbb{N}} \in \textbf{W}(A), \qquad {\rm Frob}_{A}\left(\textbf{a}\right)=\left(\textbf{F}_{n}\left(a_{0}, \ldots ,a_{n+1}\right)\right)_{n\in\mathbb{N}}$$ et 
$$ \forall \textbf{a}=\left(a_{n}\right)_{n\in \mathbb{N}} \in \textbf{W}(A), \qquad
V_{A}\left(\textbf{a}\right)=\left(0,a_{0},a_{1}, \ldots \right). $$ 
L'application ${\rm Frob}_A$ est appel\'ee l'application de Frobenius de $\textbf{W}(A)$, et l'application $V_A$ le d\'{e}calage de $\textbf{W}(A)$. \end{defi}
\begin{assertion} \label{Frobdecanaturel} Les transformations ${\rm Frob}$ et $V$ sont naturelles, c'est-{\`a}-dire que, pour tout homomorphisme d'anneaux $\rho \colon A\to B$, on a les relations : 
$$ \textbf{W}\left(\rho\right)\circ {\rm Frob}_{A}= {\rm Frob}_{B} \circ \textbf{W}\left(\rho\right) \qquad \mbox{ et } \qquad \textbf{W}\left(\rho\right)\circ V_{A}=V_{B} \circ \textbf{W}\left(\rho\right) .$$ 
\end{assertion} 
\begin{proposition} 
L'application ${\rm Frob}_A$ est un endomorphisme de l'anneau $\textbf{W}(A)$ des vecteurs de Witt sur l'anneau $A$. L'application $V_A$ est un endomorphisme du groupe additif de l'anneau $\textbf{W}(A)$. 
\end{proposition}
\begin{proposition} \label{formulesFrobdeca} Pour tout anneau commutatif $A$, on a 
$${\rm fant}_{A} \circ {\rm Frob}_{A}=f_{A} \circ \,{\rm fant}_{A} \qquad \mbox{ et } \qquad {\rm fant}_{A} \circ V_{A}=v_{A} \circ {\rm fant}_{A}. $$ 

Quelque soient les vecteurs de Witt $ \textbf{a}, \textbf{b} \in \textbf{W}(A)$, on a $$V_A \left(\textbf{a}\times {\rm Frob}_{A}\left(\textbf{b}\right)\right)=V_{A}\left(\textbf{a}\right)\times \textbf{b} \qquad \mbox{ et } \qquad V_{A}\left(\textbf{a} \right)\times V_{A}\left(\textbf{b} \right)=pV_{A}\left(\textbf{a} \times \textbf{b} \right)$$
Posons $\textbf{b}_0 = V_{A}\left(\textbf{1}\right)=\left(0,1,0,...,\right)$.  
 Pour tout  vecteur de Witt $\textbf{a} \in\textbf{W}(A)$, on a 
\begin{equation} \label{relations1} V_{A}\left({\rm Frob}_{A}\left( \textbf{a} \right)\right)=\textbf{b}_0 \times \textbf{a} \qquad \mbox{ et } \qquad {\rm Frob}_{A}(V_{A}(\textbf{a})) =p\textbf{a} \; . \end{equation}
\end{proposition} 

\begin{rem} \label{Vnideal}
Pour tout entier $m \in \mathbb{N}$, l'image $V^m \textbf{W}(A)$ de la $m$-i{\`e}me puissance de composition du d\'{e}calage $V_A$ est un id\'{e}al de $\textbf{A}$. De plus, si $\rho : A \to B$ est un morphisme d'anneaux, on a toujours $\textbf{W}(\rho)\left(V^m \textbf{W}(A) \right) \subseteq V^m \textbf{W}(B)$. 
\end{rem}

\begin{proposition}
Soit $A$ un anneau commutatif. Pour tout vecteur de Witt $\textbf{a}\in\textbf{W}\left(A\right)$, on a la congruence
${\rm Frob}_{A}\left(\textbf{a} \right)\equiv  \textbf{a}^{\times p} \pmod{p\textbf{W}\left(A\right)}$, o\`{u} $\textbf{a}^{\times p}$ est le produit dans $\textbf{W}\left(A\right)$ de $p$ \'{e}lements \'{e}gaux \`{a} $\textbf{a}$. \end{proposition} 
\begin{dem} Voir \cite.\end{dem} \findem

\begin{proposition} \label{Frobtau} Soit $A$ un anneau commutatif, et $a$ un \'{e}l\'{e}ment de $A$. On a les \'{e}galit\'{e}s 
$$ {\rm Frob}_A \left( \boldsymbol{\tau}_A(a) \right) = \boldsymbol{\tau}_A(a^{p}) = \left( \boldsymbol{\tau}_A(a) \right)^p \; .$$
\end{proposition} 
\begin{dem} Les relations {\`a} montrer sont v\'{e}rifi\'{e}es dans le cas o{\`u} $p$ n'est pas diviseur de z\'{e}ro dans $A$ parce que les vecteurs de Witt consid\'{e}r\'{e}s ont la m\^{e}me image par l'application fant\^{o}me. Dans le cas g\'{e}n\'{e}ral Voir lemme 3\cite [page AC IX.6]{Bourbaki}, il existe un morphisme surjectif $\rho : B \to A$ tel que $p$ n'est pas diviseur de z\'{e}ro dans $B$. On a donc $a = \rho(b)$ o{\`u} $b \in B$ et par naturalit\'{e} de $\boldsymbol{\tau}$ et de ${\rm Frob}$, on en tire
$$ {\rm Frob}_A \left( \boldsymbol{\tau}_A(a) \right) = {\rm Frob}_A \left( \boldsymbol{\tau}_A(\rho(b)) \right) = \textbf{W}(\rho) \left( {\rm Frob}_B \left( \boldsymbol{\tau}_B(b) \right) \right) \; .$$ Puisque $p$ n'est pas diviseur de z\'{e}ro dans $B$, on sait que ${\rm Frob}_B \left( \boldsymbol{\tau}_B(b) \right) = \boldsymbol{\tau}_B(b^{p})$, et on en d\'{e}duit 
$$ {\rm Frob}_A \left( \boldsymbol{\tau}_A(a) \right) = \textbf{W}(\rho) \left( \boldsymbol{\tau}_B(b^{p}) \right) = \boldsymbol{\tau}_A(\rho(b^{p})) = \boldsymbol{\tau}_A(\rho(b)^{p}) = \boldsymbol{\tau}_A(a^{p}) \; .$$\end{dem} \findem

\subsubsection{D\'{e}veloppement en s\'{e}rie d'un vecteur de Witt}
On munit dor\'{e}navant l'anneau $\textbf{W}(A)$ de la topologie associ\'{e}e {\`a} la filtration d\'{e}croissante $\left(V^m \textbf{W}(A) \right)_{m \in \mathbb{N}}$. On sait que ceci fait de $\textbf{W}(A)$ un anneau topologique s\'{e}par\'{e} et complet \cite[page AC IX.11]{Bourbaki}. Cette topologie va nous permettre d'exprimer tout vecteur de Witt sur $A$ comme somme d'une s\'{e}rie. 

\begin{proposition} \label{developpeWitt} Soit $\textbf{a} = (a_n)_{n \in \mathbb{N}}$ un vecteur de Witt sur l'anneau $A$. La s\'{e}rie de terme g\'{e}n\'{e}ral $V_A^{n}(\boldsymbol{\tau}(a_{n}))$ est convergente de somme $\textbf{a}$. \end{proposition}
\begin{dem} 
Voir \cite[page AC IX.11]{Bourbaki}.
\end{dem} \findem

\subsubsection{La transformation $\Delta$}
\begin{proposition} \label{DefiDelta} 
Il existe un unique morphisme fonctoriel $\Delta$ du foncteur $\textbf{W}$ dans le foncteur $\textbf{W} \circ \textbf{W}$ tel que, pour tout anneau commutatif $A$ et pour tout vecteur de Witt $\textbf{a}$ sur $A$, l'image $\Delta_A(\textbf{a})$ du vecteur $\textbf{a}$ par la composante $\Delta_A$ de $\Delta$ en $A$ v\'{e}rifie
\begin{equation}
{\rm fant}_{\textbf{W}(A)} \left( \Delta_A(\textbf{a}) \right) = \left\langle {\rm Frob}^n(\textbf{a}) \right\rangle_{n \in \mathbb{N}} \; . 
\end{equation}
\end{proposition} 
\begin{dem} Voir \cite[exercice 15, p. AC IX.44]{Bourbaki}. \end{dem} \findem

La proposition \ref{DefiDelta} ne nous sera utile que dans le cas o{\`u} $A = \mathbb{F}_p$ est le corps {\`a} $p$ \'{e}l\'{e}ments. Dans ce cas, on sait \cite{Serre} que $\textbf{W}(A)$ est isomorphe {\`a} l'anneau $\Z_p$ des entiers $p$-adiques et que le morphisme de Frobenius de $\textbf{W}(A)$ s'identifie {\`a} l'identit\'{e} de $\Z_p$. On a donc : 

\begin{proposition} \label{DefiDelta2} Il existe un unique monomorphisme $\Delta$ de $\Z_p$ dans $\textbf{W}(\Z_p)$ tel que 
$$ \forall a \in \Z_p, \qquad {\rm fant}_{\Z_p}(\Delta(a)) = \left\langle a \right\rangle_{n \in \mathbb{N}} \; .$$
\end{proposition} 

\subsection{Vecteurs de Witt et id\'{e}aux} 

\begin{notation} 
Soit $I$ un id\'{e}al d'un anneau commutatif $A$. On note $\textbf{W}(I)$ l'ensemble des vecteurs de Witt sur $A$ dont toutes les composantes sont \'{e}l\'{e}ments de $I$. 
\end{notation} 

\begin{proposition} \label{Wkermorphisme} 
Soit $\rho : A \to B$ un morphisme d'anneaux. L'id\'{e}al noyau du morphisme $\textbf{W}(\rho)$ est stable par les applications ${\rm Frob}_A$ et $V_A$, et est simplement
$$ {\rm ker} \left( \textbf{W}(\rho) \right) = \textbf{W}({\rm ker}(\rho)) \; .$$
\end{proposition} 
\begin{dem} Clair d'apr\'{e}s la naturalit\'{e} des morphismes de Frobenius et de d\'{e}calage, et d'apr\'{e}s la d\'{e}finition de $\textbf{W}(\rho)$. \end{dem} \findem

\begin{cor} \label{Wideal} Pour tout id\'{e}al $I$ de l'anneau commutatif $A$, l'ensemble $\textbf{W}(I)$ est un id\'{e}al de l'anneau $\textbf{W}(A)$ stable par les morphismes de Frobenius et de d\'{e}calage. \end{cor} 

\begin{cor} \label{Wdiff} Soit $I$ un id\'{e}al de l'anneau commutatif $A$, et $\textbf{a} = (a_n)_{n \in \mathbb{N}}$ et $\textbf{b} = (b_n)_{n \in \mathbb{N}}$ deux vecteurs de Witt sur $A$. Alors la congruence 
$$ \textbf{a} \equiv \textbf{b} \pmod{\textbf{W}(I) } $$ 
est \'{e}quivalente {\`a} la suite de congruences 
$$ \forall n \in \mathbb{N}, \qquad a_n \equiv b_n \pmod{I}. $$ \end{cor}
\begin{dem} Soit $\rho : A \to A/I$ la projection naturelle. D'apr{\`e}s la proposition \ref{Wkermorphisme}, on sait que ${\rm ker} \left( \textbf{W}(\rho) \right) = \textbf{W}(I)$. Donc la congruence  $ \textbf{a} \equiv \textbf{b} \pmod{\textbf{W}(I) } $ \'{e}quivaut {\`a} l'\'{e}galit\'{e} $\textbf{W}(\rho) (\textbf{a}) = \textbf{W}(\rho) (\textbf{b})$. Par d\'{e}finition du morphisme $\textbf{W}(\rho)$, ceci \'{e}quivaut {\`a} dire que $\rho(a_n) = \rho(b_n)$ pour tout entier naturel $n$. \end{dem} \findem

\begin{proposition} \label{idealWS} Soit $A$ un anneau et $I$  un id\'{e}al de $A$. Si un vecteur de Witt $\textit{\textbf{a}}$ sur $A$ est \'{e}l\'{e}ment de $\textbf{W}(I)$, alors la suite ${\rm fant}_A(\textit{\textbf{a}})$ est \'{e}l\'{e}ment de $\textbf{S}(I)$. R\'{e}ciproquement, si l'id\'{e}al $I$ satisfait la condition 
\begin{equation} \label{Idivisionparp} \forall a \in A, \qquad p a \in I \Rightarrow  a \in I \; ,
\end{equation}
alors le vecteur de Witt $\textit{\textbf{a}}$ sur $A$ est \'{e}l\'{e}ment de l'id\'{e}al  $\textbf{W}(I)$ d\'{e}s que ${\rm fant}_A(\textit{\textbf{a}}) \in \textbf{S}(I) $. \end{proposition} 
\begin{dem} Si $\textit{\textbf{a}} = \left( a_n \right)_{n \in \mathbb{N}}$ est \'{e}l\'{e}ment de $\textbf{W}(I)$, la condition (\ref{Idivisionparp}) permet de v\'{e}rifier par r\'{e}currence que $a_n \in I$. \end{dem} \findem

\begin{proposition} \label{Widealproduit}
Soit $I$ et $J$ deux id\'{e}aux de l'anneau commutatif $A$. On a toujours l'inclusion 
\begin{equation} \textbf{W}(I) \textbf{W}(J) \subseteq \textbf{W}(I J) . \end{equation}
\end{proposition} 
\begin{dem} D'apr\'{e}s le corollaire \ref{Wideal} appliqu\'{e} au cas de l'anneau $A_2$ de polyn\^omes $\Z \left[ (X_n)_{n \in \mathbb{N}}, (Y_n)_{n \in \mathbb{N}} \right]$ et  l'id\'{e}al $I_2$ des polyn\^omes $P$ qui s'annulent quand on substitue 0 {\`a} toutes les ind\'{e}termin\'{e}es $X_n$ pour $n \in \mathbb{N}$, on sait que le produit du vecteur $(X_n)_{n \in \mathbb{N}}$ par le vecteur $(Y_n)_{n \in \mathbb{N}}$ doit \^{e}tre \'{e}lment de l'id\'{e}al $\textbf{W}(I_2)$. Donc, pour tout entier naturel $n$, le polyn\^ome $\textbf{P}_n(X_0, \cdots, X_n, Y_0, \cdots, Y_n)$ est \'{e}l\'{e}ment de $I_2$, et est donc de la forme $$\textbf{P}_n(X_0, \cdots, X_n, Y_0, \cdots, Y_n) = \sum_{j=0}^n X_j Q_j(X_0, \cdots, X_n, Y_0, \cdots, Y_n) \, , $$ pour certains polyn\^omes $Q_j \in A_2$. Comme la multiplication de Witt est commutative, on a aussi 
$$\textbf{P}_n(X_0, \cdots, X_n, Y_0, \cdots, Y_n) = \sum_{j=0}^n Y_j Q_j(Y_0, \cdots, Y_n, X_0, \cdots, X_n) \, , $$
et donc finalement 
$$\textbf{P}_n(X_0, \cdots, X_n, Y_0, \cdots, Y_n) = \sum_{j=0}^n \sum_{k=0}^n X_j Y_k R_{j,k}(X_0, \cdots, X_n, Y_0, \cdots, Y_n) \, , $$
pour certains polyn\^omes $R_{j,k} \in A_2$. 
Soit maintenant $A$ un anneau quelconque et $I, J$ deux id\'{e}aux de $A$ : on se donne deux vecteurs de Witt sur l'anneau $A$ 
$$ \textbf{a} = (a_n)_{n \in \mathbb{N}} \qquad \mbox{ et } \qquad \textbf{b} = (b_n)_{n \in \mathbb{N}} \; , $$
tel que $\textbf{a}$ est \'{e}l\'{e}ment de $\textbf{W}(I)$ et $\textbf{b}$ est \'{e}l\'{e}ment de $\textbf{W}(J)$. Alors, la composante d'indice $n$ du produit $\textbf{a} \textbf{b}$ est donn\'{e}e par 
$$ \textbf{P}_n(a_0, \ldots, a_n, b_0, \ldots, b_n)  = \sum_{j=0}^n \sum_{k=0}^n a_j b_k R_{j,k}(a_0, \cdots, a_n, b_0, \cdots, b_n) \; , $$
ce qui met en \'{e}vidence le fait que cette composante est \'{e}l\'{e}ment de l'id\'{e}al $IJ$. \end{dem} \findem

\begin{cor} \label{Wtopologiefiltration} Soit $A$ un anneau et $I$ un id\'{e}al de $A$. La suite d'id\'{e}aux $\left( \textbf{W}(I^h) \right)_{h \in \mathbb{N}}$ est une filtration de l'anneau $\textbf{W}(A)$.  \end{cor} 

\begin{defi} \label{definitionItopologie} Soit $A$ un anneau et $I$ un id\'{e}al de $A$. La topologie d\'{e}finie par la filtration $\left( \textbf{W}(I^h) \right)_{h \in \mathbb{N}}$ de $\textbf{W}(A)$ est appel\'{e}e la $I$-topologie de $\textbf{W}(A)$. \end{defi}

\begin{proposition} Soit $A$ un anneau et $I$ un id\'{e}al de $A$. L'application fant\^ome ${\rm fant}_A : \textbf{W}(A) \to \textbf{S}(A)$ est un morphisme continu pour les $I$-topologies sur $\textbf{W}(A)$ et $\textbf{S}(A)$. \end{proposition} 
\begin{dem} En effet, la proposition \ref{idealWS} montre que, pour tout entier naturel $h$, les composantes fant\^omes d'un \'{e}l\'{e}ment de $\textbf{W}(I^h)$ appartiennent {\`a} $I^h$ . \end{dem} \findem

\begin{proposition} \label{Wtopologiecomplet} Soit $A$ un anneau et $I$ un id\'{e}al de $A$. Si l'anneau $A$ est s\'{e}par\'{e} et complet pour la topologie $I$-adique, alors l'anneau $\textbf{W}(A)$ est aussi s\'{e}par\'{e} et complet pour la $I$-topologie. \end{proposition} 
\begin{dem} Si l'anneau $A$ est s\'{e}par\'{e} pour la topologie $I$-adique, on sait que l'intersection de tous les id\'{e}aux $I^h$ lorsque $h$ d\'{e}crit l'ensemble des entiers naturels est r\'{e}duite {\`a} $\lbrace 0 \rbrace$. Par cons\'{e}quent, l'intersection des id\'{e}aux $\textbf{W}(I^h)$ quand $h$ d\'{e}crit $\mathbb{N}$ est r\'{e}duite {\`a} $\lbrace \textbf{0} \rbrace$, ce qui signifie que la $I$-topologie de $\textbf{W}(A)$ est s\'{e}par\'{e}e. Reste {\`a} montrer que toute suite de vecteurs de Witt sur $A$ qui est de Cauchy pour la $I$-topologie converge pour cette topologie. Or soit $\left( \textbf{a}_k \right)_{k \in \mathbb{N}}$ une telle suite. Ecrivons $a_{k,n}$ pour la composante d'indice $n$ du vecteur $\textbf{a}_k \in \textbf{W}(A)$. Par hypoth{\`e}se, pour tout entier naturel $h$, il existe un entier naturel $k_0$ tel que, pour tout couple $(r,s)$ d'entiers naturels tels que $\min(r,s) \ge k_0$, on ait $\textbf{a}_{s} - \textbf{a}_r \in \textbf{W}(I^h)$, c'est-{\`a}-dire que $a_{r,n} - a_{s,n} \in I^h$, de sorte que toutes les suites $(a_{k,n})_{k \in \mathbb{N}}$ obtenues en fixant l'indice $n$ sont de Cauchy pour la topologie $I$-adique de $A$. Par hypoth{\`e}se, elles sont donc convergentes pour la topologie $I$-adique de $A$. Comme $I^h$ est un ferm\'{e} de $A$ pour la topologie $I$-adique, on sait que $b_n - a_{k,n}$ est, quelque soit l'entier naturel $n$,  \'{e}l\'{e}ment de $I^h$ pour tout entier $k \ge k_0$. Soit $\textbf{b}= (b_n)_{n \in \mathbb{N}}$ le vecteur de Witt sur $A$ dont la composante $b_n$ d'indice $n$ est la limite de la suite $(a_{k,n})_{k \in \mathbb{N}}$ pour la topologie $I$-adique. Par le corollaire \ref{Wdiff}, on voit que $\textbf{b} - \textbf{a}_k$ est \'{e}l\'{e}ment de $\textbf{W}(I^h)$ pour tout entier $k \ge k_0$. Ce raisonnement appliqu\'{e} pour tout entier naturel $h$ montre la suite $\left( \textbf{a}_k \right)_{k \in \mathbb{N}}$ converge vers $\textbf{b}$ pour la $I$-topologie. \end{dem} \findem

\begin{proposition} \label{evaluation} Soit $A$ un anneau, $I$ un id\'{e}al de $A$ et $\textbf{x}$ un vecteur de Witt appartenant {\`a} $\textbf{W}(I)$. On suppose que l'anneau $A$ est s\'{e}par\'{e} et complet pour la topologie $I$-adique. Il existe un unique morphisme continu ${\rm ev}_{\textbf{x},I}$ de l'anneau $\textbf{W}(A)\lbrack \! \lbrack t \rbrack\! \rbrack $ des s\'{e}ries enti{\`e}res formelles {\`a} coefficients dans $\textbf{W}(A)$, muni de la topologie $t$-adique, sur l'anneau $\textbf{W}(A)$, muni de la topologie $I$-adique, tel que : \\
\indent - la restriction de ${\rm ev}_{\textbf{x},I}$ au sous-anneau  $\textbf{W}(A) \subset \textbf{W}(A)\lbrack \! \lbrack t \rbrack \! \rbrack$ est l'identit\'{e} ; \\ \indent - ${\rm ev}_{\textbf{x},I}(t) = \textbf{x}$.  \end{proposition}
\begin{dem} Les conditions donn\'{e}es caract\'{e}risent au plus un morphisme continu, puisque le sous-anneau $\textbf{W}(A) \lbrack t \rbrack $ de l'anneau $\textbf{W}(A)\lbrack \! \lbrack t \rbrack \! \rbrack$ engendr\'{e} par la partie $\textbf{W}(A) \cup \lbrace t \rbrace$ est dense pour la topologie $t$-adique de $\textbf{W}(A)\lbrack \! \lbrack t \rbrack \! \rbrack$. 

Pour \'{e}tablir l'existence du morphisme ${\rm ev}_{x,I}$, il suffit de construire d'abord sa  restriction au sous-anneau $\textbf{W}(A) \lbrack t \rbrack$, puis de montrer que cette restriction est uniform\'{e}ment continue pour la distance $t$-adique. La premi{\`e}re construction est \'{e}videmment r\'{e}alis\'{e}e en posant 
$$ \forall r \in \mathbb{N}, \forall (a_0, \ldots, a_r) \in \textbf{W}(A)^{r+1}, \qquad {\rm ev}_{x,I}\left(\sum_{j=0}^r a_j t^j \right) = \sum_{j=0}^r a_j x^j \in \textbf{W}(A) \; .$$ 
Pour v\'{e}rifier que cette derni{\`e}re formule d\'{e}finit une application uniform\'{e}ment continue de l'anneau $\textbf{W}(A)\lbrack t \rbrack$, muni de la distance $t$-adique, dans l'anneau $\textbf{W}(A)$, muni de la distance associ\'{e}e {\`a} la filtration de ce dernier anneau par les id\'{e}aux $\textbf{W}(I^n)$, il suffit d'observer que la validit\'{e} de l'implication 
$$ \sum_{j=0}^r a_j t^j \in t^n \textbf{W}(A) \lbrack \! \lbrack t \rbrack \! \rbrack \Rightarrow \sum_{j=0}^r a_j x^j  \in \textbf{W}(I)^n \subseteq \textbf{W}(I^n)\; $$
r\'{e}sulte du fait que $x \in \textbf{W}(I)$\end{dem} \findem

\begin{proposition} \label{diagrammecommutatif} Soit $A$ un anneau, $I$ un id\'{e}al de $A$ et $\textbf{x}$ un vecteur de Witt appartenant {\`a} $\textbf{W}(I)$. On suppose que l'anneau $A$ est s\'{e}par\'{e} et complet pour la topologie $I$-adique. On note $a= {\rm fant}_A(\textbf{x})$ et ${\rm fant}_T$ le morphisme de $\textbf{W}(A)\lbrack \! \lbrack T \rbrack \! \rbrack$ dans $\textbf{S}(A)\lbrack \! \lbrack T \rbrack \! \rbrack$ d\'{e}fini par 
$$ f(T) = \sum_{n \in \mathbb{N}} \textbf{a}_n T^n \in \textbf{W}(A)\lbrack \! \lbrack T \rbrack \! \rbrack  \mapsto {\rm fant}_T(f(T)) = \sum_{n \in \mathbb{N}} {\rm fant}_A(\textbf{a}_n) T^n \in \textbf{S}(A)\lbrack \! \lbrack T \rbrack \! \rbrack \; .$$
Le diagramme 
\setlength{\unitlength}{1mm}
\begin{center} 
\begin{picture}(50,40)
\thicklines
\put(0.7,0.3){$\textbf{W}(A)$}
\put(14,2){\vector(1,0){30}}
\put(46,0.3){$\textbf{S}(A)$}
\put(25,4){${\rm fant}_A$}
\put(5,35){\vector(0,-1){30}}
\put(50,35){\vector(0,-1){30}}
\put(-4,37){$\textbf{W}(A)\lbrack \! \lbrack T \rbrack \! \rbrack$}
\put(42,37){$\textbf{S}(A)\lbrack \! \lbrack T \rbrack \! \rbrack$}
\put(14, 38){\vector(1,0){26}}
\put(-5,19){${\rm ev}_{\textbf{x},I}$}
\put(53,19){${\rm sv}_{a,I}$} 
\put(25,40){${\rm fant}_T$}
\end{picture} 
\end{center} 
est commutatif.  \end{proposition} 
\begin{dem} Les deux applications ${\rm fant}_A \circ {\rm ev}_{\textbf{x},I}$ et ${\rm sv}_{a,I} \circ {\rm fant}_T$ sont deux morphismes continus de l'anneau $\textbf{W}(A) \lbrack \! \lbrack T \rbrack \! \rbrack$ dans $\textbf{S}(A)$. Comme elles co\"incident en tout \'{e}l\'{e}ment de $\textbf{W}(A)$ et en $T$, et que le sous-anneau $\textbf{W}(A) \lbrack T \rbrack$ de $\textbf{W}(A) \lbrack \! \lbrack T \rbrack \! \rbrack$ engendr\'{e} par $\textbf{W}(A) \cup \lbrace T \rbrace$ est dense pour la topologie $T$-adique dans l'anneau $\textbf{W}(A) \lbrack \! \lbrack T \rbrack \! \rbrack$, elles sont \'{e}gales. \end{dem} \findem

 \subsection{Vecteurs de Witt de longueur finie} 
\begin{defi} Soit $\ell \ge 1$ un entier naturel. Munissons $A^{\ell}$ des op\'{e}rations "somme" et "produit" d\'{e}finies comme suit : si $\textbf{a}=(a_{j})_{0 \le j < \ell}$ et $\textbf{b}=(b_{j})_{0 \le j < \ell}$ sont deux \'{e}l\'{e}ments de $ A^{\ell} $, on pose 
$$\quad\textbf{a}+\textbf{b}= \left(S_{0}(a_{0},b_{0}),\cdots,  S_{\ell-1}(a_{0}, \ldots, a_{\ell-1}, b_{0},\ldots, b_{\ell-1}) \right) $$
et
$$\quad \textbf{a}\times\textbf{b}= \left(P_{0}(a_{0},b_{0}),\cdots, P_{\ell-1}(a_{0}, \ldots, a_{\ell-1}, b_{0},\ldots, b_{\ell-1}) \right).$$

On note $\textbf{W}_\ell(A)$ l'ensemble $A^\ell$ muni des deux op\'{e}rations $+$ et $\times$ ainsi d\'{e}finies. Un \'{e}l\'{e}ment de cet anneau $\textbf{W}_\ell(A)$ est appel\'{e} un vecteur de Witt de longueur $\ell$ sur $A$. \end{defi}

Pour les m{\^e}mes raisons que $\textbf{W}(A)$,  $\textbf{W}_\ell(A)$ est un anneau.
\begin{proposition} \label{Wittlongueurfinie} L'application  de $\textbf{W}(A)$ vers $A^{\ell}$ qui envoie le vecteur de Witt $\textbf{a}=(a_{n})_{n \in \mathbb{N}}$ sur le vecteur de Witt $\left( a_j \right)_{0 \le j < \ell} $ de longueur $\ell$ est un morphisme surjectif d'anneaux, dont le noyau est l'id\'{e}al  $V^{\ell}\textbf{W}(A).$ Par cons\'{e}quent, l'anneau $\textbf{W}_\ell(A)$ des vecteurs de Witt de longueur $\ell$ est isomorphe {\`a} l'anneau quotient $\textbf{W}(A)/V^\ell \textbf{W}(A)$.  \end{proposition}
\begin{dem} \'Evident. \end{dem} \findem
 
 \begin{rem} \label{foncteurlongueurfinie} Si $\rho : A \to B$ est un morphisme d'anneaux, il r\'{e}sulte de la d\'{e}finition du morphisme $\textbf{W}(\rho) : \textbf{W}(A) \to \textbf{W}(B)$ que $\textbf{W}(\rho)(V^\ell \textbf{W}(A)) \subseteq V^\ell \textbf{W}(B)$ et donc ce morphisme induit par passage au quotient un morphisme d'anneaux $\textbf{W}_\ell(\rho)$ de $\textbf{W}_\ell(A)$ dans $\textbf{W}_\ell(B)$. \end{rem} 

 \section{Le morphisme d'Artin-Hasse}
 
 \subsection{S\'erie d'Artin-Hasse associ\'ee \`{a} un vecteur de Witt} 

Nous notons $AH(x)$ la fameuse s\'{e}rie d'Artin-Hasse d\'{e}finie par 
$$ AH(x) = \exp \left( \sum_{n=0}^{\infty} p^{-n} x^{p^{n}} \right) . $$
On sait \cite{Robert},\cite{Koblitz} que les coefficients de la s\'erie enti\`{e}re formelle $AH(x)$ appartiennent \`{a} l'anneau $\mathbb{Z}_p$ des entiers $p$-adiques.  On peut donc consid\'{e}rer cette s\'{e}rie comme {\`a} coefficients dans l'anneau $\Z_p \cap \mathbb{Q}$. 

Pour tout anneau $A$ muni d'un morphisme d'anneaux de $\Z_p \cap \mathbb{Q}$ dans $A$, nous noterons encore $AH(x)$ la s\'{e}rie {\`a} coefficients dans $A$ obtenue en rempla{\c cant chacun des coefficients de la s\'{e}rie d'Artin-Hasse par son image dans $A$. Pour un tel anneau $A$, on note $\mathcal{G}(A) = \lbrace g \in A[[x]] , g(0) = 0 \rbrace$ le mono\"{i}de constitu\'e des s\'eries enti\`{e}res formelles en une ind\'etermin\'ee $x$ \`{a} coefficients dans l'anneau $A$ dont le terme constant est nul, muni de l'op\'eration de composition, et $\Lambda(A)$ le groupe multiplicatif des s\'{e}ries enti\`{e}res formelles en une ind\'etermin\'ee $x$ \`{a} coefficients dans $A$ dont le terme constant vaut 1, muni de la topologie $x$-adique. Nous consid\'erons $\Lambda(A)$ comme muni de l'action $f \mapsto f \circ g$ du mono\"{i}de  $\mathcal{G}(A)$ par substitution \`{a} droite. Pour tout \'{e}l\'{e}ment $\textit{\textbf{a}} = (a_i)_{i \in \mathbb{N}}$ de l'anneau de Witt $\textbf{W}\left(A\right)$, la famille des s\'eries $(AH(a_i x^{p^i}))_{i \in \mathbb{N}}$ est multipliable dans le groupe topologique $\Lambda(A)$. \begin{notation} On note $E(\textit{\textbf{a}})$ le produit de la famille $(AH(a_i x^{p^i}))_{i \in \mathbb{N}}$.\end{notation}
 Un calcul simple montre que 
$$ E(\textit{\textbf{a}}) =\prod_{i=0}^{\infty}AH\left(a_{i}x^{p^{i}}\right) = \exp\left(\sum^{\infty}_{n=0} {\rm fant}_{n}\left(\textit{\textbf{a}}\right)p^{-n}x^{p^{n}}\right) . $$   

\begin{proposition} \label{formulesAH} Soit $A$ un anneau muni d'un morphisme d'anneaux de $\Z_p \cap \mathbb{Q}$ dans $A$. Pour tous vecteurs de Witt $\textbf{a, b}$, et pour tout $\alpha \in A$, on a les formules : 
\begin{equation} \label{AHmorphisme}
E\left(\textbf{a+b}\right) =E\left(\textbf{a}\right) E\left(\textbf{b}\right) ; 
\end{equation}
\begin{equation} \label{AHdeca}
E\left(V_A\left(\textbf{a}\right)\right) =E\left(\textbf{a}\right) \circ x^{p} ;  
\end{equation}
\begin{equation} \label{AHFrob}
E(\textbf{a})^p = E({\rm Frob}_A(\textbf{a})) \circ x^p ; 
\end{equation}
\begin{equation} \label{AHtau} 
E\left(\boldsymbol{\tau}_A(\alpha) \textbf{a} \right) = E(\textbf{a}) \circ (\alpha x) \; .
\end{equation}
\end{proposition}
En particulier, l'application $E$, appel\'ee morphisme d'Artin-Hasse, est un morphisme du groupe additif $\textbf{W}(A)$ dans le groupe multiplicatif $\Lambda(A)$. 
\begin{lem} \label{AHcontinu} Soit un anneau $A$ muni d'un morphisme d'anneaux de $\Z_p \cap \mathbb{Q}$ dans $A$, et $I$ un id\'{e}al de $A$. Si le vecteur de Witt $\textit{\textbf{a}}$ sur $A$ est \'{e}l\'{e}ment de l'id\'{e}al $\textbf{W}(I)$, alors tous les coefficients de la s\'{e}rie enti{\`e}re formelle $E(\textit{\textbf{a}}) - 1$ sont \'{e}l\'{e}ments de $I$. \end{lem}
\begin{dem} On \'{e}crit $\textit{\textbf{a}} = \left( a_i \right)_{i \in \mathbb{N}}$. Les coefficients de la s\'{e}rie enti{\`e}re formelle $AH(x)$ sont \'{e}l\'{e}ments de $A$, et son coefficient constant est $1$. Par cons\'{e}quent, pour tout $i \in \mathbb{N}$, tous les coefficients de la s\'{e}rie $AH\left(a_i x^{p^i}\right) - 1$ sont, puisque $a_i \in I$, \'{e}l\'{e}ments de l'id\'{e}al $I$.  \end{dem}\findem

\subsection{S\'eries de Lubin-Tate} 

\subsubsection{G\'{e}n\'{e}rateur de Tate}
Soit $F(T)$ une s\'erie de Lubin-Tate pour l'uniformisante $p$ de l'anneau des entiers $p$-adiques. On rappelle \cite{Lubin.Tate} que cela signifie que $F(T)$ est une s\'erie enti\`{e}re formelle \`{a} coefficients dans l'anneau $\mathbb{Z}_p$ de la forme 
\begin{equation} \label{LT} F(T) = pT + T^p + pT^2 G(T) \; , \end{equation}
o\`{u} $G(T)$ est \'el\'ement de $\mathbb{Z}_p\lbrack \! \lbrack T \rbrack \! \rbrack$. Une telle s\'erie permet de d\'efinir une suite $(K_m)_{m \in \mathbb{N}}$ d'extensions finies du corps $\mathbb{Q}_p$ des nombres $p$-adiques : $K_m$ est le corps obtenu en adjoignant \`{a} $\mathbb{Q}_p$ les z\'eros de la $(m+1)$-i\`{e}me puissance de composition $F^{\circ (m+1)}(T)$ de la s\'erie $F(T)$. Nous notons $\mathcal{O}_m$ l'anneau des entiers de $K_m$ et $\mathfrak{m}_m$ l'id\'eal maximal de $\mathcal{O}_m$. Un g\'en\'erateur de Tate $( \pi_m )_{m \in \mathbb{N}}$ (relativement \`{a} la s\'erie $F(T)$) est \cite{Christol} une suite d'\'el\'ements de l'id\'eal de valuation d'une extension alg\'ebrique de $\mathbb{Q}_p$ telle que $F(\pi_0) = 0$ et 
\begin{equation} \label{genTate} \forall m \in \mathbb{N}, \qquad F(\pi_{m+1}) = \pi_m \; . \end{equation}  
L'existence d'une infinit\'e de g\'en\'erateurs de Tate relativement \`{a} n'importe quelle s\'erie de Lubin-Tate est facilement \'{e}tablie. On constate que $\pi_m$ est une uniformisante de l'anneau de valuation discr{\`e}te $\mathcal{O}_m$. En outre, si $( \pi_m )_{m \in \mathbb{N}}$ est un g\'en\'erateur de Tate, nous conviendrons de poser $\pi_m = 0$ pour tout entier rationnel $m < 0$, de sorte que la relation (\ref{genTate}) est vraie aussi pour tout $m \in \mathbb{Z}$. 

\subsubsection{Vecteur de Witt de s\'{e}ries formelles d\'{e}fini par une s\'{e}rie de Lubin-Tate} 
\begin{lem} \label{puissancescompositionF}
{\it Si $F(T)$ est une s\'erie de Lubin-Tate pour l'uniformisante $p$ de $\mathbb{Z}_p$, alors la suite $(F^{\circ n}(T))_{n \in \mathbb{N}  }$ des puissances de composition de la s\'erie $F(T)$ est suite des composantes fant\^omes d'un unique vecteur de Witt $\textbf{w}$ ayant toutes ses composantes dans l'id\'eal $T \mathbb{Z}_p\lbrack \! \lbrack T \rbrack \! \rbrack$. }\end{lem}
\begin{dem} On consid\`{e}re sur l'anneau $\mathbb{Z}_p \lbrack \! \lbrack T \rbrack \! \rbrack$, muni de la valuation $T$-adique, l'unique endomorphisme continu $\sigma$ tel que $\sigma(T) = F(T)$. La partie de l'anneau $\mathbb{Z}_p \lbrack \! \lbrack T \rbrack \! \rbrack$ constitu\'ee des s\'eries enti\`{e}res formelles $g$ telles que $\sigma(g) - g^p$ soit \'el\'ement de l'id\'eal $ p \mathbb{Z}_p \lbrack \! \lbrack T \rbrack \! \rbrack$ est \'evidemment un sous-anneau ferm\'e de $\mathbb{Z}_p \lbrack \! \lbrack T \rbrack \! \rbrack$, qui contient l'anneau $\mathbb{Z}_p$ des constantes ainsi que l'ind\'etermin\'ee $T$, elle est donc \'egale \`{a} $\mathbb{Z}_p \lbrack \! \lbrack T \rbrack \! \rbrack$. Par cons\'equent, on peut appliquer le lemme \ref{lemmefantome}, qui montre l'existence d'un unique vecteur de Witt {\it $\textbf{w}$} sur $\mathbb{Z}_p \lbrack \! \lbrack T \rbrack \! \rbrack$ tel que fant$_n(${\it \textbf{w}}$) = F^{\circ n}(T)$ pour tout entier naturel $n$. D'apr\'{e}s la proposition \ref{idealWS},  comme toutes les composantes fant\^omes du vecteur {\it $\textbf{w}$} sont \'el\'ements de l'id\'eal $T \mathbb{Z}_p\lbrack \! \lbrack T \rbrack \! \rbrack$, il en est de m\^{e}me de  toutes ses composantes. \end{dem} \findem

\subsubsection{Sp\'{e}cialisations du vecteur d\'{e}fini par une s\'{e}rie de Lubin-Tate} 
Soit $\left( \pi_m\right)_{m\in \mathbb{N}}$ un  g\'enerateur de Tate relativement \`{a} une s\'erie de Lubin-Tate $F(T)$, et {\it $\textbf{w}$} le vecteur de Witt sur l'anneau $\mathbb{Z}_p\lbrack \! \lbrack T \rbrack \! \rbrack$ dont on a montr\'{e} l'existence au Lemme \ref{puissancescompositionF}. Pour tout entier naturel $m$, on dispose de l'unique homomorphisme $\varepsilon_{\pi_m}$ dit de "sp\'ecialisation en $\pi_m$" de $\mathbb{Z}_p \lbrack \! \lbrack T \rbrack \! \rbrack$ dans $\mathcal{O}_m$ tel que 
\begin{equation} \label{specialisation} \forall g(T) \in \mathbb{Z}_p \lbrack \! \lbrack T \rbrack \! \rbrack, \qquad \varepsilon_{\pi_m}(g(T)) = g(\pi_m). \end{equation} 
\begin{notation} On pose $\boldsymbol{\varpi}_{m} = \textbf{W}(\varepsilon_{\pi_m})(\textit{\textbf{w}})$ pour tout entier naturel $m$. \end{notation} On remarque que le vecteur de Witt $\boldsymbol{\varpi}_{m}$ est \'el\'ement de $W(\mathfrak{m}_m)$ et v\'erifie 
\begin{equation}\label{fantvarpi} 
\forall n \in \mathbb{N}, \qquad {\rm fant}_n(\boldsymbol{\varpi}_{m}) = \pi_{m-n} \; \end{equation}.
\begin{cor} \label{Frobvarpi} Soit $m$ un entier naturel. L'image du vecteur de Witt $\boldsymbol{\varpi}_m$ par le morphisme de Frobenius est \'{e}gale {\`a} $\boldsymbol{\varpi}_{m-1}$ si $m \ge 1$, et {\`a} $\textbf{0}$ si $m=0$. \end{cor}

Soit $F^\Delta$ la s\'{e}rie enti{\`e}re {\`a} coefficients dans $\textbf{W}(\Z_p)$ obtenue en rempla{\c cant chaque coefficient de la s\'{e}rie de Lubin-Tate par son image par le morphisme $\Delta$ de la proposition \ref{DefiDelta2}. Nous allons  calculer l'image de la s\'{e}rie $F^\Delta(T)$ par le morphisme d'\'{e}valuation ${\rm ev}_{\boldsymbol{\varpi}_m, \frak{m}_m} : \textbf{W}(\mathcal{O}_m)\lbrack \! \lbrack T \rbrack \! \rbrack \to \textbf{W}(\mathcal{O}_m)$ d\'{e}termin\'{e} en conformit\'{e} avec la proposition \ref{evaluation}. 

\begin{proposition} \label{evaluationLT} Soit $m$ un entier naturel. On a 
$$ {\rm ev}_{\boldsymbol{\varpi}_m, \frak{m}_m}(F^\Delta(T)) = \begin{cases} \boldsymbol{\varpi}_{m-1} & \mbox{ si } m \ge 1 \\ \textbf{0} & \mbox{ si } m =0 \end{cases}.$$ \end{proposition} 
\begin{dem} 
D'apr{\`e}s la proposition \ref{diagrammecommutatif}, on sait que 
$$ {\rm fant}_{\mathcal{O}_m} \left( {\rm ev}_{\boldsymbol{\varpi}_m, \frak{m}_m}(F^\Delta(T)) \right)  = {\rm sv}_{{\rm fant}_{\mathcal{O}_m}(\boldsymbol{\varpi}_m), \frak{m}_m} \left( {\rm fant}_T(F^\Delta(T)) \right) . $$
Or, par la d\'{e}finition \ref{DefiDelta2} et celle du morphisme $S_T$ de la proposition \ref{Sevalformule}, on a $$ {\rm fant}_T(F^\Delta(T)) = S_T(F(T)) .$$
Par la proposition \ref{Sevalformule}, on a
$$ \kappa_n({\rm sv}_{{\rm fant}_{\mathcal{O}_m}(\boldsymbol{\varpi}_m), \frak{m}_m} \left( {\rm fant}_T(F^\Delta(T)) \right) ) = \varepsilon_{\kappa_n({\rm fant}_{\mathcal{O}_m}(\boldsymbol{\varpi}_m))}(F(T)). $$
Or $\kappa_n({\rm fant}_{\mathcal{O}_m}(\boldsymbol{\varpi}_m)) = \pi_{m-n}$. Donc $$ \kappa_n({\rm sv}_{{\rm fant}_{\mathcal{O}_m}(\boldsymbol{\varpi}_m), \frak{m}_m} \left( {\rm fant}_T(F^\Delta(T)) \right) ) = F(\pi_{m-n}) = \pi_{m-1-n} \; .$$
On en tire pour $m \ge 1$ 
$$ {\rm sv}_{{\rm fant}_{\mathcal{O}_m}(\boldsymbol{\varpi}_m), \frak{m}_m} \left( {\rm fant}_T(F^\Delta(T)) \right) = \langle \pi_{m-1-n} \rangle_{n \in \mathbb{N}} = {\rm fant}_{\mathcal{O}_m}(\boldsymbol{\varpi}_{m-1}), $$
et pour $m =0$ 
$$ {\rm sv}_{{\rm fant}_{\mathcal{O}_0}(\boldsymbol{\varpi}_0), \frak{m}_0} \left( {\rm fant}_T(F^\Delta(T)) \right) = \langle \pi_{-1-n} \rangle_{n \in \mathbb{N}} = {\rm fant}_{\mathcal{O}_0}(\textbf{0}).$$

Comme $p$ n'est pas diviseur de z\'{e}ro dans $\mathcal{O}_m$, l'application ${\rm fant}_{\mathcal{O}_m}$ est injective, ce qui permet de conclure {\`a} l'\'{e}galit\'{e} d\'{e}sir\'{e}e. 
\end{dem} \findem

\subsection{Exponentielle de Robba}
\begin{defi}
Soit $\pi = (\pi_m)_{m \in \mathbb{N}}$ un g\'en\'erateur de Tate relativement \`{a} une s\'erie de Lubin-Tate $F(T)$ pour l'uniformi\-sante $p$ de $\mathbb{Z}_p$, et fixons un entier $m \in \mathbb{N}$. On associe \`{a} ces donn\'ees l'exponentielle de Robba \cite{Pulita} qui est par d\'efinition la s\'erie enti\`{e}re formelle \`{a} coefficients dans l'anneau $\mathcal{O}_m$ not\'ee  $e_{m,\pi}\left(x\right)$\ 
telle que $$e_{m,\pi}=E\left(\boldsymbol{\varpi}_{m}\right)=\exp\left(\pi_{m}x+\pi_{m-1}p^{-1}x^{p}+...+\pi_{0}p^{-m}x^{p^{m}}\right).$$
\end{defi}

\begin{proposition} \label{convRobba}
{\it Le rayon de convergence de l'exponentielle de Robba $e_{m,\pi}$ est au moins \'{e}gal {\`a} $1$}. \end{proposition}
\begin{dem} Comme le vecteur de Witt $\boldsymbol{\varpi}_m$ est \'{e}\'{e}l\'{e}ment de l'anneau $\textbf{W}(\mathcal{O}_m)$, son image $e_{m, \pi}(x)$ par le morphisme de Artin-Hasse est \'{e}l\'{e}ment du groupe $\Lambda(\mathcal{O}_m)$. Or toute s\'{e}rie {\`a} coefficients entiers converge sur le disque ouvert de $\mathbb{C}_p$ de rayon 1. \end{dem} \findem

\begin{lem} \label{lemmeChristol} 
Soit $A$ un anneau contenant $\mathcal{O}_m = \mathbb{Z}_p[\pi_{m}]$ et $\textbf{a} =\left(a_{n}\right)_{n \in \mathbb{N}}$ un vecteur de Witt de $\textbf{W}\left(A\right)$. On a l'identit\'{e}
\begin{equation} E\left(\boldsymbol{\varpi}_{m}\textbf{a}\right)=\prod_{i=0}^{m}e_{m-i,\pi} \left(a_{i}x^{p^{i}}\right). \end{equation}
\end{lem}
\begin{dem}
On sait par la proposition \ref{developpeWitt} que 
$$ \textit{\textbf{a}} = \sum_{i=0}^\infty V^i(\boldsymbol{\tau}(a_i)) \; .$$
On en d\'{e}duit que 
$$ \boldsymbol{\varpi}_m  \textit{\textbf{a}} = \sum_{i=0}^\infty \boldsymbol{\varpi}_m V^i(\boldsymbol{\tau}(a_i)) \; .$$
D'apr\'{e}s la proposition \ref{formulesFrobdeca}, on obtient 
$$ \boldsymbol{\varpi}_m  \textit{\textbf{a}} = \sum_{i=0}^\infty  V^i({\rm Frob}^i(\boldsymbol{\varpi}_m) \boldsymbol{\tau}(a_i)) \; ,$$
ce qui fournit, d'apr\'{e}s le lemme \ref{Frobvarpi}, l'identit\'{e} 
$$ \boldsymbol{\varpi}_m  \textit{\textbf{a}} = \sum_{i=0}^m  V^i(\boldsymbol{\varpi}_{m-i}  \boldsymbol{\tau}(a_i)) \; .$$
D'apr\'{e}s la proposition \ref{formulesAH}, cette derni{\`e}re identit\'{e} entre vecteurs de Witt implique que 
$$ E( \boldsymbol{\varpi}_m  \textit{\textbf{a}} ) = \prod_{i=0}^m E (V^i(\boldsymbol{\varpi}_{m-i})  \boldsymbol{\tau}(a_i)) = \prod_{i=0}^m E(\boldsymbol{\varpi}_{m-i}) \circ (a_i x^{p^i}) , $$
ce qui est la formule {\`a} d\'{e}montrer. \end{dem}\findem\\

La proposition \ref{convRobba} et le lemme \ref{lemmeChristol} entra\^inent le corollaire suivant. 
\begin{cor} \label{surconvergence} Soit $\textbf{a} = (a_n)_{n \in \mathbb{N}}$ un vecteur de Witt de $\textbf{W}(\mathbb{C}_p)$. 
Si on a $|a_{i}|<1$ pour tout indice $i\in \lbrack0, m \rbrack$, le rayon de convergence de la fonction $E\left(\boldsymbol{\varpi}_{m} \textbf{a} \right)$ est strictement plus grand que $1$.\end{cor}

\subsection{Exponentielle de Pulita}

\subsubsection{Frobenius absolu} 
Soit $\mathcal{O} = \cup_{m \in \mathbb{N}} \mathcal{O}_m$ la sous-$\Z_p$-alg{\`e}bre de $\mathbb{C}_p$ engendr\'{e}e par les racines des puissances de composition de la s\'{e}rie de Lubin-Tate $F(T)$. C'est un anneau local d'id\'{e}al maximal $\frak{m}_{\mathcal{O}} = 
\cup_{m \in \mathbb{N}} \frak{m}_m$ et de corps r\'{e}siduel $\mathbb{F}_p$. On note $K$ le corps des fractions de l'anneau $\mathcal{O}$, de sorte que $K = \cup_{m \in \mathbb{N}}K_m$. 

\begin{proposition} \label{Frobeniusabsolu}
Il existe un automorphisme continu $\varphi$ du corps valu\'{e} $\C_p$ tel qu'on ait $\varphi(x) = x$ pour tout $x$ de $K$ et $\vert \varphi(x) - x^p \vert < 1$ pour tout $x \in \C_p$ tel que $\vert x \vert \le 1$. 
\end{proposition} 
 
Nous appelons Frobenius absolu tout automorphisme de $\C_p$ satisfaisant les conditions de la proposition \ref{Frobeniusabsolu}. Un tel automorphisme n'est pas unique, mais nous en fixerons un dans la suite, constamment not\'{e} $\varphi$. Pour un vecteur de Witt $\textit{\textbf{a}}$ sur $\C_p$, nous notons $\textit{\textbf{a}}^\varphi$ le vecteur image de $\textit{\textbf{a}}$ par l'automorphisme $\textbf{W}(\varphi)$. 

\begin{lem} \label{imageFrobracine} Si $t \in \C_p$ est une racine de l'unit\'{e} d'ordre premier {\`a} $p$, alors $\varphi(t) = t^p$. \end{lem}
\begin{dem} Les \'{e}l\'{e}ments $\varphi(t)$ et $t^p$ sont deux racines de l'unit\'{e} d'ordre $r$ premier {\`a} $p$, alors que, par d\'{e}finition de $\varphi$, on sait que leurs classes dans le corps r\'{e}siduel $\widetilde{\C_p}$ co\"incident. Comme le polyn\^{o}me $T^r - 1$ n'a que des racines simples dans le corps $\widetilde{\C_p}$ de caract\'{e}ristique $p$, il en r\'{e}sulte que $\varphi(t) = t^p$.  \end{dem} \findem

En comparant la valeur absolue $\vert \cdot \vert$ de $\C_p$ {\`a} la valeur absolue $x \mapsto \vert \varphi(x) \vert$,  on voit facilement\cite[p. 4]{Artin} que le Frobenius absolu est n\'{e}cessairement une isom\'{e}trie de $\C_p$. 

\begin{lem} \label{stabiliteidealFrobabsolu} Si $I$ est un id\'{e}al de l'anneau $\C_p^\circ$, alors $\varphi(I) \subseteq I$. \end{lem}
\begin{dem} Soit $x$ un \'{e}l\'{e}ment de $I$. Puisque $\varphi$ est une isom\'{e}trie, le quotient $\frac{\varphi(x)}{x}$ a une valeur absolue \'{e}gale {\`a} 1, donc est \'{e}l\'{e}ment de l'anneau des entiers $\C_p^\circ$. Par cons\'{e}quent $\varphi(x) \in x \C_p^\circ \subseteq I$. \end{dem} \findem

La propri\'{e}t\'{e} \'{e}vidente suivante nous sera utile ult\'{e}rieurement. 
\begin{assertion} \label{Frobasolucommute} Le Frobenius absolu commute avec le d\'{e}calage et le Frobenius de $\textbf{W}(\C_p)$. \end{assertion}

\subsubsection{D\'{e}finition de la s\'{e}rie exponentielle de Pulita}

\begin{defi} \label{expoPulita} 
Soit $m$ un entier naturel et $\textit{\textbf{a}}$ un vecteur de Witt sur $\C_p$. L'exponentielle de Pulita attach\'{e}e {\`a} ces donn\'{e}es est la s\'{e}rie 
{\it $$\theta_m(\textbf{a}) =\frac{E\left(\boldsymbol{\varpi}_{m}\textbf{a}\right)\left(x\right)}{E\left(\boldsymbol{\varpi}_{m}\textbf{a}^\varphi \right)\left(x^{p}\right)}. 
$$}
\end{defi}
D'apr{\`e}s la proposition \ref{formulesAH}, l'exponentielle de Pulita s'exprime sous la forme 
{\it \begin{equation} \label{expPulita} 
\theta_m(\textbf{a}) = E\left(\boldsymbol{\varpi}_{m}\textbf{a} - V\left(\boldsymbol{\varpi}_{m}\textbf{a}^\varphi \right) \right) . 
\end{equation}}

En particulier, l'application $\theta_m$ est un morphisme du groupe additif $\textbf{W}(\C_p)$ dans le groupe multiplicatif $\Lambda(\C_p)$. 

\begin{lem} \label{Itheta} Soit $I$ un id\'{e}al de l'anneau $\C_p^\circ$, et $\textit{\textbf{a}}$ un \'{e}l\'{e}ment de $\textbf{W}(I)$. Alors tous les coefficients de la s\'{e}rie $\theta_m(\textit{\textbf{a}}) - 1$ sont \'{e}l\'{e}ments de l'id\'{e}al $\pi_m I$. \end{lem} 
\begin{dem} D'apr{\`e}s le lemme \ref{AHcontinu} et l'\'{e}galit\'{e} (\ref{expPulita}), il suffit de montrer que le vecteur $\boldsymbol{\varpi}_{m}\textit{\textbf{a}} - V\left(\boldsymbol{\varpi}_{m} \textit{\textbf{a}}^\varphi \right)$ est \'{e}l\'{e}ment de l'id\'{e}al $\textbf{W}(\pi_m I)$. D'apr{\`e}s le lemme \ref{puissancescompositionF} et l'\'{e}quation (\ref{specialisation}) qui d\'{e}finit le morphisme d'\'{e}valuation en $\pi_m$, il est clair que le vecteur $\boldsymbol{\varpi}_m$ est \'{e}l\'{e}ment du groupe $\textbf{W}(\pi_m \Z_p\lbrack \pi_m \rbrack)$, et donc de l'id\'{e}al $\textbf{W}(\pi_m \C_p^\circ)$ de l'anneau des vecteurs de Witt sur $\C_p^\circ$. D'autre part, on sait par le lemme \ref{stabiliteidealFrobabsolu} que le vecteur $\textit{\textbf{a}}^\varphi$ est aussi \'{e}l\'{e}ment de $\textbf{W}(I)$. En utilisant la proposition \ref{Widealproduit}, on en d\'{e}duit que les vecteurs $\boldsymbol{\varpi}_m \textit{\textbf{a}}$ et $\boldsymbol{\varpi}_m \textit{\textbf{a}}^\varphi$ sont \'{e}l\'{e}ments de l'id\'{e}al $\textbf{W}(\pi_m I)$. D'o{\`u} le r\'{e}sultat d\'{e}sir\'{e}. \end{dem} \findem

Nous d\'{e}sirons maintenant pr\'{e}ciser les morphismes $\theta_m \circ V_{\C_p}^k$, o\`u $k \in \mathbb{N}$. 

 \begin{lem} \label{pulitadecal}Soit un vecteur de Witt $\textit{\textbf{a}}$ sur $\C_p^\circ,$ et deux entiers naturels $m$ et $k$, alors $$\theta_m(V^k(\textit{\textbf{a}}))= \begin{cases} 1 & \mbox{ si } m < k \\ 
\theta_{m-k}(\textit{\textbf{a}})\circ x^{p^k} & \mbox{ si } m \geq k \end{cases}. $$\end{lem}
\begin{dem}   Dans le cas o{\`u} $m < k$, on a ${\rm fant}_{\C_p}(\boldsymbol{\varpi}_m) {\rm fant}_{\C_p}(V^k (\textit{\textbf{a}})) = 0$. Comme l'application fant\^{o}me ${\rm fant}_{\C_p}$ est injective, ceci montre que $\boldsymbol{\varpi}_m V^k (\textit{\textbf{a}})) = \textbf{0}$. De m{\^e}me, en utilisant l'assertion \ref{Frobasolucommute}, on voit que $\boldsymbol{\varpi}_m V^k (\textit{\textbf{a}}))^\varphi = \textbf{0}$.  Il r\'{e}sulte alors de l'\'{e}quation (\ref{expPulita}) que $\theta_m(V^k(\textit{\textbf{a}})) = 1$. 

Supposons maintenant que $m \ge k$.  La formule {\`a} montrer \'{e}tant \'{e}videmment vraie pour $k=0$, sa validit\'{e} pour tout $k \ge m$ suit par r\'{e}currence de sa validit\'{e} pour $k = 1 \ge m$ ; on suppose donc que $k = 1 \ge m$. Alors, par la formule (\ref{expPulita}), on a 
$$\theta _m(V(\textit{\textbf{a}}))=E\left(\varpi_mV(\textit{\textbf{a}})-V(\varpi_mV(\textit{\textbf{a}}^{\varphi})\right).$$
Or, utilisant la proposition \ref{formulesFrobdeca}, on obtient 
$$V(\textit{\textbf{a}}))=
 V(\rm Frob(\varpi_m)\textit{\textbf{a}})- V^2(\rm Frob(\varpi_m)\textit{\textbf{a}}^{\varphi}) = V \left(\rm Frob(\varpi_m)(\textit{\textbf{a}})-V(\rm Frob(\varpi_m)\textit{\textbf{a}}^{\varphi} \right),$$ c'est {\`a} dire, d'apr{\`e}s le lemme \ref{Frobvarpi}, $$\theta _m(V(\textit{\textbf{a}}))=
 E\left(V(\varpi_{m-1}\textit{\textbf{a}}-V(\varpi_{m-1}\textit{\textbf{a}}^{\varphi}))\right),$$ ce qui, d'apr{\`e}s la formule (\ref{AHdeca}) et la d\'{e}finition \ref{expoPulita} de $\theta_{m-1}$, conduit {\`a} l'\'{e}galit\'{e} d\'{e}sir\'{e}e $$\theta _m(V(\textit{\textbf{a}}))= 
 \theta_{m-1}(\textit{\textbf{a}})\circ x^{p}.$$ \end{dem}\findem
 
 \begin{lem} \label{thetamtaut} Soit $t$ une racine de l'unit\'{e} d'ordre premier {\`a} $p$. Alors on a $$ \theta_m(\boldsymbol{\tau}(t)) = \theta_m(\textbf{1}) \circ (tx) \; . $$ \end{lem} 
 \begin{dem} Par le lemme \ref{imageFrobracine}, on sait que $\varphi(t) = t^p$. Utilisant alors le lemme \ref{Frobtau}, on en d\'{e}duit que $\boldsymbol{\tau}(t)^\varphi = \boldsymbol{\tau}(\varphi(t)) = \boldsymbol{\tau}(t^p) = {\rm Frob} \left( \tau(t) \right)$. 
 D'o{\`u}, par la proposition \ref{formulesFrobdeca}, les \'{e}galit\'{e}s
 $$ \boldsymbol {\varpi}_m\boldsymbol{\tau}(t)-V(\boldsymbol {\varpi}_m\boldsymbol{\tau}(t)^{\varphi} ) = \boldsymbol {\varpi}_m\boldsymbol{\tau}(t) - \boldsymbol{\tau}(t) V( \boldsymbol{\varpi}_m ) = \boldsymbol{\tau}(t) \left( \boldsymbol{\varpi}_m - V (\boldsymbol{\varpi}_m ) \right)\; .$$
 Par la formule (\ref{AHtau}) de la proposition \ref{formulesAH}, on en d\'{e}duit 
 $$\theta_m(\boldsymbol{\tau}(t)) = E\left(\boldsymbol{\tau}(t) \left( \boldsymbol{\varpi}_m - V (\boldsymbol{\varpi}_m ) \right)\right) = E( \boldsymbol{\varpi}_m - V (\boldsymbol{\varpi}_m ) ) \circ (tx) , $$
 ce qui est l'\'{e}galit\'{e} voulue.

 \end{dem} \findem

\subsubsection{Rayon de convergence des s\'{e}ries de Pulita} 
\begin{lem} \label{bm} Pour tout entier naturel $m$, il existe un vecteur de Witt $\textbf{b}_m$ \'{e}l\'{e}ment de l'id\'{e}al $\textbf{W}(\frak{m}_{m+1})$ de l'anneau $\textbf{W}(\mathcal{O}_{m+1})$ tel que 
\begin{equation}
\boldsymbol{\varpi}_{m} = \boldsymbol{\varpi}_{m+1} (\textbf{b}_m + p \cdot \textbf{1}).
\end{equation}
\end{lem}
\begin{dem} 
Par d\'{e}finition, il existe une s\'{e}rie $G(T) \in \Z_p\lbrack \! \lbrack T \rbrack \! \rbrack$ telle que la s\'{e}rie de Lubin-Tate $F(T)$ est de la forme donn\'{e}e par l'\'{e}quation (\ref{LT}), et donc 
$$ F(T) = T \left( p + T^{p-1} + pT G(T) \right).$$ 
Soit $G^\Delta(T)$ la s\'{e}rie {\`a} coefficients dans l'anneau $\textbf{W}(\Z_p) \subset \textbf{W}(\mathcal{O}_{m+1})$ obtenue en rempla{\c c}ant les coefficients de la s\'{e}rie $G(T)$ par leurs images par le morphisme $\Delta$ de la d\'{e}finition \ref{DefiDelta2}. On a 
$$ F^\Delta(T) = T \left( p \cdot \textbf{1}+ T^{p-1} + (p \cdot \textbf{1} ) T G^\Delta(T) \right).$$ 
Par la proposition \ref{evaluationLT}, on a donc 
$$ \boldsymbol{\varpi}_m = \boldsymbol{\varpi}_{m+1} (p \cdot {1} + \textbf{b}_m), $$
en posant  $$\textbf{b}_m = \boldsymbol{\varpi}_{m+1}^{p-1} + (p \cdot \textbf{1}) \boldsymbol{\varpi}_{m+1} {\rm ev}_{\boldsymbol{\varpi}_{m+1}, \frak{m}_{m+1}}(G^\Delta(T)), $$
qui est \'{e}l\'{e}ment de l'id\'{e}al $\textbf{W}(\frak{m}_{m+1})$ puisqu'il est divisible par $\boldsymbol{\varpi}_{m+1}$. \end{dem} \findem

\begin{lem} \label{Frobeniusideal} 
Pour tout vecteur de Witt $\textbf{a}$ sur l'anneau $\C_p^\circ$ des entiers de $\C_p$, les composantes du vecteur ${\rm Frob}_{\C_p}(\textbf{a}) - \textbf{a}^\varphi$ sont \'{e}l\'{e}ments de l'id\'{e}al maximal $\C_p^\vee = \lbrace x \in \C_p , \; \vert x \vert < 1 \rbrace$. 
\end{lem}
\begin{dem} 
Consid\'{e}rons le morphisme  $\rho : \C_p^\circ \to \widetilde{\C_p} = \C_p^\circ / C_p^\vee$ de r\'{e}duction modulo l'id\'{e}al maximal $\C_p^\vee$. On note $\widetilde{\varphi}$ l'automorphisme $\widetilde{x} \mapsto \widetilde{x}^p$ de Frobenius du corps r\'{e}siduel $\widetilde{\C_p}$. On a l'\'{e}galit\'{e} $\widetilde{\varphi} \circ \rho = \rho \circ \varphi$ par d\'{e}finition de l'automorphisme de Frobenius $\varphi$ de $\C_p$. On a donc les deux diagrammes commutatifs   
\setlength{\unitlength}{1mm}
\begin{center} 
\begin{picture}(70,40)
\thicklines
\put(0.7,0.3){$\textbf{W}(\C_p^\circ)$}
\put(14,2){\vector(1,0){30}}
\put(46,0.3){$\textbf{W}(\widetilde{\C_p})$}
\put(25,4){${\textbf{W}(\rho)}$}
\put(5,35){\vector(0,-1){30}}
\put(50,35){\vector(0,-1){30}}
\put(0.7,37){$\textbf{W}(\C_p^\circ)$}
\put(42,37){$\textbf{W}(\widetilde{\C_p})$}
\put(14, 38){\vector(1,0){26}}
\put(-8,19){${\rm Frob}_{\C_p^\circ}$}
\put(51,19){${\rm Frob}_{\widetilde{\C_p}}$} 
\put(25,40){${\textbf{W}(\rho)}$}
\end{picture} 
\end{center} 
et 
\setlength{\unitlength}{1mm}
\begin{center} 
\begin{picture}(70,40)
\thicklines
\put(0.7,0.3){$\textbf{W}(\C_p^\circ)$}
\put(14,2){\vector(1,0){30}}
\put(46,0.3){$\textbf{W}(\widetilde{\C_p})$}
\put(25,4){${\textbf{W}(\rho)}$}
\put(5,35){\vector(0,-1){30}}
\put(50,35){\vector(0,-1){30}}
\put(0.7,37){$\textbf{W}(\C_p^\circ)$}
\put(42,37){$\textbf{W}(\widetilde{\C_p})$}
\put(14, 38){\vector(1,0){26}}
\put(-8,19){${\textbf{W}(\varphi)}$}
\put(51,19){${\textbf{W}(\widetilde{\varphi})}$} 
\put(25,40){${\textbf{W}(\rho)}$}
\end{picture} . 
\end{center} 
Soit $\textit{\textbf{a} }$ un vecteur de Witt sur $\C_p^\circ$. Par la commutativit\'{e} du premier diagramme, on a $${\rm Frob}_{\widetilde{\C_p}}(\textbf{W}(\rho)(\textit{\textbf{a} })=\textbf{W}(\rho)({\rm Frob}_{\C_p^\circ}(\textit{\textbf{a} })).$$
D'autre part, en vertu de la remarque \ref{1.8}, et comme $\widetilde{\C_p}$ est un corps de caract\'{e}ristique $p$, on sait que l'automorphisme $\textbf{W}(\widetilde{\varphi})$ co\"incide avec le Frobenius ${\rm Frob}_{\widetilde{\C_p}}$.  Par cons\'{e}quent, la commutativit\'{e} du deuxi\'{e}me diagramme entra\^ine l'\'{e}galit\'{e} $${\rm Frob}_{\widetilde{\C_p}}(\textbf{W}(\rho)(\textit{\textbf{a} })=\textbf{W}(\rho)(\textbf{W}(\varphi)(\textit{\textbf{a} }) ).$$ On en d\'{e}duit que $$\textbf{W}(\varphi)(\textit{\textbf{a} })-{\rm Frob}_{\C_p^\circ}(\textit{\textbf{a} })\in {\rm ker}(\textbf{W}(\rho))=\textbf{W}(\ker \rho)=\textbf{W}(\C_p^\vee). $$

\end{dem}\findem
\begin{cor} \label{diffpV} Pour tout vecteur de Witt $\textbf{a}$ sur l'anneau des entiers $\C_p^\circ$ de $\C_p$, la diff\'{e}rence 
$p\textbf{a}-V_{\C_p}(\textbf{a}^\varphi)$ est \'{e}l\'{e}ment de l'id\'{e}al $\textbf{W}(\C_p^\vee)$, o\`u $\C_p^\vee$ est l'id\'{e}al maximal de $\C_p^\circ$. 
\end{cor}
\begin{dem} Il suffit de remplacer dans l'\'{e}nonc\'{e}  du lemme \ref{Frobeniusideal} le vecteur $\textit{\textbf{a}}$ par le vecteur $V_{\C_p}(\textit{\textbf{a}})$ en utilisant la relation (\ref{relations1}).
\end{dem}\findem
\begin{cor} \label{diffpvn} Pour tout vecteur de Witt $\textbf{a}$ sur l'anneau des entiers $\C_p^\circ$ de $\C_p$, et pour tout entier naturel $m$, la diff\'{e}rence 
$p^m \textbf{a}-V_{\C_p}^m (\textbf{a}^{\varphi^m})$ est \'{e}l\'{e}ment de l'id\'{e}al $\textbf{W}(\C_p^\vee)$. \end{cor} 
\begin{dem} R\'{e}currence sur $m$. \end{dem} \findem

\begin{proposition}\label{thmsurconvpul}
Pour tout entier naturel $m$ et pour tout \'el\'ement $\textbf{a}$ de l'anneau $\textbf{W}(\C_p^\circ)$, la s\'{e}rie exponentielle de Pulita $\theta_m(\textbf{a})$ a un rayon de convergence strictement sup\'{e}rieur {\`a} 1. \end{proposition}
\begin{dem} On sait par la relation (\ref{expPulita}) que la s\'{e}rie $\theta_m(\textit{\textbf{a}})$ est l'image par le morphisme $E$ d'Artin-Hasse du vecteur $\boldsymbol{\varpi}_m \textit{\textbf{a}} - V(\boldsymbol{\varpi}_m \textit{\textbf{a}}^\varphi)$. Or, d'apr{\`e}s le corollaire \ref{Frobvarpi} et le lemme \ref{bm}, on a 
$$ \boldsymbol{\varpi}_m \textit{\textbf{a}} - V(\boldsymbol{\varpi}_m \textit{\textbf{a}}^\varphi) = \boldsymbol{\varpi}_{m+1} \left( \textit{\textbf{b}}_m \textit{\textbf{a}} + p \textit{\textbf{a}} - V(\textit{\textbf{a}}^\varphi) \right). $$
pour un certain vecteur $\textit{\textbf{b}}_m$ appartenant {\`a} l'id\'{e}al $\textbf{W}\left(\frak{m}_\mathcal{O}\right)$ de l'anneau $\textbf{W}(\mathcal{O})$. 
Par le corollaire \ref{diffpV}, on voit que le vecteur $\textit{\textbf{b}}_m \textit{\textbf{a}} + p \textit{\textbf{a}} - V(\textit{\textbf{a}}^\varphi)$ est \'{e}l\'{e}ment de $\textbf{W}(\C_p^\vee)$. D'o{\`u} le r\'{e}sultat en vertu du corollaire 
\ref{surconvergence}. 
\end{dem}\findem

\subsection{Morphismes de Pulita} 

Nous allons maintenant \'{e}tudier des morphismes  de $\textbf{W}(\C_p^\circ)$ dans $\Lambda(\C_p^\circ)$ construits {\`a} l'aide des morphismes $\theta_m$. Ces morphismes seront introduits {\`a} partir d'une g\'{e}n\'{e}ralisation des exponentielles de Pulita. 

\subsubsection{S\'{e}ries exponentielles de Pulita g\'{e}n\'{e}ralis\'{e}es} 
\begin{notation} \label{thetams} Pour tout vecteur de Witt $\textit{\textbf{a}}$ sur $\C_p$, et pour tout couple $(m,s)$ d'entiers naturels tel que $s \ge 1$, on note $\theta_{m,s}(\textit{\textbf{a}})$ la s\'{e}rie 
\begin{equation} \label{factorisationthetams} \theta_{m,s}(\textit{\textbf{a}}) = \prod_{i=0}^{s-1} \theta_m(\textit{\textbf{a}}^{\varphi^i}) \circ x^{p^i} \; . \end{equation} 
\end{notation}
Comme le facteur $\theta_m(\textit{\textbf{a}}^{\varphi^i}) \circ x^{p^i}$ figurant dans le membre de droite de l'identit\'{e} (\ref{factorisationthetams}) est l'image par le morphisme d'Artin-Hasse du vecteur $V^i\left(\boldsymbol{\varpi}_m \textit{\textbf{a}}^{\varphi^i} - V\left(\boldsymbol{\varpi}_m \textit{\textbf{a}}^{\varphi^{i+1}} \right) \right)$, on voit que 
$$ \theta_{m,s}(\textit{\textbf{a}}) = E\left(\boldsymbol{\varpi}_m \textit{\textbf{a}} - V^s\left(\boldsymbol{\varpi}_m \textit{\textbf{a}}^{\varphi^s} \right) \right) \; .$$
En particulier, l'application $\theta_{m,s}$ est un morphisme du groupe additif $\textbf{W}(\C_p)$ dans le groupe multiplicatif $\Lambda(\C_p)$. On remarque que $\theta_{m,1}=\theta_{m}$, ce qui justifie le nom d'exponentielle de Pulita g\'{e}n\'{e}ralis\'{e}e. 

\begin{lem} \label{Ithetams} Soit $m$ et $s$ deux entiers naturels tels que $s \ge 1$, $I$ un id\'{e}al de l'anneau $\C_p^\circ$, et $\textit{\textbf{a}}$ un \'{e}l\'{e}ment de $\textbf{W}(I)$. Alors tous les coefficients de la s\'{e}rie $\theta_{m,s}(\textit{\textbf{a}}) - 1$ sont \'{e}l\'{e}ments de l'id\'{e}al $\pi_m I$. \end{lem} 
\begin{dem} R\'{e}sulte de l'identit\'{e} (\ref{factorisationthetams}) et des lemmes \ref{stabiliteidealFrobabsolu} et \ref{Itheta}. \end{dem} \findem

\begin{lem}\label{pulita V} Pour tout vecteur de Witt $\textit{\textbf{a}}$ sur $\C_p$, et pour tout triplet $(m,s,k)$ d'entiers naturels tel que $s \ge 1$ et $k \le m$, on a $\theta_{m,s}(V^k(\textit{\textbf{a}}))=\theta_{m-k,s}(\textit{\textbf{a}})\circ x^{p^k}.$ 
 \end{lem}
 \begin{dem}
 Par d\'{e}finition $$\theta_{m,s}(V^k(\textit{\textbf{a}}))=\prod_{i=0}^{s-1} \theta_{m}(V^k\textit{\textbf{a}}^{\varphi^{i}})\circ x^{p^{i}},$$
  alors, d'apr{\`e}s le lemme \ref{pulitadecal} $$\theta_{m,s}(V^k(\textit{\textbf{a}}))=\prod_{i=0}^{s-1} \theta_{m-k}(\textit{\textbf{a}}^{\varphi^{i}})\circ x^{p^{k}}\circ x^{p^{i}},$$
  et comme $x^{p^{k}}\circ x^{p^{i}}=x^{p^{i}}\circ x^{p^{k}},$ on a donc$$\theta_{m,s}(V^k(\textit{\textbf{a}}))=\prod_{i=0}^{s-1} \theta_{m-k}(\textit{\textbf{a}}^{\varphi^{i}})\circ x^{p^{i}}\circ x^{p^{k}},$$
on obtient le r\'{e}sultat d\'{e}sir\'{e}:
  $$\theta_{m,s}(V^k(\textit{\textbf{a}}))=\theta_{m-k,s}(\textit{\textbf{a}})\circ x^{p^{k}}.$$\end{dem}\findem

\begin{lem}\label{thetataums} Soit $m$ et $s$ deux entiers naturels tels que $s \ge 1$. Si $t \in \C_p$ est une racine de l'unit\'{e} d'ordre premier {\`a} $p$, alors 
$$ \theta_{m,s}(\boldsymbol{\tau}(t) ) = \theta_{m,s}(\textbf{1}) \circ (tx) \; .$$ \end{lem} 
\begin{dem} \end{dem} 
R\'{e}sulte de la d\'{e}finition du morphisme $\theta_{m,s}$ et des lemmes \ref{imageFrobracine} et \ref{thetamtaut}. \findem 

Lorsque $s$ varie, les s\'{e}ries $\theta_{m,s}$ satisfont une propri\'{e}t de transitivit\'{e}, exprim\'{e}e par le lemme suivant. 


\begin{lem}\label{thetasr_fonction_thetas} Soit $\textit{\textbf{a}}$ un vecteur de Witt sur $\C_p$, et $m, r \ge 1, s \ge 1$ trois entiers naturels. On a la factorisation 
$$\theta_{m, sr}(\textit{\textbf{a}}) = \prod_{j=0}^{r-1}\theta_{m,s}(\textit{\textbf{a}}^{\varphi^{js}})\circ x^{p^{sj}} \; . $$\end{lem}
\begin{dem} Comme le facteur $\theta_{m,s}(\textit{\textbf{a}}^{\varphi^{js}}) \circ x^{p^{sj}}$ figurant dans le membre de droite de l'identit\'{e} {\`a} montrer est l'image par le morphisme d'Artin-Hasse $E$ du vecteur $V^{sj}\left(\boldsymbol{\varpi}_m \textit{\textbf{a}}^{\varphi^{js}} - V^s\left(\boldsymbol{\varpi}_m \textit{\textbf{a}}^{\varphi^{(j+1)s}} \right) \right)$, on voit que leur produit est l'image par $E$ de la somme 
\begin{align*} \sum_{j=0}^{r-1} V^{sj}\left(\boldsymbol{\varpi}_m \textit{\textbf{a}}^{\varphi^{js}} - V^s\left(\boldsymbol{\varpi}_m \textit{\textbf{a}}^{\varphi^{(j+1)s}} \right) \right) & = \sum_{j=0}^{r-1} \left( V^{sj}\left(\boldsymbol{\varpi}_m \textit{\textbf{a}}^{\varphi^{js}} \right) - V^{s(j+1)}\left(\boldsymbol{\varpi}_m \textit{\textbf{a}}^{\varphi^{(j+1)s}} \right) \right)  = \\ & = \boldsymbol{\varpi}_m \textit{\textbf{a}} - V^{sr}\left(\boldsymbol{\varpi}_m \textit{\textbf{a}}^{\varphi^{rs}} \right) \end{align*}
Or on sait que justement
$$ \theta_{m,sr}(\textit{\textbf{a}}) = E\left(\boldsymbol{\varpi}_m \textit{\textbf{a}} - V^{sr}\left(\boldsymbol{\varpi}_m \textit{\textbf{a}}^{\varphi^{sr}} \right) \right) \; .$$

\end{dem}\findem

\begin{proposition}\label{surconvergencethetams}  Pour tout couple $(m,s)$ d'entiers naturels tel que $s \ge 1$, et pour tout vecteur de Witt $\textit{\textbf{a}}$ sur l'anneau $\C_p^\circ$ des entiers de $\C_p$, la s\'{e}rie  $\theta_{m,s}(\textit{\textbf{a}})$ a un rayon de convergence $> 1$.  \end{proposition} 
\begin{dem} D'apr{\`e}s la proposition \ref{thmsurconvpul}, tous les facteurs du membre de droite de l'identit\'{e} (\ref{factorisationthetams}) sont des s\'{e}ries de rayon de convergence $> 1$. \end{dem} \findem

\subsubsection{Les morphismes $\overline{\theta}_{m,s}$} 
\begin{lem} \label{noyautheta} L'id\'{e}al $V^{m+1}\textbf{W}(\C_p)$ est contenu dans le noyau du morphisme $\theta_{m,s}$. \end{lem}
\begin{dem} Soit $\textit{\textbf{a}}$ un vecteur de Witt \'{e}l\'{e}ment de $V^{m+1}\textbf{W}(\C_p)$. On peut \'{e}crire $\textit{\textbf{a}}= \left( a_n \right)_{n \in \mathbb{N}}$ avec $ a_i=0$ pour tout indice $i \le m$, et donc ${\rm fant}_n(\textit{\textbf{a}}) = 0$ pour tout entier $n \in \lbrack 0 .. m \rbrack$. Comme ${\rm fant}_n(\boldsymbol{\varpi}_m) = \pi_{m-n} = 0$ pour tout entier $n > m$, on en d\'{e}duit que toutes les composantes fant\^{o}mes du vecteur $\boldsymbol{\varpi}_m \textit{\textbf{a}}$ sont nulles. Comme $p$ est inversible dans le corps $\C_p$, l'application fant\^{o}me ${\rm fant}_{\C_p}$ est bijective. Donc $\boldsymbol{\varpi}_m \textit{\textbf{a}} = \textbf{0}$. D'autre part, en vertu de la remarque \ref{Vnideal}, le vecteur $\textit{\textbf{a}}^{\varphi^s}$ est aussi \'{e}l\'{e}ment de $V^{m+1}\textbf{W}(\C_p)$, donc on a aussi $\boldsymbol{\varpi}_m \textit{\textbf{a}}^{\varphi^s} = \textbf{0}$.

Finalement, on voit que $\boldsymbol{\varpi}_m \textit{\textbf{a}} - V_{\C_p}^s \left( \boldsymbol{\varpi}_m \textit{\textbf{a}}^{\varphi^s} \right) = \textbf{0}$, et donc $\theta_{m,s}(\textit{\textbf{a}}) = 1$. 
\end{dem} \findem

\begin{notation} On d\'{e}finit $\Lambda^*(\C_p^\circ)$ comme le sous-groupe de $\Lambda(\C_p)$ constitu\'{e} des s\'{e}ries enti{\`e}res formelles {\`a} coefficients \'{e}l\'{e}ments de $\C_p^\circ$ dont le terme constant est \'{e}gal {\`a} 1, et  dont le rayon de convergence, ainsi que celui de leur inverse, est $>1$. \end{notation} 

\begin{lem} \label{imagetheta} L'image par le morphisme $\theta_{m,s}$ du sous-anneau $\textbf{W}(\C_p^\circ)$ est contenue dans le sous-groupe $\Lambda^*(\C_p^\circ)$. \end{lem}
\begin{dem}  Soit un vecteur de Witt $\textit{\textbf{a}}$ sur $\C_p^\circ$. D'apr{\`e}s la proposition \ref{surconvergencethetams}, la s\'{e}rie $\theta_{m,s}(\textit{\textbf{a}})$ a un rayon de convergence strictement sup\'{e}rieur {\`a} $1$. D'autre part son inverse est $\theta_{m,s}(-\textit{\textbf{a}}) $ qui a un rayon de convergence $> 1$ car $\textit{\textbf{-a}}$ est aussi un vecteur de Witt sur $\C_p^{\circ}$.\end{dem} \findem

\begin{proposition} \label{morphismePulita} 
Le morphisme $\theta_{m,s} : \textbf{W}(\C_p) \to \Lambda(\C_p)$ induit un unique morphisme 
$$ \overline{\theta}_{m,s} : \textbf{W}_{m+1}(\C_p^\circ) \to \Lambda^*(\C_p^\circ) \; .$$
\end{proposition}
\begin{dem} D'apr{\`e}s le lemme \ref{imagetheta}, la restriction du morphisme $\theta_{m,s}$ au sous-anneau $\textbf{W}(\C_p^\circ)$ induit un unique morphisme de $\textbf{W}(\C_p^\circ)$ dans $\Lambda^*(\C_p^\circ)$. Par  le lemme \ref{noyautheta}, et comme l'anneau $\textbf{W}_{m+1}(\C_p^\circ)$ est isomorphe au quotient $\textbf{W}(\C_p^\circ) / V^{m+1} \textbf{W}(\C_p^\circ)$ d'apr{\`e}s la proposition \ref{Wittlongueurfinie}, le r\'{e}sultat suit.\end{dem}\findem

\section{Expression des caract{\`e}res de $W_{\ell}(\mathbb{F}_{q})$}

Dans cette section, on fixe un entier naturel $s \ge 1$, et on pose $q = p^s$. 

\subsection{Caract\`ere de Teichm\"uller d'un corps fini}

On note $\mu_{q-1}$ le groupe des racines de l'unit\'{e} dans $\C_p$ d'ordre divisant $q-1$. Le sous-anneau $\Z_p\left\lbrack \mu_{q-1} \right\rbrack$ de $\C_p$ est l'anneau des entiers de l'extension $\Q_p \left(\mu_{q-1} \right)$ non ramifi\'{e}e du corps $\Q_p$ des nombres $p$-adiques, de sorte que son corps r\'{e}siduel est d'ordre $q$. Il existe donc au moins un (et en fait exactement $s$) morphismes d'anneaux, dits de r\'{e}duction, de $\Z_p\left\lbrack \mu_{q-1} \right\rbrack$ sur le corps $\F_q$. On supposera qu'on en a choisi un, not\'{e} $x \mapsto \widetilde{x}$. 
Ceci nous permet d'introduire un caract{\`e}re de Teichm\"uller du corps $\F_q$ comme suit. 

\begin{notation} 
Pour $x \in \F_q$, on note $\Te(x)$ l'unique \'{e}l\'{e}ment de $\Z_p\left\lbrack \mu_{q-1} \right\rbrack$ tel que $\Te(x)^{q} = \Te(x)$ et $\widetilde{\Te(x)} = x$. 
\end{notation} 
L'existence et l'unicit\'{e} de l'\'{e}l\'{e}ment $\Te(x)$ r\'{e}sulte imm\'{e}diatement du lemme de Hensel. On a ainsi d\'{e}fini une application 
$$ \Te : \F_q \to \C_p$$
qui est un g\'{e}n\'{e}rateur du groupe des caract{\`e}res multiplicatifs de $\F_q$. C'est cette application qui est appel\'{e} caract{\`e}re de Teichm\"uller. 

Une propri\'{e}t\'{e} utile du caract{\`e}re de Teichm\"uller est la suivante. 
\begin{proposition} \label{FrobTeich} Soit $x \in \F_q$. Alors on a l'identit\'{e} 
$$ {\rm Frob}_{\mathcal{O} \lbrack \mu_{q-1} \rbrack}(\boldsymbol{\tau}(\Te(x)) = \boldsymbol{\tau}(\Te(x^p)) \; .$$ \end{proposition} 
\begin{dem} Cas particulier de la proposition \ref{Frobtau}. \end{dem} \findem

\begin{lem} \label{Frobeniuss} Soit $K$ le corps des fractions de l'anneau $\mathcal{O} = \cup_{m \in \mathbb{N}} \mathcal{O}_m$, $\varphi$ un automorphisme de Frobenius absolu satisfaisant les conditions de la proposition \ref{Frobeniusabsolu}. Alors, pour tout \'{e}l\'{e}ment $x \in \F_q$, on a $\varphi(\Te(x)) = \Te(x^p)$. De plus la puissance $\varphi^s$ d'ordre $s$ de $\varphi$ fixe tous les \'{e}l\'{e}ments du corps $K\lbrack \mu_{q-1} \rbrack$. \end{lem}
\begin{dem} L'\'{e}l\'{e}ment $\varphi(\Te(x))$ appartient {\`a} l'anneau $\Z_p\lbrack \mu_{q-1} \rbrack$, est \'{e}gal {\`a} sa puissance d'ordre $q$, et sa r\'{e}duction est \'{e}gale {\`a} $x^p$, puisqu'on sait par d\'{e}finition de $\varphi$ que $\vert \varphi(\Te(x)) - \Te(x)^p \vert < 1$. Ceci montre que $\varphi(\Te(x)) = \Te(x^p)$. Par d\'{e}finition, on sait que $\varphi$ fixe les points de $K$, il suffit donc de v\'{e}rifier que $\varphi^s(t) = t$ pour tout \'{e}l\'{e}ment $t$ de $\mu_{q-1}$. Or  on a $t = \Te(x)$, o{\`u} $x = \widetilde{t}$ est la r\'{e}duction de $t$ dans $\F_q$.  Il en r\'{e}sulte que $\varphi^s(t) = \varphi^s(\Te(x)) = \Te(x^{p^s}) = \Te(x) = t$. \end{dem} \findem

 \subsection{Evaluation des s\'{e}ries de Pulita} 

Pour toute s\'{e}rie $G \in \C_p \lbrack \! \lbrack x \rbrack \! \rbrack$ dont le rayon de convergence est $> 1$, et pour tout $z \in \C_p^\circ$, on note $G(z)$ la valeur de $G$ en $z$. On obtient ainsi un morphisme $G \mapsto G(z)$ du groupe $\Lambda^*(\C_p^\circ)$ dans le groupe $\C_p^\times$, que nous appelons morphisme d'\'{e}valuation en $z$. En particulier, si $\textit{\textbf{a}}$ est un vecteur de Witt sur $\C_p^\circ$, on peut calculer $\theta_m(\textit{\textbf{a}})(z)$ en tout point $z$ de l'anneau des entiers de $\C_p$. 
\begin{lem} \label{2.1} Soit $s \ge 1$ un entier naturel et $t \in \mathbb{C}_p$ tel que $t^q = t$, o{\`u} on a pos\'{e} $q=p^s$. Si $\textit{\textbf{a}}$ est un vecteur de Witt dont toutes les composantes sont \'{e}l\'{e}ments de l'id\'{e}al maximal $\frak{m}_{\mathcal{O}\lbrack \mu_{q-1} \rbrack}$ de l'anneau $\mathcal{O}\lbrack \mu_{q-1} \rbrack = \cup_{m \in \mathbb{N}} \mathcal{O}_m\lbrack \mu_{q-1} \rbrack$, alors, pour tout entier $\ell \ge 1$, on a  $ \theta_{\ell -1, s}(\textbf{a})(t) = 1 \; .$ \end{lem}
\begin{dem} Soit $\textit{\textbf{a}}\in \mathbf{W}(\frak{m}_{\mathcal{O}\lbrack \mu_{q-1} \rbrack})$. Comme l'automorphisme $\varphi^s$ laisse fixe tout \'{e}l\'{e}ment du corps $K\lbrack \mu_{q-1} \rbrack$ par le lemme \ref{Frobeniuss},  on a $\textit{\textbf{a}}^{\varphi^s} =\textit{\textbf{a}}$. On sait par le corollaire \ref{surconvergence}  que $E\left( \boldsymbol{\varpi_{\ell -1}} \textit{\textbf{a}}\right )=E\left( \boldsymbol{\varpi_{\ell -1}} \textit{\textbf{a}}^{\varphi^s}\right )$ a un rayon de convergence $>1$. Puisque $t=t^{p^s}$, les deux nombres $E\left( \boldsymbol{\varpi_{\ell -1}} \textit{\textbf{a}}\right )(t)$ et $E\left( \boldsymbol{\varpi_{\ell -1}} \textit{\textbf{a}}^{\varphi^s}\right )(t^{p^s})$ sont \'{e}gaux, d'o{\`u} on d\'{e}duit que leur rapport $\theta_{\ell -1,s}(\textit{\textbf{a}})(t)$ est 1. 
\end{dem} \findem
\begin{lem}\label{2.2} Soit $z$ un \'{e}l\'{e}ment de $\C_p^\circ$ et $\ell \ge 1$ un entier naturel. La valeur en $z$ de la s\'{e}rie de Pulita $\theta_{\ell-1}(\textit{\textbf{1}}) \;$  est congrue {\`a} $1 + \pi_{\ell-1} z$ modulo l'id\'{e}al $\pi_{\ell-1}^2 \C_p^\circ$.\end{lem}
\begin{dem} Consid\'{e}rons l'anneau $A=\mathbb{Z}_p \lbrack \! \lbrack T \rbrack \! \rbrack$, son id\'{e}al $T^2\mathbb{Z}_p[[T]]$ et le vecteur de Witt  $\textit{\textbf{w}}-V(\textit{\textbf{w}})-\boldsymbol{\tau}(T)$ de $\textbf{W}(A)$, o{\`u}  $\textit{\textbf{w}}$ est l'unique vecteur de Witt sur l'anneau $\mathbb{Z}_p\lbrack \! \lbrack T \rbrack \! \rbrack$ dont la suite des composantes fant\^omes est la suite $\langle F^{\circ n}(T) \rangle_{n \in \mathbb{N}}$ des puissances de composition de la s\'erie de Lubin-Tate $F(T)$ (Lemme \ref{puissancescompositionF}). Utilisons le lemme suivant. 

\begin {lem} \label{2.3} Toutes les composantes fant\^{o}mes du vecteur $\textbf{w}-V(\textbf{w})-\boldsymbol{\tau}(T)$ sont \'{e}l\'{e}ments de l'id\'{e}al $T^2\mathbb{Z}_p[[T]]$.  \end{lem}\begin{dem} On a 
$$ \forall j \in \mathbb{N}, \qquad {\rm fant}_j\left( \textbf{w}-V(\textbf{w})-\boldsymbol{\tau}(T) \right) = \begin{cases} F^{\circ j}(T) - p F^{\circ(j-1)}(T) - T^{p^j} & \mbox{ si } j \ge 1 \\ 
0 & \mbox{ si } j = 0 \end{cases}. $$
Comme $F(T)$ est une s\'{e}rie de Lubin-Tate, il existe une s\'{e}rie $G(T) \in \Z_p \lbrack \! \lbrack T \rbrack \! \rbrack$ telle que $F(T) = pT + T^p + pT^2 G(T)$, donc, pour $j \ge 1$, on a 
$$ F^{\circ j}(T) - p F^{\circ(j-1)}(T) - T^{p^j} = F^{\circ(j-1)}(T)^p + p F^{\circ(j-1)}(T)^2 G(F^{\circ(j-1)}(T)) - T^{p^j} \; , $$ qui appartient {\`a} $T^2 \Z_p \lbrack \! \lbrack T \rbrack \! \rbrack $ puisque la s\'{e}rie enti{\`e}re formelle$F^{\circ(j-1)}(T)$ est divisible par $T$.\end{dem}\findem.

On d\'{e}duit du lemme \ref{2.3} et de la proposition \ref{idealWS} que le vecteur de Witt $\textit{\textbf{w}}-V(\textit{\textbf{w}})-\boldsymbol{\tau}(T)$ est \'{e}l\'{e}ment de l'id\'{e}al $\mathbf{W}(T^2 \Z_p \lbrack \! \lbrack T \rbrack \! \rbrack)$. Comme l'homomorphisme $\varepsilon_{\pi_{\ell-1}} : \Z_p \lbrack \! \lbrack T \rbrack \! \rbrack \to \mathcal{O}_{\ell-1} = \Z_p\lbrack \pi_{\ell-1} \rbrack$ de sp\'{e}cialisation en $\pi_{\ell-1}$ d\'{e}fini par la relation (\ref{specialisation}) envoie l'id\'{e}al $T^2 \Z_p \lbrack \! \lbrack T \rbrack \! \rbrack$ dans l'id\'{e}al $\pi_{\ell-1}^2 \mathcal{O}_{\ell-1}$, l'image par le morphisme $\textbf{W}(\varepsilon_{\pi_{\ell-1}})$ du vecteur $\textit{\textbf{w}}-V(\textit{\textbf{w}})-\boldsymbol{\tau}(T)$, c'est-{\`a}-dire  $\boldsymbol{\varpi}_{\ell-1}-V(\boldsymbol{\varpi}_{\ell-1})-\boldsymbol{\tau}(\pi_{\ell-1})$, appartient {\`a} l'id\'{e}al $\textbf{W}(\pi_{\ell-1}^{2}\mathcal{O}_{\ell-1})$. On peut donc \'{e}crire $\boldsymbol{\varpi}_{\ell-1}-V(\boldsymbol{\varpi}_{\ell-1})$ sous la forme $\boldsymbol{\tau}(\pi_{\ell-1})+\textit{\textbf{y}}$, avec $\textit{\textbf{y}}\in \textbf{W}(\pi_{\ell-1}^{2}\mathcal{O}_{\ell-1})$, de sorte que $E\left(\textit{\textbf{y}}\right)=1+\pi_{\ell-1}^2H(x)$ pour une certaine s\'{e}rie $H\in \mathcal{O}_{\ell-1} \lbrack \! \lbrack x \rbrack \! \rbrack$. D'autre part $E\left(\boldsymbol{\tau}(\pi_{\ell-1})\right)=AH(\pi_{\ell-1} x)=1+\pi_{\ell-1}x +\pi_{\ell-1}^2H_1(x),$ o{\`u}
$H_1(x)$ est une s\'{e}rie enti{\`e}re formelle \'{e}l\'{e}ment de $\mathcal{O}_{\ell-1} \lbrack \! \lbrack x \rbrack \! \rbrack$. On a donc 
$$ E\left(\boldsymbol{\tau}(\pi_{\ell-1}) +\textit{\textbf{y}}\right) = \left(1+\pi_{\ell-1}x +\pi_{\ell-1}^2H_1(x) \right) \left( 1+\pi_{\ell-1}^2H(x) \right) = 1 + \pi_{\ell-1} x + \pi_{\ell-1}^2 H_2(x) , $$ avec 
$H_2(x) = H(x) + H_1(x) + \pi_{\ell-1}^2 H(x) H_1(x) $ qui appartient {\`a} $\mathcal{O}_{\ell-1}\lbrack \! \lbrack x \rbrack \! \rbrack$. On voit donc que $\theta_{\ell-1}(\textbf{1}) = E\left(\boldsymbol{\varpi}_{\ell-1} - V(\boldsymbol{\varpi}_{\ell-1}) \right) = 1 + \pi_{\ell-1} x + \pi_{\ell-1}^2 H_2(x)$. 
Comme la s\'{e}rie $\theta_{\ell-1}(\textbf{1})$ a un rayon de convergence $> 1$, il en est de m{\^e}me de la s\'{e}rie $H_2(x)$, et on a 
$\theta_{\ell-1}(\textbf{1})(z) = 1 + \pi_{\ell-1} z + \pi_{\ell-1}^2 H_2(z)$. 
Ce qui implique le r\'{e}sultat escompt\'{e}.
 \end{dem}\findem
 \begin{cor}\label{congthetamsen1} Soit $z$ un \'{e}l\'{e}ment de $\C_p^\circ$ et $(\ell, s)$ un couple d'entiers naturels non nuls. La valeur en $z$ de la s\'{e}rie $\theta_{\ell-1,s}(\textbf{1})$ est congrue {\`a} $1 + \pi_{\ell-1} \left( \sum_{j=0}^{s-1} z^{p^j} \right)$ modulo l'id\'{e}al $\pi_{\ell-1}^2 \C_p^\circ$. \end{cor}
 \begin{dem} Par d\'{e}finition $$\theta_{\ell-1,s}(\textit{\textbf{1}})=\prod_{i=0}^{s-1}(\theta_{\ell-1}(\textit{\textbf{1}})\circ x^{p^{i}}),$$ 
  d'o{\`u}, vu le lemme \ref{2.2}, l'\'{e}galit\'{e}
  $$\theta_{\ell-1,s}(\textit{\textbf{1}})(z)\equiv \prod_{i=0}^{s-1}\left(1+\pi_{\ell-1}z^{p^{i}}\right)\, \pmod{\pi_{\ell-1}^2 \C_p^{\circ} }\; .$$
En d\'{e}veloppant le dernier produit, on obtient la congruence {\`a} montrer. 
\end{dem}\findem
 
 \begin{proposition} \label{evaluationthetamstaua} Soit $\ell$ et $s$ deux entiers naturels non nuls, et $a, z$ deux \'{e}l\'{e}ments de l'anneau $\C_p^\circ$. Alors il existe un entier naturel $r$, et un \'{e}l\'{e}ment $t$ de $\C_p^\circ$ tel que $t^{p^r} - t = 0$  et $\vert a - t \vert < 1$. En posant $b = a - t$ si $\vert a - t \vert > \vert \pi_{\ell-1} \vert$, et $b = \pi_{\ell - 1}$ dans le cas contraire, on a la congruence 
$$ \theta_{\ell-1,s}(\boldsymbol{\tau}(a))(z) \equiv 1 + \pi_{\ell-1} \left( \sum_{j=0}^{s-1} a^{p^j} z^{p^j} \right) \pmod{ \pi_{\ell-1} b \C_p^\circ} \; .$$
 \end{proposition}
 \begin{dem} La classe r\'{e}siduelle $\widetilde{a}$ de $a$ dans le corps r\'{e}siduel $ \widetilde{\C_p}$ de $\C_p $ est alg\'{e}brique sur $\F_p$, donc appartient {\`a} un corps fini isomorphe {\`a} $\F_{p^r}$ pour un certain entier rationnel $r \ge 1$, d'o{\`u} l'existence de l'entier naturel $r$ tel que $\widetilde{a}^{p^r} - \widetilde{a} = 0$. D'apr{\`e}s le lemme de Hensel, le polyn\^{o}me $x^r - x$ admet une unique racine $t  \in \C_p^\circ$ telle que $\vert a - t \vert < 1$. 
 Passons maintenant au calcul de $ \theta_{\ell-1,s}(\boldsymbol{\tau}(a))$.
Comme on a dans tous les cas l'in\'{e}galit\'{e} $\vert a - t \vert \le \vert b \vert$, le corollaire \ref{Wdiff} entra\^ine la congruence  $\boldsymbol{\tau}(a) \equiv \boldsymbol{\tau}(t) \pmod{\textbf{W}(b\C_p^\circ)}$. D'apr{\`e}s le lemme \ref{Ithetams}, on en d\'{e}duit que la s\'{e}rie $\theta_{\ell - 1,s}(\boldsymbol{\tau}(a) - \boldsymbol{\tau}(t)) - 1$ est \'{e}l\'{e}ment de l'id\'{e}al $ \pi_{\ell-1} b x \C_p^\circ\lbrack \! \lbrack x \rbrack \! \rbrack$.  Par cons\'{e}quent 
  $$\theta_{\ell-1,s}(\boldsymbol{\tau}(a)) = \theta_{\ell-1,s}(\boldsymbol{\tau}(a)-\boldsymbol{\tau}(t)) \theta_{\ell -1, s}(\boldsymbol{\tau}(t) ) \equiv \theta_{\ell-1,s}(\boldsymbol{\tau}(t)) \pmod{\pi_{\ell -1} b x \C_p^\circ \lbrack \! \lbrack x \rbrack \! \rbrack },$$ 
d'o{\`u} l'on tire la congruence 
\begin{equation} \label{equcongr0}
 \theta_{\ell - 1, s}(\boldsymbol{\tau}(a))(z) \equiv \theta_{\ell - 1, s}(\boldsymbol{\tau}(t))(z) \pmod{\pi_{\ell -1} b \C_p^\circ} \end{equation}
  D'apr{\`e}s  le lemme \ref{thetataums}, on a $$\theta_{\ell-1,s}(\boldsymbol{\tau}(t))= \theta_{\ell-1,s}(\textbf{1}) \circ (tx) \; . $$  Puisque $\pi_{\ell - 1} b $ est toujours un diviseur de $\pi_{\ell - 1}^2$ dans l'anneau $\C_p^\circ$, on en d\'{e}duit d'apr{\`e}s le corollaire \ref{congthetamsen1}
\begin{equation} \label{equcongr1}
\theta_{\ell-1,s}(\boldsymbol{\tau}(t))(z) \equiv 1 + \pi_{\ell-1} \left( \sum_{j=0}^{s-1} t^{p^j}z^{p^j} \right) \pmod{\pi_{\ell-1} b \C_p^\circ} \; .\end{equation}
Le r\'{e}sultat d\'{e}sir\'{e} suit des congruences (\ref{equcongr0}) et (\ref{equcongr1}) puisque $t \equiv a \pmod{b \C_p^\circ}$. 
 \end{dem}\findem

 \begin{proposition}\label{valeurthetals a}Soit $\textit{\textbf{a}}=(a_i)_{i\geq 0}$ un vecteur de Witt sur la cl\^{o}ture de l'extension non ramifi\'{e}e maximale de $K_{\ell -1} = \Q_p\lbrack \pi_{\ell-1} \rbrack$ et $z$ un entier de $\C_p^\circ $. Alors la valeur $\theta_{\ell-1,s}(\textit{\textbf{a}})(z)$ est congrue {\`a} $1+\pi_{\ell-1} \sum_{j=0}^{s-1} {a_0}^{p^j} z^{p^j}$ modulo l'id\'{e}al $\pi^2_{\ell-1}\C_p^\circ.$ \end{proposition} 
 \begin{dem} 
D'apr\'{e}s la proposition \ref{developpeWitt}, on a  $$\theta_{{\ell -1},s}(\textit{\textbf{a}})=\theta_{{\ell -1},s}\left(\sum_{i=0}^\infty V^{i}(\boldsymbol{\tau} (a_i) ) \right)=\theta_{{\ell -1},s}\left(\sum_{i=0}^{\ell-1}V^{i}(\boldsymbol{\tau} (a_i)\right)+\theta_{{\ell -1},s}\left(\sum_{i=\ell}^\infty V^{i}(\boldsymbol{\tau} (a_i) ) \right), $$ 
ce qui implique,  en vertu du lemme \ref{noyautheta} : 
 $$ \theta_{{\ell -1},s}(\textit{\textbf{a}})= \theta_{{\ell -1},s}\left(\sum_{i=0}^{\ell-1}V^{i}(\boldsymbol{\tau} (a_i) ) \right)$$ ce qui, selon la proposition  \ref{thetataums}, conduit {\`a}   
  $$ \theta_{{\ell -1},s}(\textit{\textbf{a}})=\prod_{i=0}^{\ell-1}\theta_{{\ell -1-i},s}(\boldsymbol{\tau} (a_i)) \circ x^{p^{i}}$$
Utilisons la proposition \ref{evaluationthetamstaua} : pour tout entier $i \in \lbrack 0 .. \ell-1 \rbrack$, il existe un entier naturel $r_i$ et un \'{e}l\'{e}ment $t_i \in \C_p^\circ$ tel que $t_i^{p^{r_i}} - t_i = 0$ et $\vert a_i - t_i \vert < 1$. Posons $b_i = \pi_{\ell-1-i}$ si $\vert a_i - t_i \vert \le \vert \pi_{\ell-1-i}\vert$ et $b_i = a_i - t_i$ dans le cas contraire. Suivant la proposition \ref{evaluationthetamstaua}, on a alors 
$$ \theta_{\ell-1-i,s}(\boldsymbol{\tau}(a_i))(z) \equiv 1 + \pi_{\ell-1-i} \left( \sum_{j=0}^{s-1} a_i^{p^j} z^{p^j} \right) \pmod{\pi_{\ell-1-i} b_i \C_p^\circ} \; .$$ Observons que l'extension $\Q_p(t_i)/\Q_p$ est non ramifi\'{e}e, ce qui entra\^ine que l'extension $\Q_p(t_i, a_i)/K_{\ell-1}$ est non ramifi\'{e}e sur $K_{\ell-1}$. Par cons\'{e}quent $\vert t_i - a_i \vert \le \vert \pi_{\ell-1} \vert$, et donc $\vert b_i \vert = \max( \vert t_i - a_i \vert, \vert \pi_{\ell-1-i} \vert) \le \vert \pi_{\ell-1} \vert$. 
 Pour $i > 0$, on a $\pi_{\ell-1-i}\in \pi_{\ell-1}^2\C_p^{\circ}$, et on obtient  
 $$\theta_{{\ell -1-i},s}(\boldsymbol{\tau}(a_i))(z) \equiv  1 \pmod{ \pi_{\ell-1}^2\C_p^\circ} \;. $$
On conclut que 
$$ \theta_{\ell-1,s}(\textit{\textbf{a}})(z) \equiv 1 + \pi_{\ell-1} \left( \sum_{j=0}^{s-1} a_0^{p^j} z^{p^j} \right) \pmod{\pi_{\ell-1}^2 \C_p^\circ}. $$ 
\end{dem} \findem

\subsection{Construction d'un caract{\`e}re additif de $W_\ell(\mathbb{F}_q)$}

\setlength{\parindent}{5mm}
\begin{proposition} \label{2.4} Soit $t \in \mathbb{Z}_p\lbrack \mu_{q-1} \rbrack$ tel que $t^{p^{s}} = t$. Le compos\'{e} du morphisme $\overline{\theta}_{\ell -1,s}$ suivi du morphisme ${\rm eval}_t :\Lambda^*(\C_p^\circ) \to \C_p^\times$ d'\'{e}valuation en $t$, induit un unique morphisme $\psi_{\ell, s, t}$  de $\textbf{W}_\ell(\F_q)$ dans $\C_p^\times$. De plus, si ${\rm Tr}_{\Q_p(\mu_{q-1})/\Q_p}(t) \not\in p \Z_p$, alors : \\
\indent 1. pour tout caract{\`e}re additif $\psi : \textbf{W}_\ell(\F_q) \to \C_p^\times$, il existe un unique $\textit{\textbf{a}} \in \textbf{W}_\ell(\F_q)$ tel qu'on ait $\psi(\textit{\textbf{y}}) = \psi_{\ell, s, t}(\textit{\textbf{ay}})$ pour tout vecteur $\textit{\textbf{y}} \in \textbf{W}_\ell(\F_q)$ ; \\
\indent 2. l'image de $\psi_{\ell, s, t}$ est exactement le groupe $\mu_{p^\ell}$ des racines de l'unit\'{e} d'ordre divisant $p^\ell$. \end{proposition}

\begin{sloppypar} \begin{dem} La restriction de l'application ${\rm eval}_t \circ \overline{\theta}_{\ell-1,s}$ {\`a} l'anneau $\textbf{W}_\ell(\mathcal{O}\lbrack \mu_{q-1}\rbrack)$ est un morphisme de groupes. D'apr{\`e}s le lemme \ref{2.1}, le noyau de ce morphisme contient l'id\'{e}al $\textbf{W}_\ell\left(\frak{m}_{\mathcal{O}\lbrack \mu_{q-1} \rbrack}\right)$. On en d\'{e}duit par passage au quotient un morphisme du groupe $\textbf{W}_\ell(\mathcal{O}\lbrack \mu_{q-1}\rbrack)/\textbf{W}_\ell\left(\frak{m}_{\mathcal{O}\lbrack \mu_{q-1} \rbrack}\right)$ dans $\C_p^\times$. Or l'extension $K/\Q_p$ \'{e}tant totalement ramifi\'{e}e, il en est de m{\^e}me de $K\lbrack \mu_{q-1} \rbrack/\Q_p\lbrack \mu_{q-1} \rbrack$. Donc le corps r\'{e}siduel de $K\lbrack \mu_{q-1} \rbrack$ qui, par d\'{e}finition, est le quotient $\mathcal{O}\lbrack \mu_{q-1}\rbrack/\frak{m}_{\mathcal{O}\lbrack \mu_{q-1} \rbrack}$, est le m{\^e}me que le corps r\'{e}siduel de $\Q_p\lbrack \mu_{q-1} \rbrack$, ce qui permet de prolonger le morphisme de r\'{e}duction de $\Z_p\lbrack \mu_{q-1} \rbrack$ {\`a} $\F_q$ en un unique morphisme $\mathcal{O}\lbrack \mu_{q-1} \rbrack \to \F_q$. 
Il en r\'{e}sulte que le groupe $\textbf{W}_\ell(\mathcal{O}\lbrack \mu_{q-1}\rbrack)/\textbf{W}_\ell\left(\frak{m}_{\mathcal{O}\lbrack \mu_{q-1} \rbrack}\right)$ est isomorphe {\`a} $\textbf{W}_\ell(\F_q)$, ce qui prouve l'existence et l'unicit\'{e} du morphisme $\psi_{\ell, s, t}$.

 1. On a  montr\'{e} la surjectivit\'{e} du morphisme du groupe $\textbf{W}_\ell(\F_q)$ dans son groupe dual, qui, {\`a} un vecteur $\textit{\textbf{a}}$,  associe le caract{\`e}re de $ \textbf{W}_\ell(\F_q)$ qui envoie  tout vecteur $\textit{\textbf{y}} $ sur  $\psi_{\ell, s, t}(\textit{\textbf{a}}\textit{\textbf{y}}) $. Or, comme  $\textbf{W}_\ell(\F_q) $ et son dual sont deux groupes finis de m{\^e}me ordre, il suffit de prouver que cette m\^eme application est injective, c'est-{\`a}-dire que son noyau se r\'{e}duit {\`a} $\textit{\textbf{0}}. $ 
 
 Par contraposition, il s'agit de montrer que, pour tout vecteur $\textit{\textbf{a}}\neq\textit{\textbf{0}}$ de $ \textbf{W}_\ell(\F_q)$, il existe un vecteur $\textit{\textbf{y}}\in  \textbf{W}_\ell(\F_q),$ tel que l'image $ \psi_{\ell, s, t}(\textit{\textbf{a}}\textit{\textbf{y}})$ est non triviale. 
 
On utilise le fait \cite[proposition 8, p. AC IX.16]{Bourbaki} que l'anneau $\textbf{W}(\F_q)$ est un anneau de valuation discr{\`e}te d'uniformisante $V(\textit{\textbf{1}})$.  Par cons\'{e}quent, il existe un entier $n < \ell$ tel que l'id\'{e}al non nul $\textit{\textbf{a}} \textbf{W}_\ell(\F_q)$ de l'anneau $ \textbf{W}_\ell(\F_q) \simeq \textbf{W}(\F_q)/V^\ell\textbf{W}(\F_q)$ soit engendr\'{e} par la classe $V^n(\textit{\textbf{1}}) + V^\ell\textbf{W}(\F_q)$. Il existe donc un vecteur $\textit{\textbf{y}}$ de longueur $\ell$ tel que $\textit{ \textbf{a} \textbf{y}} = V^n(\textit{\textbf{1}}) + V^\ell\textbf{W}(\F_q)$. 
  
   Par d\'{e}finition du morphisme $\psi_{\ell, s, t}$, puis en utilisant le lemme \ref{pulita V}, on a alors $$\psi_{\ell, s, t} (\textit{\textbf{a}}\textit{\textbf{y}})=\theta_{\ell-1,s}(V^n \textit{\textbf{1}})(t) = \theta_{\ell-n-1,s}(\textit{\textbf{1}})(t^{p^n}) \; .$$ 
   
Or, d'apr{\`e}s la proposition \ref{valeurthetals a}, on sait que $$\theta_{\ell-n-1}(\textit{\textbf{1}})(t^{p^n})\equiv 1+\pi_{\ell-n-1}\left(\sum_{j=0}^{s-1}t^{{p^{n+j}}} \right)  \pmod{ \pi_{\ell-n-1}^2\C_p^\circ},$$
donc $\theta_{\ell-n-1}(\textit{\textbf{1}})(t^{p^n})\neq1$, d\'{e}s que $\sum_{j=0}^{s-1}t^{{p^{n+j}}}$ n'est pas un \'{e}l\'{e}ment de $\pi_{\ell-n-1} \C_p^\circ$. Or, l'extension  $\Q_p(\mu_{q-1})/\Q_p$ est galoisienne, non ramifi\'{e}e, et son groupe de Galois est isomorphe au groupe de Galois de l'extension r\'{e}siduelle, l'isomorphisme \'{e}tant d\'{e}fini en associant {\`a} tout automorphisme $\sigma$ de $\Q_p(\mu_{q-1})$ l'automorphisme $\widetilde{x} \mapsto \widetilde{\sigma(x)}$ du corps r\'{e}siduel \cite[proposition 2.11, p. 27]{Iwasawa}. Par cons\'{e}quent, le groupe de Galois de l'extension $\Q_p(\mu_{q-1})/\Q_p$ est cyclique, d'ordre $s$, engendr\'{e} par la restriction {\`a} $\Q_p(\mu_{q-1})$ du Frobenius absolu $\varphi$. On en d\'{e}duit, en utilisant le lemme \ref{imageFrobracine}, que 
$${\rm Tr}_{\Q_p(\mu_{q-1})/\Q_p}(t) = \sum_{j=0}^{s-1} \varphi^j(t) = \sum_{j=0}^{s-1} t^{p^j} \; .$$
Comme $t^{p^s} = t$, la valeur de $t^{p^j}$ ne d\'{e}pend que de la classe de l'entier $j$ modulo $s$, et donc $\sum_{j=0}^{s-1} t^{p^{n+j}} = \sum_{j=0}^{s-1} t^{p^j} = {\rm Tr}_{\Q_p(\mu_{q-1})/\Q_p}(t)$.

On conclut que $\psi_{\ell, s, t} (\textit{\textbf{a}}\textit{\textbf{y}})=\theta_{\ell-n-1}(\textit{\textbf{1}})(t^{p^n})\neq1$ si ${\rm Tr}_{\Q_p(\mu_{q-1})/\Q_p}(t)$ n'est pas \'{e}l\'{e}ment de $\pi_{n-\ell-1} \C_p^\circ$, ce qui revient {\`a} dire que ${\rm Tr}_{\Q_p(\mu_{q-1})/\Q_p}(t) \not\in p \Z_p$. 

2. On sait que le groupe ab\'{e}lien fini $\textbf{W}_\ell(\F_q)$ est isomorphe {\`a} son dual. Or, dans le groupe $ \textbf{W}_\ell(\F_q)$, identifi\'{e} au quotient $\textbf{W}(\F_q)/V^\ell \textbf{W}(\F_q)$,  l'\'{e}l\'{e}ment $\textit{\textbf{1}} + V^\ell \textbf{W}(\F_q)$ est d'ordre exactement  $p^\ell$, il existe donc un caract{\`e}re $\psi : \textbf{W}_\ell(\F_q) \to \C_p^\times$ qui est exactement d'ordre $p^\ell$, ce qui n\'{e}cessite que $\psi(\textbf{W}_\ell(\F_q))$ est le groupe $\mu_{p^\ell}$ tout entier. Ainsi il existe $\textit{\textbf{y}} \in W_\ell(\F_q)$ tel que $\psi(\textit{\textbf{y}})$ est une racine primitive d'ordre $p^\ell$ de l'unit\'{e}. D'apr{\`e}s le point 1. qui vient d'\^{e}tre d\'{e}montr\'{e}, il existe un vecteur $\textit{\textbf{a}}$ de $\textbf{W}_\ell(\F_q)$ tel que $\psi(\textit{\textbf{z}}) = \psi_{\ell, s, t} (\textit{\textbf{a} \textbf{z}})$ pour tout vecteur $\textit{\textbf{z}}$ de $\textbf{W}_\ell(\F_q)$. Par cons\'{e}quent le groupe $\mu_{p^\ell}$ est contenu dans l'image de $\psi_{\ell, s, t}$. Comme tout \'{e}l\'{e}ment de $\textbf{W}_\ell(\F_q)$ est d'ordre divisant $p^\ell$, l'inclusion r\'{e}ciproque est imm\'{e}diate.  
 \end{dem}\findem \end{sloppypar} 
 \begin{rem}  Dans le cas o{\`u} $s=1$, A. Pulita \cite[Theorem 2.7, p. 521]{Pulita} a voulu caract\'{e}riser l'image de $\textit{\textbf{1}} + V^\ell \textbf{W}(\F_p)$ par le morphisme $\psi_{\ell, 1, t}$, o{\`u} $t$ satisfait la relation $t^p = t$, {\`a} l'aide de l'unique racine $\xi_{\ell-1}$ primitive d'ordre $p^\ell$ de l'unit\'{e} telle que $ \vert t \pi_{\ell-1} - (\xi_{\ell-1} - 1) \vert <   \vert \pi_{\ell-1} \vert$. Malheureusement, cette racine $\xi_{\ell-1}$ n'est pas unique, sauf dans le cas $\ell=1$ d\'{e}j{\`a} trait\'{e} par Dwork \cite{Dwork2}.  Montrons en effet que le nombre d'\'{e}l\'{e}ments de l'ensemble $E_{t, \ell}$ des racines de l'unit\'{e} 
 $\xi$ d'ordre divisant $p^{\ell}$ qui v\'{e}rifient 
  $\vert (1+t\pi_{\ell-1} )-\xi \vert<\vert \pi_{\ell-1}\vert$ est \'{e}gal {\`a} $p^{\ell-1}$, donc n'est \'{e}gal {\`a} $1$ que pour  $\ell=1$. On se base sur la propri\'{e}t\'{e} que, si $\zeta \in \C_p$ est une racine primitive de l'unit\'{e} d'ordre $p^r$, o{\`u} $r \ge 1$ est un entier naturel, alors $\vert \zeta-1 \vert = \vert \pi_{r-1} \vert$. On en d\'{e}duit que, si $\zeta$ et $\zeta'$ sont deux \'{e}l\'{e}ments de $E_{t, \ell}$, alors $\zeta'/\zeta$ est d'ordre divisant $p^{\ell -1}$. Donc on a au plus $p^{\ell-1}$ \'{e}l\'{e}ments dans $E_{t, \ell}$. D'autre part, le corps r\'{e}siduel de $\Q_p(\mu_{p^\ell})$ \'{e}tant isomorphe {\`a} $\F_p$, on sait que, pour tout \'{e}l\'{e}ment $\zeta \in \mu_{p^\ell}$, il existe un \'{e}l\'{e}ment $t \in \Z_p$ tel que $t^p = t$ et $\left\vert \frac{\zeta-1}{\pi_{\ell-1}} - t \right\vert < 1$, c'est-{\`a}-dire $\vert (1+t \pi_{\ell-1}) - \zeta \vert < \vert \pi_{\ell-1} \vert$. Par cons\'{e}quent l'ensemble $\mu_{p^\ell}$ de cardinal $p^\ell$ est r\'{e}union des $p$ parties $E_{t, \ell}$ ayant chacune au plus $p^{\ell-1}$ \'{e}l\'{e}ments, donc chacune de ces parties doit avoir exactement $p^{\ell-1}$ \'{e}l\'{e}ments. 
  
\end {rem}

\subsection{Formule de transitivit\'{e}} 

Nous allons maintenant \'{e}noncer et d\'{e}montrer une relation entre les caract{\`e}res additifs $\psi_{\ell, sr, t}$ de $\textbf{W}_\ell(\F_{q^r})$ et les caract{\`e}res $\psi_{\ell, s, t}$ de $\textbf{W}_\ell(\F_q)$. Pour ce faire, nous aurons {\`a} utiliser l'application trace 
$$ {\rm Tr}_{\textbf{W}_\ell(\F_{q^r})/\textbf{W}_\ell(\F_q)}  : \textbf{W}_\ell(\F_{q^r}) \to \textbf{W}_\ell(\F_q) $$ 
d\'{e}finie comme suit. On note $\widetilde{\varphi}_r$ l'automorphisme de Frobenius qui {\`a} tout \'{e}l\'{e}ment $y$ de $\F_{q^r}$ associe $y^p$ ; alors on pose 
\begin{equation} \label{definition_trace} 
\forall \textit{\textbf{y}} \in \textbf{W}_\ell(\F_{q^r}) \qquad {\rm Tr}_{\textbf{W}_\ell(\F_{q^r})/\textbf{W}_\ell(\F_q)} \left( \textit{\textbf{y}}\right) = \sum_{i=0}^{r-1} \textbf{W}_\ell({\widetilde{\varphi}_r}^{is}) \left( \textit{\textbf{y}}\right) \; . 
\end{equation}
\begin{rem} 
L'automorphisme $\widetilde{\varphi}_r$ est induit sur le corps $\F_{q^r}$ par l'automorphisme du corps r\'{e}siduel $\widetilde{\C_p}
$ que nous avons not\'{e} $\widetilde{\varphi}$ dans la d\'{e}monstration du lemme \ref{Frobeniusideal}. Ceci explique la pr\'{e}sence de l'indice $r$ dans la notation. 
\end{rem} 

Nous aurons \'{e}galement besoin de l'application {\it rel{\`e}vement } ${\rm Te}$ de $\textbf{W}_\ell(\F_{q^r})$ dans $\textbf{W}(\Z_p\lbrack \mu_{q^r-1} \rbrack ) $ associant au vecteur de Witt $\textit{\textbf{y}}= (y_j)_{0 \le j < \ell}$ de longueur $\ell$ sur $\F_{q^r}$ le vecteur de Witt de longueur infinie sur $\Z_p\lbrack \mu_{q^r-1} \rbrack$ d\'{e}fini par 
\begin{equation} \label{relevement} 
{\rm Te}(\textit{\textbf{y}}) = \sum_{j=0}^{\ell-1} V^j(\boldsymbol{\tau}(\Te(y_j))) \; .
\end{equation}
Ainsi ${\rm Te}(\textit{\textbf{y}})$ n'est autre que le vecteur de Witt dont les $\ell$ premi{\`e}res composantes sont les images par le caract{\`e}re de Teichm\"uller des composantes du vecteur $\textit{\textbf{y}}$, les autres composantes \'{e}tant nulles. 

\begin{lem} \label{additif_relevement}
Si $\textit{\textbf{y}}$ et $\textit{\textbf{z}}$ sont deux vecteurs de Witt de longueur $\ell$ sur $\F_{q^r}$, alors on a la congruence 
$$ {\rm Te}\left(\textit{\textbf{y}} + \textit{\textbf{z}} \right) \equiv {\rm Te}(\textit{\textbf{y}}) + {\rm Te} (\textit{\textbf{z}}) \pmod{\textbf{W}(p\Z_p\lbrack \mu_{q^r-1} \rbrack) + V^\ell \textbf{W}(\Z_p\lbrack \mu_{q^r-1} \rbrack )} \; .$$
\end{lem} 
\begin{dem} Soit $\rho : \Z_p\lbrack \mu_{q^r-1} \rbrack \to \F_{q^r}$ un morphisme de r\'{e}duction. Puisque l'extension $\Q_p(\mu_{q^r-1})/\Q_p$ est non ramifi\'{e}e, le noyau du morphisme $\rho$ est $p \Z_p\lbrack \mu_{q^r-1} \rbrack$. En vertu de la proposition \ref{Wkermorphisme}, la congruence {\`a} montrer est \'{e}quivalente {\`a} l'identit\'{e} 
$$ \textbf{W}_\ell(\rho) \left( {\rm Te}\left(\textit{\textbf{y}} + \textit{\textbf{z}} \right) \right) = \textbf{W}_\ell(\rho)\left( {\rm Te}(\textit{\textbf{y}}) + {\rm Te} (\textit{\textbf{z}})\right) \; .$$ Or il est facile de voir, en utilisant la d\'{e}finition du caract{\`e}re de Teichm\"uller, que $\textbf{W}_\ell(\rho) \circ {\rm Te}$ est l'application identique de $\textbf{W}_\ell(\F_{q^r})$. Par cons\'{e}quent $\textbf{W}_\ell(\rho) \left( {\rm Te}\left(\textit{\textbf{y}} + \textit{\textbf{z}} \right) \right) = \textit{\textbf{y}} + \textit{\textbf{z}}$. Et, puisque $\textbf{W}_\ell(\rho)$ est un morphisme d'anneaux, on a aussi $\textbf{W}_\ell(\rho)\left( {\rm Te}(\textit{\textbf{y}}) + {\rm Te} (\textit{\textbf{z}})\right) = \textbf{W}_\ell(\rho)\left( {\rm Te}(\textit{\textbf{y}}) \right) + W_\ell(\rho) \left( {\rm Te} (\textit{\textbf{z}})\right) = \textit{\textbf{y}} + \textit{\textbf{z}}$. \end{dem} \findem


\begin{proposition} \label{transitivepsi} Soit $t \in \Z_p\lbrack \mu_{q-1} \rbrack$ tel que $t^q = t$. On a la formule de transitivit\'{e}
$$ \psi_{\ell, sr, t} = \psi_{\ell, s, t} \circ {\rm Tr}_{\textbf{W}_\ell(\F_{q^r})/\textbf{W}_\ell(\F_q)} \; .$$
\end{proposition} 
\begin{dem} Soit un vecteur de Witt $\textit{\textbf{y}} =(y_j)_{0 \le j < \ell}$ de longueur $\ell$ sur $\F_{q^r}$.  
Par d\'{e}finition,   $$ \psi_{\ell, sr, t} \left(\textit{\textbf{y}}\right)=\theta_{\ell-1, sr}\left( {\rm Te}(\textit{\textbf{y}})\right)(t) \; .$$
D'apr{\`e}s le lemme \ref{thetasr_fonction_thetas}, on en d\'{e}duit que $$\psi_{\ell, sr, t} \left(\textit{\textbf{y}}\right)=\prod_{i=0}^{r-1}\theta_{\ell-1,s}( {\rm Te}(\textit{\textbf{y}})^{\varphi^{is}})(t^{p^{si}}) \; .$$  
Or on sait que $t^q = t^{p^s} = t$ ; par cons\'{e}quent : 
$$\psi_{\ell, sr, t} \left(\textit{\textbf{y}}\right)=\prod_{i=0}^{r-1}\theta_{\ell-1,s}( {\rm Te}(\textit{\textbf{y}})^{\varphi^{is}})(t) \; .$$  
On utilise alors l'identit\'{e} 
$$ {\rm Te}(\textit{\textbf{y}})^\varphi = {\rm Te}(\textbf{W}_\ell(\widetilde{\varphi}_r)( \textit{\textbf{y}} ) ) \; $$ 
qui se d\'{e}duit du lemme \ref{imageFrobracine}. Donc 
$$\psi_{\ell, sr, t} \left(\textit{\textbf{y}}\right)=\prod_{i=0}^{r-1}\theta_{\ell-1,s}( {\rm Te}(\textbf{W}_\ell({\widetilde{\varphi}_r}^{is})(\textit{\textbf{y}}))(t) \; .$$   
Puisque $\theta_{\ell-1,s}$ est un morphisme, on en d\'{e}duit 
$$\psi_{\ell, sr, t} \left(\textit{\textbf{y}}\right)=\theta_{\ell-1,s} \left( \sum_{i=0}^{r-1} {\rm Te}(\textbf{W}_\ell({\widetilde{\varphi}_r}^{is})(\textit{\textbf{y}}) )\right)(t) \; .$$
Or par l'\'{e}quation \eqref{definition_trace} et par le lemme \ref{additif_relevement}, il existe un vecteur $\textit{\textbf{u}}$ \'{e}l\'{e}ment de l'id\'{e}al $\textbf{W}(p \Z_p\lbrack \mu_{q^r-1} \rbrack) + V^\ell \textbf{W}(\Z_p\lbrack \mu_{q^r-1} \rbrack)$ tel que 
$$ \sum_{i=0}^{r-1} {\rm Te}(\textbf{W}_\ell({\widetilde{\varphi}_r}^{is})(\textit{\textbf{y}}) = {\rm Te} \left( {\rm Tr}_{\textbf{W}_\ell(\F_{q^r})/\textbf{W}_\ell(\F_q)}(\textit{\textbf{y}}) \right) + \textit{\textbf{u}} \; .$$
Comme la valeur de la s\'{e}rie $\theta_{\ell-1,s}(\textbf{\textit{u}})$ en $t$ est \'{e}gale {\`a} 1 en vertu des lemmes \ref{noyautheta} et \ref{2.1}, on en d\'{e}duit
$$ \psi_{\ell, sr, t} \left(\textit{\textbf{y}}\right) = \theta_{\ell-1,s}\left( {\rm Te} \left( {\rm Tr}_{\textbf{W}_\ell(\F_{q^r})/\textbf{W}_\ell(\F_q)}(\textit{\textbf{y}}\right) )\right)(t) \; ,$$
\'{e}quation o{\`u} on reconna\^it dans le membre de droite $\psi_{\ell, s, t}\left( {\rm Tr}_{\textbf{W}_\ell(\F_{q^r})/\textbf{W}_\ell(\F_q)}(\textit{\textbf{y}}) \right)$. 
 \end{dem}\findem
 
 \subsection{Fonctions scindantes}

On va maintenant chercher {\`a} exprimer analytiquement le caract{\`e}re additif $\psi_{\ell,s, t}$ du groupe $\textbf{W}_\ell(\F_q)$ qui a \'{e}t\'{e} d\'{e}fini par la proposition \ref{2.4}. Remarquons que, par d\'{e}finition de $\psi_{\ell,s,t}$, on doit avoir $t^q = t$. 

Soit un vecteur de Witt $\textit{\textbf{y}}=(y_j)_{0 \le j < \ell}$ de longueur $\ell$ sur $\mathbb{F}_{q}$.  
On a vu pr\'{e}c\'{e}demment comment lui associer son rel{\`e}vement  ${\rm Te}({\textit{\textbf{y}}})$, vecteur de Witt de longueur infinie sur $\Z_p\lbrack \mu_{q-1} \rbrack$, d\'{e}fini par $${\rm  Te}({\textit{\textbf{y}}}) =  \sum_{j=0}^{\ell-1} V^j(\boldsymbol{\tau}(\Te(y_j)) \; .$$  Par d\'{e}finition du caract{\`e}re $\psi_{\ell, s, t}$, on a 
$$ \psi_{\ell, s, t}(\textit{\textbf{y}}) = \theta_{\ell-1,s}({\rm Te}(\textit{\textbf{y}}))(t) \; .$$ 
Or les lemmes \ref{pulita V} et \ref{thetataums} conduisent {\`a} l'\'{e}galit\'{e}  
$$ \theta_{\ell-1,s}({\rm Te}(\textit{\textbf{y}})) = \prod_{j=0}^{\ell-1} \theta_{\ell-j-1,s}(\textit{\textbf{1}}) \circ (\Te(y_j) x^{p^j} ) \; .$$
Par cons\'{e}quent
\begin{equation} \label{explicitecaractere} \psi_{\ell, s, t}(\textit{\textbf{y}}) = \prod_{j=0}^{\ell-1} \theta_{\ell-j-1,s}(\textit{\textbf{1}})(t^{p^j} \Te(y_j)) \; . \end{equation}

Introduisons la "fonction" $\Omega_{\ell, s, t}$ qui est la s\'{e}rie enti{\`e}re en $\ell$ ind\'{e}termin\'{e}es $x_0, x_1, \ldots, x_{\ell-1}$ {\`a} coefficients dans l'anneau $\mathcal{O}\lbrack \mu_{q-1} \rbrack$ d\'{e}finie par 
\begin{equation}
\label{definitionfonctionscindante} 
\Omega_{\ell, s, t}(x_0, x_1, \ldots, x_{\ell-1}) = \prod_{j=0}^{\ell-1} \theta_{\ell-j-1,s}(\textit{\textbf{1}})(t^{p^j} x_j) 
\end{equation} 

Avec cette notation $\Omega_{\ell,s,t}$, la formule \eqref{explicitecaractere} s'\'{e}crit comme suit. 

\begin{proposition} \label{expression_analytique_caractere_additif} 
Soit $t \in \Z_p \lbrack \mu_{q-1} \rbrack$ tel que $t^q = t$. Alors, pour tout vecteur $\textit{\textbf{y}} = (y_i)_{0 \le i < \ell}$ de longueur $\ell$ sur $\F_q$, on a 
$$ \psi_{\ell, s, t} ( \textit{\textbf{y}}) = \Omega_{\ell, s, t}(\Te(y_0), \ldots, \Te(y_{\ell-1}) )\; . $$

\end{proposition} 

\begin{proposition} \label{factorisation_fonction_scindante} 
Soit $t \in \Z_p \lbrack \mu_{q-1} \rbrack$ tel que $t^q = t$. Si $r \ge 1$ est un entier naturel, alors on a  l'identit\'{e}
$$ \Omega_{\ell, sr, t}(x_0, \ldots, x_{\ell-1}) = \prod_{i=0}^{r-1} \Omega_{\ell, s, t}(x_0^{q^i}, \ldots, x_{\ell-1}^{q^i}) \; . $$
\end{proposition}
 \begin{dem}  
 On a par l'\'{e}quation  \eqref{definitionfonctionscindante} $$ \Omega_{\ell, sr, t}(x_0, \ldots, x_{\ell-1})=\prod_{j=0}^{\ell-1} \theta_{\ell-j-1,sr}(\textit{\textbf{1}})(t^{p^j} x_j). $$
  D'apr{\`e}s le lemme \ref{thetasr_fonction_thetas} $$ \Omega_{\ell, sr, t}(x_0, \ldots, x_{\ell-1})=\prod_{j=0}^{\ell-1}\prod_{i=0}^{r-1}\left(\theta_{\ell-j-1,s}(\textit{\textbf{1}})\circ x^{p^{si}}\right)(t^{p^j} x_j) \; ; $$
donc on obtient $$ \Omega_{\ell, sr, t}(x_0, \ldots, x_{\ell-1})=\prod_{i=0}^{r-1}\prod_{j=0}^{\ell-1}\theta_{\ell-j-1,s}(\textit{\textbf{1}}){((t^{p^j}) x_j)}^{p^{si}} \; , $$ et par la suite, on a $$\Omega_{\ell, sr, t}(x_0, \ldots, x_{\ell-1}) = \prod_{i=0}^{r-1} \prod_{j=0}^{\ell-1} \theta_{\ell-j-1,s}( \textit{\textbf{1}}){((t^{p^j})^{p^{si}} (x_j)^{p^{si}}} \; , $$
 d 'o{\`u} $$\Omega_{\ell, sr, t}(x_0, \ldots, x_{\ell-1})=\prod_{i=0}^{r-1} \prod_{j=0}^{\ell-1}\theta_{\ell-j-1,s}( \textit{\textbf{1}}){(t^{p^j} x_j^{q^{i}}}) \; , $$
 ce qui donne le r\'{e}sultat cherch\'{e} en utilisant l'\'{e}quation \eqref{definitionfonctionscindante} exprimant la s\'{e}rie $\Omega_{\ell, s, t}$.  
\end{dem}\findem

Nous rappelons la d\'{e}finition des fonctions scindantes donn\'{e}e par \cite{Blache} en l'\'{e}tendant du corps premier $\F_p$ au corps fini $\F_q$. 
\begin{defi} \label{definition_fonction_scindante} Soit $\ell \ge 1$ un entier naturel, et $q=p^s$ une puissance de $p$. Une {\it fonction $s$-scindante de niveau $\ell$}, est une s\'{e}rie enti{\`e}re formelle $\Omega(x_1, \ldots, x_{\ell})$, {\`a} coefficients dans $ \mathcal{O}$, en les $\ell$ ind\'{e}termin\'{e}es $x_1, \ldots, x_{\ell}$, qui converge sur un polydisque de la forme $ D(0, r_1) \times \cdots \times D(0, r_\ell)$,  dont les rayons  
$r_1,  \ldots , r_\ell$ sont tous sup\'{e}rieurs {\`a} 1, et v\'{e}rifiant de plus les deux conditions suivantes. \\
\indent 1. La  fonction $\psi$ de $\textbf{W}_\ell(\F_q)$ dans $\C_p$ qui au vecteur de Witt $\textit{\textbf{y}} = (y_j)_{0 \le j < \ell}$ associe $\Omega(\Te(y_0), \Te(y_1) , \ldots, \Te(y_{\ell-1}))$ est un caract{\`e}re additif d'ordre $p^\ell$ de $\textbf{W}_\ell(\F_q)$ dans $\C_p^*$.\\ 
\indent 2. Pour chaque entier $r \ge 1$, le caract{\`e}re additif de $\textbf{W}_\ell(\F_{q^r})$  obtenu en composant $\psi$ avec l'application trace ${\rm Tr}_{\textbf{W}_\ell(\F_q^r)/\textbf{W}_\ell(\F_q)}$ s'exprime sous la forme 
\begin{equation} \label{splitting_function} 
\psi\left( {\rm Tr}_{\textbf{W}_\ell(\F_q^r)/\textbf{W}_\ell(\F_q)} (\textit{\textbf{y}}) \right) = \prod_{i =0}^{r-1} \Omega\left(\Te(y_0)^{q^i}, \ldots, \Te(y_{\ell-1})^{q^i} \right) \; \end{equation} 
pour tout vecteur de Witt $ \textit{\textbf{y}} = (y_j)_{0 \le j < \ell} \in \textbf{W}_\ell(\F_{q^r}).$
\end{defi} 
\begin{thm} Soit $\ell \ge 1, s \ge 1$ deux entiers naturels, et notons $q=p^s$. Soit $t \in \Z_p\lbrack \mu_{q-1} \rbrack$ tel que $t^{q}=t$ et ${\rm Tr}_{\Q_p(\mu_{q-1})/\Q_p}(t) \not\in p\Z_p$. Alors la s\'{e}rie $\Omega_{\ell, s, t}$ est une fonction $s$-scindante de niveau $\ell$. 
\end{thm} 
\begin{dem} Il est \'{e}vident que la s\'{e}rie  $\Omega_{\ell, s, t}(x_0, x_1, \ldots, x_{\ell-1})$,   introduite par l'\'{e}quation \eqref{definitionfonctionscindante}  comme \'{e}tant le produit des s\'{e}ries $\theta_{\ell-j-1,s}(\textit{\textbf{1}})(t^{p^j} x_j) $,
est {\`a} coefficients dans l'anneau $ \mathcal{O}$, puisque chaque facteur du produit est {\`a} coefficients dans $ \mathcal{O}$. Puisqu'on sait par  la proposition \ref{surconvergencethetams} que chacune de ces s\'{e}ries $ \theta_{\ell-j-1,s}(\textit{\textbf{1}}) $ a un rayon de convergence strictement sup\'{e}rieur {\`a} 1, il en est de m{\^e}me des m{\^e}mes s\'{e}ries appliqu\'{e}es {\`a} $t^{p^j} x_j$, puisque $t^{p^j}$ est de module 1. Par cons\'{e}quent la s\'{e}rie  $\Omega_{\ell, s, t}(x_0, x_1, \ldots, x_{\ell-1})$ est convergente dans un produit de disques de rayons strictement sup\'{e}rieurs {\`a} 1.
 D'autre part, d'apr{\`e}s les propositions \ref{2.4} et   \ref{expression_analytique_caractere_additif}, la serie $\Omega_{\ell, s, t}$ v\'{e}rifie la condition 1, et, par le biais des propositions \ref{transitivepsi}, \ref{expression_analytique_caractere_additif} et \ref{factorisation_fonction_scindante}, la m\^{e}me s\'{e}rie  $\Omega_{\ell, s, t}$ v\'{e}rifie la condition 2 de la d\'{e}finition \ref{definition_fonction_scindante} ; ce qui ach{\`e}ve la d\'{e}monstration. 
\end{dem}\findem

\subsection{Expression analytique d'un caract{\`e}re multiplicatif de $W_{2}(\mathbb{F}_{q})$} Dans le cas $\ell=1$, on a vu que le groupe des caract{\`e}res du groupe $\F_q^*$ \'{e}tait le groupe cyclique engendr\'{e} par le caract{\`e}re $\Te$ de Teichm\"uller. Pour simplifier, on va limiter ici l'\'{e}tude des caract{\`e}res multiplicatifs au cas $\ell=2$. On note $W_2(\F_q)^\star$ le groupe des \'{e}l\'{e}ments inversibles de l'anneau des vecteurs de Witt sur $\F_q$ de longueur 2. Un vecteur de Witt $\textit{\textbf{z}} = (z_0, z_1)$ est \'{e}l\'{e}ment de $\textbf{W}_\ell(\F_q)^\star,$ si et seulement si $z_0\neq 0.$ 
Soit $(z_0, z_1)\in W_{2}\left(\mathbb{F}_{q}\right)^{*}.$
 Exprimons le vecteur $\textit{\textbf{z}}$ sous la forme $(z_0, z_1) = (z_0, 0) (1, y)$, autrement dit $y$ v\'{e}rifie l'identit\'{e} $(z_0, z_1) = (z_0, z_0^py)$, i.e. $y=  \frac{z_1}{z_0^p}=z_1 z_0^{p(q-2)}$. Pour tout caract{\`e}re multiplicatif $\chi$, on a alors 
 $$ \chi(\textit{\textbf{z}}) = \chi(z_0,0) \chi(1,y) \; ;$$
 or $z_0 \mapsto \chi(z_0,0)$ est un caract{\`e}re multiplicatif de $\F_q$ comme compos\'{e} des morphismes multiplicatifs $\chi$ et $z_0 \mapsto (z_0, 0)$ ; par cons\'{e}quent, il existe un entier $m \in \lbrack 0 \ldots q-2 \rbrack$ tel qu'on ait 
 $$ \forall z_0 \in \F_q, \qquad \chi(z_0, 0) = \Te(z_0)^m.$$ 
 Par ailleurs, $y \mapsto \chi(1,y)$ est un caract{\`e}re additif de $\F_q$ car le produit dans $\textbf{W}_2(\F_q)$ des vecteurs $(1,y)$ et $(1,y')$ est \'{e}gal {\`a} $(1, y+y')$. Par cons\'{e}quent, pour tout $t \in \Z_p\lbrack \mu_{q-1} \rbrack$ tel que $t^q = t$ et ${\rm Tr}_{\Q_p(\mu_{q-1})/\Q_p}(t) \not\in p\Z_p$, il existe un \'{e}l\'{e}ment $b \in \F_q$ tel que 
 $$ \forall y \in \F_q, \qquad \chi(1,y) = \psi_{1,s,t}(by) = \Omega_{1,s,t}(\Te(b) \Te(y)) \; .$$
 
 \begin{proposition} \label{expression_chi} 
 Soit $t \in \Z_p\lbrack \mu_{q-1} \rbrack$ tel que $t^q = t$ et ${\rm Tr}_{\Q_p(\mu_{q-1})/\Q_p}(t) \not\in p\Z_p$, et $\chi$ un caract{\`e}re multiplicatif de $W_2(\F_q)$. Alors il existe un entier $m \in \lbrack 0 .. q-2 \rbrack$, et un \'{e}l\'{e}ment $b$ de $\F_q$ tel que 
 $$ \forall \textit{\textbf{z}} = (z_0, z_1) \in \textbf{W}_2(\F_q)^\star, \qquad \chi(\textit{\textbf{z}}) = \Te(z_0)^m \Omega_{1,s,t}(\Te(b)\Te(z_1) \Te(z_0)^{p(q-2)}) \; .$$
 \end{proposition} 
\section{La formule de trace }

G\'{e}n\'{e}ralisons la formule de trace obtenue par Dwork \cite{Dwork} pour exprimer les sommes de Gauss classiques, au cas de sommes d\'{e}finies sur des anneaux de vecteurs de Witt de longueur finie sur un corps fini $\F_q$, o{\`u} $q=p^s$ est une puissance du nombre premier $p$. Commen{\c c}ons par traiter le cas particulier des vecteurs de longueur $\ell=2$. 

\subsection{Expression analytique de la somme de Gauss} 
 \begin{defi} \label{definition_somme_de_Gauss} Soit $\psi$ un caract{\`e}re additif, et $\chi$ un caract{\`e}re multiplicatif, de $\textbf{W}_2(\F_q)$. On leur associe la {\it somme de Gauss} 
  $$g(\psi, \chi) = -\sum_{ \mathit{\mathbf{y}}\in \mathbb{W}_2(\F_q)^\star}\psi( \mathit{\mathbf{y}})\chi( \mathit{\mathbf{y}}). $$
\end{defi} 
Gr\^{a}ce {\`a} la proposition \ref{2.4}, et si $t \in \mathbb{Z}_p\lbrack \mu_{q-1} \rbrack$ v\'{e}rifie $t^{p^{s}} = t$ et ${\rm Tr}_{\Q_p(\mu_{q-1})/\Q_p}(t) \not\in p \Z_p$, il existe un vecteur de Witt $ \mathit{\mathbf{a}}$ de longueur $2 $ sur $\F_q$  tel que $\psi$ s'\'{e}crive comme compos\'{e}e du caract{\`e}re $\psi_{2,s,t}$ avec l'homoth\'{e}tie de rapport $\textit{\textbf{a}}$ de l'anneau $\textbf{W}_2(\F_q)$, de sorte que la somme de Gauss s'exprime comme suit $$g(\psi, \chi) = -\sum_{ \mathit{\mathbf{y}}\in \mathbb{W}_2(\F_q)^\star}\psi( \mathit{\mathbf{y}})\chi( \mathit{\mathbf{y}})=-\sum_{ \mathit{\mathbf{y}}\in \mathbb{W}_2(\F_q)^\star}\psi_{2, s, t}(  \mathit{\mathbf{a}}\mathit{\mathbf{y}})\chi( \mathit{\mathbf{y}}) \; .$$
On se restreint au cas o{\`u} le caract{\`e}re additif $\psi$ est d'ordre exactement $p^2$. Dans ce cas $\mathit{\mathbf{a}}\in \mathbb{W}_2(\F_q)^\star$, et on peut r\'{e}\'{e}crire la somme de Gauss comme suit : $$g(\psi, \chi)=-\sum_{ \mathit{\mathbf{z}}\in \mathbb{W}_2(\F_q)^\star}\psi_{2, s, t}(\mathit{\mathbf{z}})\chi(\mathit{\mathbf{a^{-1}}} \mathit{\mathbf{z}})\; \qquad ,$$ ce qui est \'{e}quivaut {\`a} dire que$$g(\psi, \chi) = -\chi(\mathit{\mathbf{a^{-1}}})\sum_{ \mathit{\mathbf{z}}\in \mathbb{W}_2(\F_q)^\star}\psi_{2, s, t}(\mathit{\mathbf{z}}) \chi(\mathit{\mathbf{z}})=\chi(\mathit{\mathbf{a^{-1}}})g(\psi_{2, s, t}, \chi) \; \qquad .$$ 
 Ainsi, {\`a} condition que le caract{\`e}re additif $\psi$ soit d'ordre maximal, l'\'{e}tude de la somme de caract{\`e}res $g(\psi, \chi)$ se ram{\`e}ne {\`a} celle de $g(\psi_{2, s, t}, \chi)$.  Donnons maintenant une expression analytique de la somme de caract{\`e}res $g(\psi_{2, s, t}, \chi).$ En \'{e}crivant $ \mathit{\mathbf{z}} = (z_0, z_1)$, on a :
\begin{equation} \label{sommede Gauss1} g(\psi_{2, s, t}, \chi)=-\sum_{ z_0\in \F_q^\star}\sum_{ z_1 \in \F_q}\psi_{2, s, t}(z_0, z_1 ) \chi (z_0, z_1)\; . \end{equation}
D'apr{\`e}s la proposition \ref{expression_chi}, il existe un entier naturel $m \in \lbrack 0 .. q-2 \rbrack$ et un \'{e}l\'{e}ment $b \in \F_q$ tels que, pour tout $(z_0, y) \in \F_q^\star \times \F_q$, on ait $$\chi(z_0, z_1) = \Te(z_0)^m \Omega_{1,s,t}(\Te(b) \Te(z_1) \Te(z_0)^{p(q-2)}) \; .$$ D'autre part, par la proposition \ref{expression_analytique_caractere_additif},  on sait que $\psi_{2,s,t}(z_0, z_1)$ se met sous la forme $\Omega_{2,s,t}(\Te(z_0), \Te(z_1))$. Substituant ces expressions dans \eqref{sommede Gauss1}, on obtient, en posant 
\begin{equation} \label{noyau} 
\widehat{H}(x_0, x_1) = - \Omega_{2,s,t}(x_0, x_1) x_0^m \Omega_{1,s,t}(\Te(b) x_1 x_0^{p(q-2)}), 
\end{equation}   
l'expression de la somme de Gauss sous la forme 
\begin{equation} 
\label{expression_analytique_somme_de _Gauss}  g(\psi_{2, s, t}, \chi) = \sum_{ z_0\in \F_q^\star}\sum_{ z_1 \in \F_q} \widehat{H}(\Te(z_0), \Te(z_1)). \end{equation} 


\subsection{L'espace des fonctions analytiques sur un polydisque ferm\'{e}} 
Pour $\beta > 0$ un nombre r\'{e}el, notons $E_\beta$ l'espace vectoriel sur le corps $K_1\lbrack \mu_{q-1} \rbrack$ dont les \'{e}l\'{e}ments sont les s\'{e}ries enti{\`e}res {\`a deux variables \begin{equation} \label{serie2variables} G(x_0, x_1) = \sum_{(n_0, n_1) \in \mathbb{N}^2} a_{n_0, n_1} x_0^{n_0} x_1^{n_1} \end{equation} {\`a coefficients dans $K_1\lbrack \mu_{q-1} \rbrack$, telles que la famille $ \left( \vert a_{n_0, n_1} \vert \beta^{n_0 + n_1} \right)_{(n_0, n_1) \in \mathbb{N}^2}$ tende vers z\'{e}ro selon le filtre des compl\'{e}mentaires de parties finies de $\mathbb{N}^2$. Pour une s\'{e}rie $G(x_0, x_1)$ \'{e}l\'{e}ment de $E_\beta$ \'{e}crite sous la forme \eqref{serie2variables}, on pose 
$$ \Vert G \Vert_\beta = \sup \lbrace \vert a_{n_0, n_1} \vert \beta^{n_0 + n_1} , (n_0, n_1) \in \mathbb{N}^2 \rbrace \; , $$ ce qui fait de $E_\beta$ un espace vectoriel norm\'{e}. 

On peut interpr\'{e}ter l'espace $E_\beta$ comme la somme directe de Banach \cite[p. 184]{Robert} de la famille d'espaces de Banach $\left( e_{n_0, n_1}(\beta) \right)_{(n_0, n_1) \in \mathbb{N}^2}$ d\'{e}finis comme suit : $e_{n_0, n_1}(\beta)$ est l'espace vectoriel de dimension un $K_1 \lbrack \mu_{q-1} \rbrack$ norm\'{e} en posant 
$$\forall a \in K_1 \lbrack \mu_{q-1} \rbrack, \qquad \Vert a \Vert_{n_0, n_1, \beta} = \vert a \vert \beta^{n_0 + n_1} \; .$$ 
Ceci entra\^ine que, pout tout r\'{e}el $\beta > 0$,  l'espace vectoriel norm\'{e} $E_\beta$ est un espace de Banach \cite[Theorem p. 184]{Robert}. 

Il est imm\'{e}diat que,  si $0 < \beta < \beta'$, alors $E_{\beta'}$ est un sous-espace de $E_\beta$. 

     \begin{defi} Soient $E$ et $F$ deux espaces de Banach, et $f : E\to F$ une application lin\'{e}aire. On dit que $f$ est compl{\`e}tement continue si elle est la limite d'une suite d'applications lin\'{e}aires continues de rang fini  \end{defi}
     \begin{proposition} Soit $\beta$ et $\beta'$ deux nombres r\'{e}els tels que $0 < \beta < \beta'$. L'application injection canonique (ou insertion) $i_{\beta', \beta} : E_{\beta'}\to E_{\beta}$ est compl{\`e}tement  continue et son image est dense. 
      \end {proposition}
      \begin{dem} Pour un \'{e}l\'{e}ment $G = \sum_{(n_0, n_1)} a_{n_0, n_1} x_0^{n_0} x_1^{n_1}$ de $E_\beta$ et pour un entier naturel $d$, notons $G_d$ le polyn\^{o}me 
     $$ G_d =  \sum_{n_0+n_1\le d}
 a_{n_0, n_1} x_0^{n_0}x_1^{n_1} \; .$$
 Pour tout $G = \sum
 a_{n_0, n_1} x_0^{n_0}x_1^{n_1}$ \'{e}l\'{e}ment de $E_{\beta}$, on a 
\begin{equation} \label{normereste} \Vert G - G_d \Vert_\beta = \sup_{n_0 + n_1 > d} \vert a_{n_0, n_1} \vert \beta^{n_0+n_1} \; . \end{equation}
 Puisque $G$ appartient {\`a} $E_\beta$, on en d\'{e}duit que la norme $\Vert G - G_d \Vert_\beta$ converge vers z\'{e}ro. Puisque les $G_d$ sont des polyn\^{o}mes, donc des \'{e}l\'{e}ments de $E_{\beta'}$, ceci suffit {\`a} montrer que l'image de $i_{\beta', \beta}$ est dense dans $E_\beta$. 
 
Pour $d \in \mathbb{N}$, consid\'{e}rons l' application $i_d : E_{\beta'}\to E_\beta $ qui {\`a} toute s\'{e}rie $G$  associe la s\'{e}rie tronqu\'{e}e $i_d(G) = G_d$. L'image de $i_d$ \'{e}tant contenue dans le sous-espace engendr\'{e} par les mon\^{o}mes de degr\'{e} total $\le d$, sous-espace qui est de dimension finie, on voit que l'application lin\'{e}aire $i_d$ est de rang fini pour tout $d \in \mathbb{N}$. De plus, on a pour tout $G = \sum
 a_{n_0, n_1} x_0^{n_0}x_1^{n_1}$ \'{e}l\'{e}ment de $E_{\beta'}$, 
 $$ \scriptstyle{\Vert G_d \Vert_\beta = \sup_{n_0 + n_1 \le d} \vert a_{n_0, n_1} \vert \beta^{n_0 + n_1} \le \sup_{n_0 + n_1 \le d} \vert a_{n_0, n_1} \vert \beta' \, ^{n_0 + n_1} \le \sup_{(n_0, n_1) \in \mathbb{N}^2} \vert a_{n_0, n_1} \vert \beta' \, ^{n_0 + n_1} = \Vert G \Vert_{\beta'} \; , }$$
 ce qui prouve que l'application lin\'{e}aire $i_d$ est continue de $E_{\beta'}$ dans $E_\beta$. 
De plus, l'\'{e}quation \eqref{normereste} entra\^ine la suite de relations 
 $$ \Vert G - G_d \Vert_\beta \le \Vert G \Vert_{\beta'} \sup_{n_0 + n_1 > d} \beta'^{-n_0 - n_1} \beta^{n_0 + n_1} = \Vert G \Vert_{\beta'} \left( \frac{\beta}{\beta'} \right)^{d+1},  $$     
 ce qui prouve que $\Vert i_{\beta', \beta} - i_d \Vert \le \left( \frac{\beta}{\beta'} \right)^{d+1}$, et donc, puisque $\frac{\beta}{\beta'} < 1$, que $\lim_{d \to + \infty} i_d = i_{\beta', \beta}$, c'est-{\`a}-dire que $i_{\beta', \beta}$ est compl{\`e}tement continue.  
     \end{dem}\findem

 
 \subsection{Op\'{e}rateur de Dwork} 
Rappelons que $q$ d\'{e}signe l'entier $p^s$. 
 \begin{notation} 
 Si $G = \sum_{(n_0, n_1) \in \mathbb{N}^{2}} a_{n_0, n_1} x_0^{n_0} x_1^{n_1}$ est une s\'{e}rie enti{\`e}re en deux variables {\`a} coefficients dans $\C_p$, on note ${\rm Dw}_q(G)$ la s\'{e}rie $ \sum_{(n_0, n_1) \in \mathbb{N}^{2}} a_{q n_0, q n_1} x_0^{n_0} x_1^{n_1} . $
 \end{notation} 
 \begin{lem} \label{opDwork} 
 Soit $\beta > 0$ un nombre r\'{e}el. Si la s\'{e}rie enti{\`e}re $G$ est \'{e}l\'{e}ment de l'espace de Banach $E_\beta$, alors la s\'{e}rie ${\rm Dw}_q(G)$ est \'{e}l\'{e}ment de $E_{\beta^{q}}$. De plus, on a l'in\'{e}galit\'{e} 
 $$ \Vert {\rm Dw}_q(G) \Vert_{\beta^{q}} \le \Vert G \Vert_\beta \; .$$
 \end{lem} 
 \begin{dem}   Soit $G(x_0, x_1)=\sum
 a_{n_0, n_1} x_0^{n_0}x_1^{n_1}$ une s\'{e}rie de $ E_{\beta}$. 
 Comme l'application de l'ensemble $\mathbb{N}^2$ dans lui-m\^{e}me envoyant $(n_0, n_1)$ sur $(q n_0, q n_1)$ est injective, l'image r\'{e}ciproque par cette application du compl\'{e}mentaire d'une partie finie est compl\'{e}mentaire d'une partie finie. Par cons\'{e}quent, la famille $ \left( \vert a_{qn_0, qn_1} \vert \beta^{q(n_0+n_1)} \right)_{(n_0, n_1) \in\mathbb{N}^{2}}$, obtenue par composition avec cette application de la famille $ \left( \vert a_{n_0, n_1} \vert \beta^{n_0+n_1} \right)_{(n_0, n_1) \in\mathbb{N}^{2}}$ qui tend vers z\'{e}ro, est elle aussi une famille tendant vers z\'{e}ro. Ceci montre que la s\'{e}rie ${\rm Dw}_q(G)$ est \'{e}l\'{e}ment de $E_{\beta^q}$. De plus, l'ensemble des nombres r\'{e}els de la forme $\vert a_{q n_0, q n_1} \vert \beta^{q(n_0 + n_1)} $ \'{e}tant une partie de l'ensemble des nombres r\'{e}els de la forme $\vert a_{n_0, n_1} \vert \beta^{n_0 + n_1}$, est major\'{e} par la borne sup\'{e}rieure de ce deuxi{\`e}me ensemble, not\'{e}e $\Vert G \Vert_\beta$. Ceci montre que la norme dans $E_{\beta^q}$ de ${\rm Dw}_q(G)$ est major\'{e}e par $\Vert G \Vert_\beta$. \end{dem} \findem
 \begin{defi} 
  Soit $\beta > 0$ un nombre r\'{e}el. On appelle op\'{e}rateur de Dwork l'application ${\rm Dw}_q$ de $E_\beta$ dans $E_{\beta^{q}}$.  
 \end{defi} 
 \begin{proposition}\label{operateurDwork} 
L'op\'{e}rateur de Dwork ${\rm Dw}_q :  E_{\beta} \to  E_{\beta^q}$   est continu pour tout $\beta > 0$. \end{proposition}
\begin{dem}   Par le lemme \ref{opDwork}, on sait que ${\rm Dw}_q$ est une application lipschitzienne, elle est donc continue.
\end{dem}\findem 

\subsection{Op\'{e}rateurs de multiplication} 
L'ensemble des s\'{e}ries enti{\`e}res formelles en deux variables {\`a} coefficients dans $\C_p$ est muni du produit de Cauchy. 
\begin{proposition} \label{Ebetastable} Soient $G = \sum a_{n_0, n_1}x_0^{n_0}x_1^{n_1}$ et $H = \sum b_{n_0, n_1}x_0^{n_0}x_1^{n_1}$ deux s\'{e}ries appartenant {\`a} l'espace $E_\beta$. Alors leur produit de Cauchy $G H$ est \'{e}l\'{e}ment de $E_\beta$, et on a
\begin{equation}
\label{normeproduit} 
\Vert G H \Vert_\beta \le \Vert G \Vert_\beta \Vert H \Vert_\beta
\end{equation}
\end{proposition} 
\begin{dem}  On sait que $G H = \sum c_{n_0, n_1} x_0^{n_0} x_1^{n_1}$, avec, pour tout $(n_0, n_1) \in \mathbb{N}^2$ : 
\begin{equation} \label{definition produit Cauchy} c_{n_0, n_1} = \sum_{k_0 = 0}^{n_0} \sum_{k_1 = 0}^{n_1} a_{k_0, k_1} b_{n_0 - k_0, n_1 - k_1} \; . \end{equation}                  
Il n'y rien {\`a} montrer si $G = H = 0$. Supposons donc qu'au moins l'une des deux s\'{e}ries $G$ et $H$ est non nulle, c'est-{\`a}-dire que $\max( \Vert G \Vert_\beta, \Vert H \Vert_\beta ) > 0 $. Pour montrer que le produit $G H$ est \'{e}l\'{e}ment de $E_\beta$, on choisit un re\'{e}l arbitraire $\epsilon > 0$ pour justifier l'existence d'une partie finie $C(\epsilon)$ de $\mathbb{N}^2$ telle que $\vert c_{n_0, n_1} \vert \beta^{n_0 + n_1} < \epsilon$ pour tout $(n_0, n_1) \in \mathbb{N}^2 \setminus C(\epsilon)$. 

Puisque $G$ appartient {\`a} $E_\beta$, il existe une partie finie $A(\epsilon)$ de $\mathbb{N}^{2}$ telle qu'on ait l'implication 
\begin{equation} \label{minorationa} 
\vert a_{k_0, k_1} \vert \beta^{k_0 + k_1} \ge \frac{\epsilon}{\max( \Vert G \Vert_\beta , \Vert H \Vert_\beta )} \Rightarrow (k_0, k_1) \in A(\epsilon) \; . 
\end{equation}
De m{\^e}me, puisque $H$ appartient {\`a} $E_\beta$, il existe une partie finie $B(\epsilon)$ de $\mathbb{N}^{2}$ telle qu'on ait l'implication 
\begin{equation} \label{minorationb} 
\vert b_{l_0, l_1} \vert \beta^{l_0 + l_1} \ge \frac{\epsilon}{\max( \Vert G \Vert_\beta , \Vert H \Vert_\beta )} \Rightarrow (l_0, l_1) \in B(\epsilon) \; . 
\end{equation}
Si $\vert c_{n_0, n_1} \vert \beta^{n_0 + n_1} \ge \epsilon$, il r\'{e}sulte de la d\'{e}finition \eqref{definition produit Cauchy} du produit de Cauchy des deux s\'{e}ries $G$ et $H$ qu'il existe des \'{e}l\'{e}ments $(k_0, k_1)$ et $(l_0, l_1)$ de $\mathbb{N}^2$ tels que $n_0 = k_0 + l_0$, $n_1 = l_1 + k_1$ et $\vert a_{k_0, k_1} \vert \beta^{k_0 + k_1} \vert b_{l_0, l_1} \vert \beta^{l_0 + l_1} \ge \epsilon.$  Par la d\'{e}finition des normes de $G$ et de $H$, on en d\'{e}duit qu'on doit alors avoir $\vert a_{k_0, k_1} \vert \beta^{k_0 + k_1} \ge \frac{\epsilon}{\max( \Vert G \Vert_\beta , \Vert H \Vert_\beta )}$ et $\vert b_{l_0, l_1} \vert \beta^{l_0 + l_1} \ge \frac{\epsilon}{\max( \Vert G \Vert_\beta , \Vert H \Vert_\beta )}$. En utilisant les implications \eqref{minorationa} et \eqref{minorationb}, on en d\'{e}duit que $(k_0, k_1) \in A(\epsilon)$ et $(l_0, l_1) \in B(\epsilon)$. Il suffit donc de prendre $C(\epsilon) = A(\epsilon) + B(\epsilon)$ pour trouver une partie finie de $\mathbb{N}^2$ satisfaisant la propri\'{e}t\'{e} d\'{e}sir\'{e}e. Par cons\'{e}quent $G H $ est \'{e}l\'{e}ment de $E_\beta$, et 
$$ \Vert G H \Vert_\beta = \sup_{(n_0, n_1) \in \mathbb{N}^2} \vert c_{n_0, n_1} \vert \beta^{n_0 + n_1} \le \Vert G \Vert_\beta \Vert H \Vert_\beta \; .$$

\end{dem} \findem

\begin{proposition}\label{operateur multiplication parH } Pour un nombre r\'{e}el $\beta > 0$, soit $H=\sum b_{n_0, n_1}x_0^{n_0}x_1^{n_1}$ une s\'{e}rie appartenant {\`a} l'espace  $   E_\beta$.  L'op\'{e}rateur multiplication $ {\rm mult}_H$ : $E_{\beta} \to  E_{\beta}$ qui envoie une s\'{e}rie $G$ sur $G H$ est continu.  \end{proposition}
 \begin{dem} D'apr{\`e}s la proposition \ref{Ebetastable}, l'application lin\'{e}aire ${\rm mult}_H$ est bien d\'{e}finie de $E_\beta$ dans lui-m{\^e}me, et lipschitzienne, donc continue. \end{dem}\findem
     
 \subsection{Op\'{e}rateur alpha} 
 On rappelle que $q = p^s$. 
 \begin{defi} 
 Soit $\beta > 1$ un nombre r\'{e}el, et $H(x_0, x_1)$ une s\'{e}rie appartenant {\`a} l'espace $E_{\beta^\frac{1}{q}}$. On appelle {\it op\'{e}rateur alpha} associ\'{e} {\`a} $H$ sur $E_\beta$ l'application compos\'{e}e 
 $$ \alpha = {\rm Dw}_q \circ {\rm mult}_H \circ i_{\beta, \beta^\frac{1}{q}} \; , $$
 o{\`u} ${\rm Dw}_q : E_{\beta^\frac{1}{q}} \to E_\beta$ est l'op\'{e}rateur de Dwork, ${\rm mult}_H : E_{\beta^\frac{1}{q}} \to E_{\beta^\frac{1}{q}}$ est l'op\'{e}rateur de multiplication par $H$, et $i_{\beta, \beta^\frac{1}{q}} : E_{\beta} \to E_{\beta^\frac{1}{q}}$ est l'injection canonique. 
 \end{defi} 
 
 \begin{proposition} Soit un r\'{e}el $\beta>1$, et $H$ un \'{e}l\'{e}ment de $E_{\beta^\frac{1}{q}}$; l'op\'{e}rateur alpha associ\'{e} {\`a} $H$ sur $E_\beta$ est un endomorphisme compl{\`e}tement continu de l'espace $E_\beta$.
 \end{proposition}
 \begin{dem}
 Pour $d \in \mathbb{N}$, consid\'{e}rons l'application $i_d : E_{\beta} \to E_{\beta^\frac{1}{q}}$ qui \'{a} toute s\'{e}rie $G$  associ\'{e}e  \'{a} la s\'{e}rie tronqu\'{e}e $i_d(G) = G_d =  \sum_{n_0+n_1\le d}
 a_{n_0, n_1} x_0^{n_0}x_1^{n_1} \;   $. L'image de $i_d$ \'{e}tant contenue dans le sous-espace engendr\'{e} par les mon\^{o}mes de degr\'{e} total $\le d$, sous-espace qui est de dimension finie, on voit que l'application lin\'{e}aire $i_d$ est de rang fini pour tout $d \in \mathbb{N}$. De plus, on a pour tout $G = \sum
 a_{n_0, n_1} x_0^{n_0}x_1^{n_1}$ \'{e}l\'{e}ment de $E_{\beta}$, 
 $$ \scriptstyle{\Vert G_d \Vert_{\beta^\frac{1}{q}}= \sup_{n_0 + n_1 \le d} \vert a_{n_0, n_1} \vert {\beta^{\frac{1}{q}(n_0 + n_1)}} \le \sup_{n_0 + n_1 \le d} \vert a_{n_0, n_1} \vert \beta \, ^{n_0 + n_1} \le \sup_{(n_0, n_1) \in \mathbb{N}^2} \vert a_{n_0, n_1} \vert \beta \, ^{n_0 + n_1} = \Vert G \Vert_{\beta} \; , }$$
 ce qui prouve que l'application lin\'{e}aire $i_d$ est continue de $E_{\beta}$ dans $E_{\beta^{\frac{1}{q}}}$. \\\\
  $$ \Vert G - G_d \Vert_{\beta^{\frac{1}{q}} } \le \Vert G \Vert_{\beta} \sup_{n_0 + n_1 > d} \beta^{-n_0 - n_1}{\beta^\frac{1}{q}} ^{(n_0 + n_1)} = \Vert G \Vert_{\beta} \left( \frac{{\beta^\frac{1}{q}} }{\beta} \right)^{d+1},  $$     
 ce qui prouve que $\Vert i_{\beta, {\beta^\frac{1}{q}} } - i_d \Vert \le  \left( \frac{{\beta^\frac{1}{q}} }{\beta} \right)^{d+1},$ et donc, puisque $\frac{{\beta^\frac{1}{q}} }{\beta} < 1$, que $\lim_{d \to + \infty} i_d = i_{\beta, {\beta^\frac{1}{q}} }$, c'est-\`{a}-dire que $i_{\beta, {\beta^\frac{1}{q}} 
 }$ est compl\`{e}tement continue.  

 Comme les applications ${\rm Dw}_q : E_{\beta^\frac{1}{q}} \to E_\beta$ , ${\rm mult}_H : E_{\beta^\frac{1}{q}} \to E_{\beta^\frac{1}{q}}$, sont continues d'apr\`{e}s les propositions  \ref{operateurDwork}, \ref{operateur multiplication parH }, alors  l'application compos\'{e}e 
 $ \alpha = {\rm Dw}_q \circ {\rm mult}_H \circ i_d \;$ est continue.
   \end{dem} \findem
 
 \begin{proposition}  \label{matricealpha } 
 La famille $(x_0^{n_0} x_1^{n_1})_{(n_0, n_1) \in \mathbb{N}^2} $ est une base de Schauder\cite{Bresis} orthogonale de l'espace $E_\beta$. La matrice dans cette base de l'op\'{e}rateur $\alpha$ associ\'{e} à la s\'{e}rie $H(x_0, x_1) = \sum_{(n_0, n_1) \in \mathbb{N}^2} b_{(m_0, m_1)} x_0^{m_0} x_1^{m_1} \in E_{\beta^\frac{1}{q}}$ est donn\'{e}e par 
 $$ \alpha(x_0^{n_0} x_1^{n_1} ) = \sum_{(m_0, m_1) \in \mathbb{N}^2, qm_0 \ge n_0, qm_1 \ge n_1} b_{q m_0 - n_0, q m_1 - n_1} x_0^{m_0} x_1^{m_1} \; . $$ 
 \end{proposition}
 \begin{dem}$$ \alpha(x_0^{n_0} x_1^{n_1} )= {\rm Dw}_q \circ {\rm mult}_H \circ i_{\beta,\beta^\frac{1}{q} }(x_0^{n_0} x_1^{n_1} )=$$ 
 $${\rm Dw}_q \circ {\rm mult}_H(x_0^{n_0} x_1^{n_1} )={\rm Dw}_q\left(\sum_{(m_0, m_1) \in \mathbb{N}^2} b_{(m_0, m_1)} x_0^{m_0} x_1^{m_1} (x_0^{n_0} x_1^{n_1} )\right)$$
 $$={\rm Dw}_q\left(\sum_{(m_0, m_1) \in \mathbb{N}^2} b_{m_0, m_1}  x_0^{n_0+m_0} x_1^{n_1+m_1}\right)$$
 $$={\rm Dw}_q\left(\sum_{(m_0, m_1) \in \mathbb{N}^2} b_{m_0, m_1}  x_0^{n_0+m_0} x_1^{n_1+m_1}\right)  $$
 $$={\rm Dw}_q\left(\sum_{(m_0, m_1) \in \mathbb{N}^2} b_{m_0-n_0, m_1-n_1}  x_0^{m_0} x_1^{m_1}\right) $$
 $$ =\sum_{(m_0, m_1)  \in \mathbb{N}^2, qm_0 \ge n_0, qm_1 \ge n_1}  b_{m_0q-n_0, m_1q-n_1}  x_0^{m_0} x_1^{m_1}$$ ce qui ach\`{e}ve la d\'{e}monstration.\end{dem} \findem

 \begin{proposition}\label{tracealpha } Soit un r\'{e}el $\beta>1$, et $H(x_0, x_1) = \sum_{(n_0, n_1) \in \mathbb{N}^2} b_{  bm_0, m_1} x_0^{m_0} x_1^{m_1}$ un \'{e}l\'{e}ment de $E_{\beta^\frac{1}{q}}$; l'op\'{e}rateur alpha associ\'{e} \'{a} $H$ est 
 un op\'{e}rateur nucl\'{e}aire de trace 
 $${\rm Tr}(\alpha) = \sum_{(n_0, n_1) \in \mathbb{N}^2} b_{(q-1)n_0, (q-1)n_1} \; .$$ \end{proposition} 
 \begin{dem} D'apr\`{e}s la  proposition \ref {matricealpha } la trace de la matrice associ\'{e}e {\`a} l'op\'{e}rateur alpha est la somme des \'{e}l\'{e}ments de la diagonale de cette matrice; autrement dit $${\rm Tr}(\alpha) = \sum_{(m_0, m_1) \in \mathbb{N}^2} b_{(q-1)n_0, (q-1)n_1} \; .$$
 
 \end{dem} \findem
 \subsection{Formule de traces}
 \begin{proposition} Soit $s \ge 1$ un entier naturel, $q = p^s$ et $t$ un \'{e}l\'{e}ment de $\Z_p\lbrack \mu_{q-1} \rbrack$ tel que $t^q = t$. Pour tout caract{\`e}re multiplicatif $\chi$ de $\textbf{W}_2(\F_q)$, il existe un r\'{e}el $\beta_0 > 1$ tel que,  pour tout $\beta  \in \rbrack 1, \beta_0 \lbrack$, on a l'\'{e}galit\'{e}
 $$ g(\psi_{2, s, t}, \chi) = (q-1)^2 {\rm Tr} ( \alpha) \; , 
 $$
o{\`u} $\alpha$ est l'op\'{e}rateur $\alpha : E_\beta \to E_\beta$ associ\'{e} {\`a}  la s\'{e}rie $\widehat{H}$ d\'{e}finie par \eqref{noyau}. 
 \end{proposition}
 \begin{dem}posons $$\widehat{H}(x_0, x_1) = \sum_{(m_0, m_1) \in \mathbb{N}^2} b_{ m_0, m_1} x_0^{m_0} x_1^{m_1}, $$ 
 D'apr\`{e}s la  proposition \ref {tracealpha } la trace de la matrice associ\'{e}e {\`a} l'op\'{e}rateur alpha associ\'{e} {\`a}  la s\'{e}rie $\widehat{H}$ est $${\rm Tr}(\alpha) = \sum_{(m_0, m_1) \in \mathbb{N}^2} b_{(q-1)m_0, (q-1)m_1} \; .$$
D'autre part, on a d'apr{\`e}s l'identit\'{e}  \ref{expression_analytique_somme_de _Gauss}  $$g(\psi_{2, s, t}, \chi) = \sum_{ z_0\in \F_q^\star}\sum_{ z_1 \in \F_q}
\widehat{H}(\Te(z_0), \Te(z_1)),$$ en posant$ \left(x_0,x_1\right)=\left((\Te(z_0), \Te(z_1)\right),$ d'o{\`u} $$g(\psi_{2, s, t}, \chi)=\sum_{ x_0\in \mu_{q-1}}\sum_{ x_1 \in \mu_{q-1}}
 \sum_{(n_0, n_1) \in \mathbb{N}^2} b_{  m_0, m_1} x_0^{m_0} x_1^{m_1},$$ et en utilisant l'identit\'{e} suivante:
 $$ \forall i=1,2...,n, \qquad \sum_{x_i \in \mu_{q-1}}x^{n_i}_i= \begin{cases} q-1  & \mbox{ si } q-1 /n_i \\ 
0 & \mbox{ sinon }  \end{cases}, $$on trouve alors que $$g(\psi_{2, s, t}, \chi)=(q-1)^2\sum_{(m_0, m_1) \in \mathbb{N}^2} b_{(q-1)n_0, (q-1)n_1} \; ,  $$ i.e $$g(\psi_{2, s, t}, \chi)=(q-1)^2{\rm Tr}(\alpha)$$

  \end{dem}\findem
 
 \begin{defi} 
 L'espace des s\'{e}ries surconvergentes est l'espace $\mathcal{H}^\dagger = \cup_{\beta > 1} E_\beta $ muni de la topologie de limite inductive des espaces vectoriels norm\'{e}s $E_\beta$. 
 \end{defi} 
 
 \begin{proposition} Sous les hypoth{\`e}ses pr\'{e}c\'{e}dentes,  la limite inductive $\alpha^\dagger$ des op\'{e}rateurs $\alpha$ associ\'{e}s sur $E_\beta$ {\`a} la s\'{e}rie $\widehat{H}$ d\'{e}finie par \eqref{noyau} est un op\'{e}rateur nucl\'{e}aire sur  $\mathcal{H}^\dagger$, et on a 
 $$ g(\psi_{2, s, t}, \chi) = (q-1)^2 {\rm Tr} ( \alpha^\dagger) \; .
 $$ \end{proposition}

\end{document}